\documentclass[10pt]{article}

\title{From Subfactors to Categories and Topology I. \\
       {\large Frobenius Algebras in and Morita Equivalence of Tensor Categories}\footnote{AMS subject classification: 18D10, 18D05; 46L37}}
\author{Michael M\"uger\thanks{Supported by EU through the TMR Networks ``Noncommutative Geometry''
and  ``Algebraic Lie Representations'', by MSRI through NSF grant DMS-9701755 and by NWO.} \\
Korteweg-de Vries Institute, Amsterdam, Netherlands \\ email: {\tt mmueger@science.uva.nl}}

\newlength{\dinwidth}
\newlength{\dinmargin}
\usepackage{latexsym, amsmath, amssymb, theorem}
\usepackage{tan-gle}

\input diagrams         

\def\1#1{{\bf #1}}
\def\2#1{{\cal #1}}
\def\5#1{{\sf #1}}
\def\6#1{{\mathfrak #1}}
\def\7#1{{\mathbb #1}}

\theoremheaderfont{\scshape}
\newtheorem{defin}{Definition}[section]
\newtheorem{lemma}[defin]{Lemma}
\newtheorem{prop}[defin]{Proposition}
\newtheorem{theorem}[defin]{Theorem}
\newtheorem{coro}[defin]{Corollary}
\newtheorem{conj}[defin]{Conjecture}
\newtheorem{defprop}[defin]{Definition-Proposition}
\theorembodyfont{\rmfamily}
\newtheorem{example}[defin]{Example}
\newtheorem{rema}[defin]{Remark}

\newcommand{\bdefin}{\begin{defin}}
\newcommand{\edefin}{\end{defin}}
\newcommand{\blemma}{\begin{lemma}}
\newcommand{\elemma}{\end{lemma}}
\newcommand{\bprop}{\begin{prop}}
\newcommand{\eprop}{\end{prop}}
\newcommand{\btheor}{\begin{theorem}}
\newcommand{\etheor}{\end{theorem}}
\newcommand{\bcoro}{\begin{coro}}
\newcommand{\ecoro}{\end{coro}}
\newcommand{\bconj}{\begin{conj}}
\newcommand{\econj}{\end{conj}}
\newcommand{\brem}{\begin{rema}}
\newcommand{\erem}{\hfill{$\Box$}\medskip\end{rema}}
\newcommand{\bexam}{\begin{example}}
\newcommand{\eexam}{\hfill{$\Box$}\medskip\end{example}}

\newcommand{\prf}{\noindent{\it Proof. }}
\newcommand{\qed}{\ \hfill $\blacksquare$\vspace{1pt}\\}

\def\be{\begin{equation}}
\def\ee{\end{equation}}
\newcommand{\ba}{\begin{array}}
\newcommand{\ea}{\end{array}}
\newcommand{\bea}{\begin{eqnarray}}
\newcommand{\eea}{\end{eqnarray}}
\newcommand{\bean}{\begin{eqnarray*}}
\newcommand{\eean}{\end{eqnarray*}}

\newcommand{\nn}{\nonumber}
\newcommand{\ve}{\varepsilon}

\newcommand{\impl}{\Rightarrow}
\newcommand{\rarr}{\rightarrow}
\newcommand{\restr}{\upharpoonright}
\newcommand{\ol}{\overline}
\newcommand{\ul}{\underline}
\newcommand{\mcirc}{\,\circ\,}

\newcommand{\oj}{{\overline{J}}}
\newcommand{\DS}{\displaystyle}

\newcommand{\obj}{\mbox{Obj}}
\newcommand{\Hom}{\mbox{Hom}}
\newcommand{\HOM}{\mbox{HOM}}
\newcommand{\End}{{\mbox{End}}}
\newcommand{\END}{\mbox{END}}
\newcommand{\id}{\mbox{id}}
\renewcommand{\mod}{\mbox{mod}}
\newcommand{\op}{{\mbox{\scriptsize op}}}
\newcommand{\cop}{{\mbox{\scriptsize cop}}}

\def\mobj#1{\raise .4\unitlens\hbox{\put(0,0){$#1$}}}
\def\mmobj#1{\raise .7\unitlens\hbox{\put(0,0){$#1$}}}

\newcommand{\mcm}[3]{\newcommand{#1}[#2]{{\ensuremath{#3}}}}

\newlength{\gwidth}     
\newlength{\gvert}      
\newlength{\gdrop}      
\newlength{\gbaredrop}  
\newlength{\goffset}    
\newlength{\gtemp}      
\newcommand{\present}[1]{%
\makebox[1\gwidth]{%
\rule[-1\gdrop]{0ex}{1\gvert}%
\raisebox{-1\gbaredrop}{#1}}}

\newcommand{\cinitdims}[2]{%
\setlength{\unitlength}{1em}%
\setlength{\goffset}{.35\unitlength}%
\setlength{\gwidth}{#1\unitlength}%
\setlength{\gvert}{#2\unitlength}%
\setlength{\gdrop}{.5\gvert}%
\addtolength{\gdrop}{-1\goffset}%
\setlength{\gbaredrop}{1\gdrop}%
\addtolength{\gvert}{.6\unitlength}%
\addtolength{\gdrop}{.3\unitlength}}    

\newcommand{\abovepic}[1]{%
\settoheight{\gtemp}{\ensuremath{#1}}%
\addtolength{\gvert}{1\gtemp}%
\settodepth{\gtemp}{\ensuremath{#1}}%
\addtolength{\gvert}{1\gtemp}}

\newcommand{\belowpic}[1]{%
\settoheight{\gtemp}{\ensuremath{#1}}%
\addtolength{\gvert}{1\gtemp}%
\addtolength{\gdrop}{1\gtemp}%
\settodepth{\gtemp}{\ensuremath{#1}}%
\addtolength{\gvert}{1\gtemp}%
\addtolength{\gdrop}{1\gtemp}}

\newcommand{\cell}[4]{\put(#1,#2){\makebox(0,0)[#3]{\ensuremath{#4}}}}

\newcommand{\prectwo}[3]%
{\begin{picture}(4.2,3.4)(-0.1,-0.2)%
\cell{2}{3.2}{b}{#1}%
\cell{2}{-0.2}{t}{#2}%
\cell{2.2}{1.5}{l}{#3}%
\qbezier(0,2)(2,4)(4,2)%
\qbezier(0,1)(2,-1)(4,1)%
\put(4,2){\vector(1,-1){0}}%
\put(4,1){\vector(1,1){0}}%
\put(2,2.5){\vector(0,-1){2}}%
\end{picture}}

\mcm{\ctwo}{3}{%
\cinitdims{4.2}{3.4}%
\abovepic{#1}%
\belowpic{#2}%
\present{\prectwo{#1}{#2}{#3}}}

\newcommand{\precthree}[5]{%
\begin{picture}(4.2,5.4)(-0.1,-0.2)%
\cell{2}{5.2}{b}{#1}%
\cell{1}{2.7}{b}{#2}%
\cell{2}{-.2}{t}{#3}%
\cell{2.2}{3.75}{l}{#4}%
\cell{2.2}{1.25}{l}{#5}%
\qbezier(0,3)(2,7)(4,3)%
\qbezier(0,2)(2,-2)(4,2)%
\put(0,2.5){\vector(1,0){4}}%
\put(2,4.5){\vector(0,-1){1.5}}%
\put(2,2){\vector(0,-1){1.5}}%
\put(4,3){\vector(1,-3){0}}%
\put(4,2){\vector(1,3){0}}%
\end{picture}}

\mcm{\cthree}{5}{%
\cinitdims{4.2}{5.4}%
\abovepic{#1}%
\belowpic{#3}%
\present{\precthree{#1}{#2}{#3}{#4}{#5}}}

\begin{document}
\maketitle\noindent

\numberwithin{equation}{section}

\abstract{We consider certain categorical structures that are implicit in subfactor
theory. Making the connection between subfactor theory (at finite index) and category
theory explicit sheds light on both subjects. Furthermore, it allows various
generalizations of these structures, e.g.\ to arbitrary ground fields, and the proof
of new results about topological invariants in three dimensions.

The central notion is that of a Frobenius algebra in a tensor category $\2A$,
which reduces to the classical notion if $\2A=\7F$-Vect, where $\7F$ is a
field. An object $X\in\2A$ with two-sided dual $\ol{X}$ gives rise to a
Frobenius algebra in $\2A$, and under weak additional conditions we prove
a converse: There exists a bicategory $\2E$ with $\obj\2E=\{\6A,\6B\}$ 
such that $\End_\2E(\6A)\stackrel{\otimes}{\simeq}\2A$ and such that there are 
$J, \oj: \ \6B\rightleftharpoons\6A$ producing the given Frobenius algebra. 
Many properties (additivity, sphericity, semisimplicity,\ldots) of $\2A$ carry
over to the bicategory $\2E$. 

We define weak monoidal Morita equivalence of tensor categories, denoted
$\2A\approx\2B$, and establish a correspondence between Frobenius algebras in
$\2A$ and tensor categories $\2B\approx\2A$. While considerably weaker than
equivalence of tensor categories, weak monoidal Morita equivalence
$\2A\approx\2B$ has remarkable consequences: $\2A$ and $\2B$ have equivalent (as
braided tensor categories) quantum doubles (`centers') and (if $\2A, \2B$ are semisimple
spherical or $*$-categories) have equal dimensions and give rise the same state sum
invariant of closed oriented 3-manifolds as recently defined by Barrett and Westbury. An 
instructive example is provided by finite dimensional semisimple and cosemisimple Hopf
algebras, for which we prove $H-\mod\approx\hat{H}-\mod$.

The present formalism permits a fairly complete analysis of the center of a semisimple
spherical category, which is the subject of the companion paper {\tt math.CT/0111205}.}


\section{Introduction}
Since tensor categories (or monoidal categories), in particular symmetric ones, have
traditionally been part and parcel of the representation theory of groups it is hardly
surprising that they continue to keep this central position in the representation theory 
of quantum groups, loop groups and of conformal field theories. See, e.g.\ \cite{cp,ka}. 
The main new ingredient in these applications is the replacement of the symmetry by a
braiding \cite{js1} which suggests connections with topology. Braided tensor categories
have in fact served as an input in new constructions of invariants of links and
3-manifolds and of topological quantum field theories \cite{t,ks}. (Recently it turned out
\cite{bw1,gk} that a braiding is not needed for the construction of the triangulation or
`state sum' invariant of 3-manifolds.)

A particular r\^{o}le in this context has been played by subfactor theory, see e.g.\
\cite{vfr1,ocn1,vfr-su,ek2}, which has led to the discovery of Jones' polynomial invariant
for knots \cite{vfr3}. Since the Jones polynomial was quickly reformulated in more
elementary terms, and due to the technical difficulty of subfactor theory, the latter
seems to have lost some of the attention of the wider public. This is deplorable, since
operator algebraists continue to generate ideas whose pertinence extends beyond subfactor
theory, e.g.\ in \cite{ocn3,ek3,yamag3,iz2}. The series of papers of which this is the
first aims at extracting the remarkable categorical structure which is inherent in
subfactor theory, generalizing it and putting it to use for the proof of new results in
categorical algebra and low dimensional topology. As will be evident to experts, the
series owes much to the important contributions of A.\ Ocneanu who, however, never
advocated a categorical point of view. 
We emphasize that our works will not assume any familiarity with subfactor theory
-- they are in fact also meant to convey the author's understanding of what (the more
algebraic side of) the theory of finite index subfactors is about. 

The present paper is devoted to the proof of several results relating two-sided adjoint
1-morphisms in 2-categories and Frobenius algebras in tensor categories. Here we are in
particular inspired by A.\ Ocneanu's notion of `paragroups' and R.\ Longo's
description of type III subfactors in terms of `Q-systems'. In order to make the series
accessible to readers with different backgrounds we motivate the constructions of the
present paper by considerations departing from classical Frobenius theory, from category
theory and from subfactor theory. While we will ultimately be interested in semisimple
spherical categories, a sizable part of our considerations holds in considerably greater
generality.


\subsection{Classical Frobenius Algebras}
One of the many equivalent criteria for a finite dimensional algebra $A$ over a field
$\7F$ to be a Frobenius algebra is the existence of a linear form $\phi: A\rarr\7F$ for
which the bilinear form $b(a,b)=\phi(ab)$ is non-degenerate. (For a nice exposition of the
present state of Frobenius theory we refer to \cite{kad}.) 
Recent results of Quinn and Abrams \cite{q,abrams0,abrams} provide the following
alternative characterization: A Frobenius algebra is a quintuple $(A,m,\eta,\Delta,\ve)$,
where $(A,m,\eta)$ and $(A,\Delta,\ve)$ are a finite dimensional algebra and coalgebra,
respectively, over $\7F$, subject to the condition
\[ m\otimes\id_A\mcirc \id_A\otimes\Delta = \Delta\mcirc m=
   \id_A\otimes m\mcirc\Delta\otimes\id_A. \]
In Section \ref{s-examples} we will say a bit more about the relation between these two
definitions. Since an $\7F$-algebra (coalgebra) is just a monoid (comonoid) in the tensor
category $\7F$-Vect, it is clear that the second definition of a Frobenius algebra makes
sense in any tensor category. A natural problem therefore is to obtain examples of
Frobenius algebras in categories other than $\7F$-Vect and to understand their
significance. The aim of the next two subsections will be to show how Frobenius algebras
arise in category theory and subfactor theory.


\subsection{Adjoint Functors and Adjoint Morphisms}
We assume the reader to be conversant with the basic definitions of categories, functors
and natural transformations, \cite{cwm} being our standard reference. (In the next section
we will recall some of the relevant definitions.) As is well known, the concept of adjoint
functors is one of the most important ones not only in category theory itself but also in
its applications to homological algebra and algebraic geometry. Before we turn to the
generalizations which we will need to consider we recall the definition. Given categories
$\2C, \2D$ and functors $F: \2C\rarr\2D,\ G: \2D\rarr\2C$, $F$ is a left adjoint of $G$,
equivalently, $G$ is a right adjoint of $F$, iff there are bijections 
\[ \phi_{X,Y}: \Hom_\2D(FX,Y)\cong\Hom_\2C(X,GY), \quad\quad X\in\obj\2C, Y\in\obj\2D, \]
that are natural w.r.t.\ $X$ and $Y$. This is denoted $F\dashv G$.
More convenient for the purpose of generalization is the equivalent
characterization, according to which $F$ is a left adjoint of $G$ iff there are natural
transformations $r: \id_\2C\rarr GF$ and $ s: FG\rarr\id_\2D$ satisfying 
\[ \id_G\otimes s \mcirc r\otimes \id_G = \id_G, \quad
    s\otimes\id_F\mcirc\id_F\otimes r = \id_F. \]

Given an adjunction $F\dashv G$, the composite functor $T=GF$ is an object in 
the (strict tensor) category $\2X=Fun(\2C)$ of endofunctors of $\2C$. (When seen as an
object of $\2X$, we denote the identity functor of $\2C$ by $\11$.) By the above
alternative characterization of adjoint pairs there are $\eta\equiv r\in\Hom_\2X(\11,T)$
and $m\in\Hom_\2X(T^2,T)$ defined by $m\equiv\id_G\otimes s\otimes\id_F$. These morphisms 
satisfy 
\bean m\mcirc m\otimes\id_T &=& m\mcirc\id_T\otimes m, \\
   m\mcirc \eta\otimes\id_T &=& m\mcirc\id_T\otimes \eta \ =\ \id_T; \eean
thus $(T, m, \eta)$ is a monoid in $\2X=Fun(\2C)$, equivalently, a monad in $\2C$.
Similarly, one finds that $(U, \Delta, \ve)\equiv(FG, \id_F\otimes r\otimes\id_G, s)$
is a comonoid in $Fun(\2D)$ or a comonad in $\2D$.

Now, given a monad in a category $\2C$ one may ask whether it arises from an adjunction as
above. This is always the case, there being two canonical solutions given by Kleisli and 
Eilenberg/Moore, respectively. See \cite[Chapter VI]{cwm} for the definitions and
proofs. In fact, considering an appropriate category of all adjunctions yielding the given 
monad, the above particular solutions are initial and final objects, respectively.

So far, we have been considering the particular 2-category ${\cal CAT}$ of small
categories. Thanks to the second definition of adjoints most of the above considerations 
generalize to an arbitrary 2-category $\2F$ (or even a bicategory). See 
\cite[Chapter 12]{cwm} and 
\cite{ks, gray} for introductions to 2- and bi-categories. Given objects (0-cells) $\2C,
\2D$ and 1-morphisms (1-cells) $F: \2C\rarr\2D, G: \2D\rarr\2C$ we say $G$ is a right
adjoint of $F$ iff there are 2-morphisms (2-cells) $r, s$ with the above properties. An
important special case pertains if the 2-category $\2F$ has only one object, say $\2C$. By
the usual `dimension shift' argument it then is the same as a (strict) tensor category,
and the adjoint 1-morphisms become dual objects in the usual sense. 

Now, a monad (comonad) in a general 2-category $\2F$ is most naturally defined \cite{str1}
as an object $\2C$ in $\2F$ (the basis) together with a monoid (resp.\ comonoid) in the
tensor category $\END_\2F(\2C)$. It is clear that an adjoint pair 
$F: \2C\rarr\2D, G: \2D\rarr\2C$ again gives rise to a monad (comonad) in $\2F$ with basis
objects $\2C$ ($\2D$). Again, the natural question arises whether every monad is
produced by an adjoint pair of 1-morphisms. Without restrictive assumptions on $\2F$ this
will in general not be true. (See \cite{str1} for a property of a 2-category which
guarantees that every monad arises from an adjoint pair of 1-morphisms.) The aim of the
present work is to explore a different aspect of this problem for which we need another
preparatory discussion. 

Assuming a 1-morphism in a 2-category $\2F$ (or an object in a tensor category) has both
left and right adjoints, there is in general no reason why they should be isomorphic
(i.e.\ related by an invertible 2-morphism). Yet there are tensor categories where every
object has a two-sided dual, in particular rigid symmetric categories (=closed
categories \cite{cwm}), rigid braided ribbon categories (=tortile categories \cite{js2}), 
$*$-categories \cite{lro}, and most generally, pivotal categories. (Functors, i.e.\
1-morphisms in $\2C\2A\2T$, which have a two-sided adjoint are occasionally called
`Frobenius functors'.) If a 1-morphism $F$ happens to have a simultaneous left and right
dual $G$, then not only $GF$ gives rise to a monad and $FG$ to a comonad, but also vice
versa. But in fact, there is more structure. Let $p: \id_\2D\rarr FG$ and 
$q: GF\rarr\id_\2C$ be the 2-morphisms associated with the adjunction $G\dashv F$ and
denote $T=GF$, $\ve=q$, $\Delta=\id_G\otimes p\otimes\id_F$. Then $(T, m, \eta, \Delta,
\ve)$ is a monoid and a comonoid in $\End(\2C)$, but in addition we have
\[ \id_T\otimes m\mcirc \Delta\otimes\id_T = \Delta\mcirc m = 
  m\otimes\id_T\mcirc \id_T\otimes\Delta. \]
The first half of this equation is proved diagrammatically by
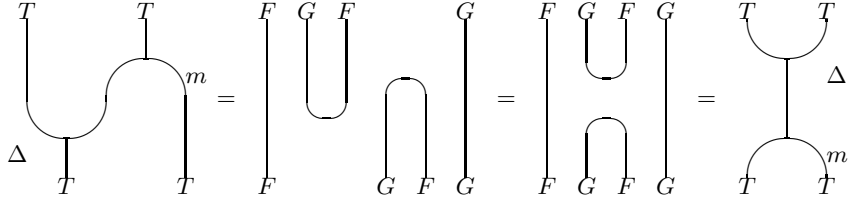
\[
\begin{tangle}
\object{T}\step[3]\object{T}\\
\id\Step\cd\mobj{m} \\
\step[-.5]\mobj{\Delta}\step[.5]\cu\Step\id\\
\step\object{T}\step[3]\object{T}
\end{tangle}
\quad = \quad
\begin{tangle}
\object{F}\step\object{G}\step\object{F}\step[3]\object{G} \\
\id\step\id\step\id\step\hcoev\step\id \\
\id\step\hev\step\id\step\id\step\id\\
\object{F}\step[3]\object{G}\step\object{F}\step\object{G}
\end{tangle}
\quad = \quad
\begin{tangle}
\object{F}\step\object{G}\step\object{F}\step\object{G} \\
\hh\id\step\id\step\id\step\id\\
\hh\id\step\ev\step\id\\
\hh\id\step\coev\step\id\\
\hh\id\step\id\step\id\step\id\\
\object{F}\step\object{G}\step\object{F}\step\object{G}
\end{tangle}
\quad = \quad
\begin{tangle}
\object{T}\Step\object{T}\\
\cu\mobj{\Delta}\\
\cd\mobj{m}\\
\object{T}\Step\object{T}
\end{tangle}
\]
and the other half similarly. Thus a 1-morphism $F:\6A\rarr\6B$ with a two-sided dual
gives rise to Frobenius algebras in the tensor categories $\End(\6A)$ and $\End(\6B)$.
Again, one may ask whether there is a converse, i.e.\ if every Frobenius algebra in a
tensor category arises as above. Already in the category $\2C=\7F$-Vect one finds
Frobenius algebras $\5Q=(Q,v,v',w,w')$ that are not of the form $Q=X\ol{X}$ for $X\in\2C$.
As one of our main results (Theorem \ref{main0}), we will see that under certain conditions
on a Frobenius algebra $\5Q$ in a tensor category $\2A$ there is a solution if one
embeds $\2A$ as a corner into a bicategory $\2E$. I.e., there are a bicategory $\2E$,
objects $\6A,\6B\in\obj\2E$ and a 1-morphism $J: \6B\rarr\6A$ with two-sided dual $\oj$
such that $\2A\simeq\END_\2E(\6A)$ and $\5Q$ arises via $Q=J\oj$ etc. If $\2A$ is
(pre)additive, abelian, $\7F$-linear, (semisimple) spherical or a $*$-category then $\2E$
will have the same properties. Under certain conditions, the bicategory $\2E$ is equivalent
to the bicategory $\2F$ containing the 1-morphism $J:\6B\rarr\6A$ which gave rise to the
Frobenius algebra.


\subsection{Subfactors}
For some basic definitions concerning subfactors we refer to Section \ref{ss-subfact}. For
the purposes of this introduction it is sufficient to know that a factor is a complex
unital $*$-algebra with center $\7C\11$, usually infinite dimensional. (Every finite
dimensional factor is isomorphic to a matrix algebra $M_n(\7C)$ for some $n\in\7N$.) 
A factor `has separable predual' if it admits a faithful continuous representation on a
separable Hilbert space. We restrict ourselves to such factors, which we (abusively) call 
separable. In his seminal work \cite{vfr1}, Vaughan Jones introduced a notion of index 
$[M:N]\in [1,\infty]$, defined for every inclusion $N\subset M$ of type II$_1$ factors. 
The index was soon generalized to factors of arbitrary type by Kosaki \cite{kos}. 
The index shares some basic properties with the one for groups: $[M:N]=1$ iff
$M=N$, $[P:M]\cdot [M:N]=[P:N]$ whenever $N\subset M\subset P$, etc. Yet the index is not 
necessarily an integer, in fact every $\lambda\in[4,\infty]$ can occur, whereas in the
interval $[1,4]$ only the countable set $\{ 4 \cos^2(\pi/n),\ n=3, 4, \ldots \}$ is
realized. The Jones index was soon generalized to factors of arbitrary types. Given
factors $P, Q$ (in fact, for general von Neumann algebras), there is a notion of $P-Q$
bimodule, and for bimodules ${}_P\2H_Q, {}_Q\2K_R$ there is a relative tensor product
${}_P\2H_Q\otimes{}_Q\2K_R$. This gives rise to a bicategory ${\cal BIM}$ whose objects
are factors, the 1-morphisms being bimodules and the 2-morphisms being the bimodule 
homomorphisms. Another bicategory ${\cal MOR}$ arising from factors has as objects all
factors, as 1-morphisms the  continuous unital $*$-algebra homomorphisms $\rho: P\rarr Q$
and as 2-morphisms the intertwiners. Thus if $\rho, \sigma: P\rarr Q$ then
\[ \Hom_{{\cal MOR}}(\rho,\sigma)= \{ x\in Q\ | \ x\rho(y)=\sigma(y)x \quad\forall y\in P \}. \]
Whereas the definition of the tensor product of bimodules is technically involved (in
particular in the non-type II$_1$ case, cf. \cite[Appendix V.B]{connes}), the composition
of 1-morphisms $\rho: P\rarr Q,\ \sigma: Q\rarr R$ in ${\cal MOR}$ is just their
composition $\sigma\circ\rho$ as maps and the unit 1-morphisms are the identity maps. Note
that the composition of 1-morphisms in ${\cal MOR}$ is strictly associative, thus ${\cal
MOR}$ is a 2-category. Every subfactor $N\subset M$ gives rise to a distinguished
1-morphism $\iota_{N,M}: N\rarr M$, the embedding map $N\hookrightarrow M$. Both ${\cal BIM}$ 
and ${\cal MOR}$ are $*$-bicategories, i.e.\ come with antilinear involutions on the
2-morphisms which reverse the direction. (For
$s\in\Hom_{\cal BIM}({}_P\2H_Q, {}_P\2K_Q)\subset \2B(\2H,\2K)$ and 
$t\in\Hom_{\cal MOR}(\rho,\sigma)\subset Q$, where $\rho, \sigma: P\rarr Q$, $s^*$ is 
given by the adjoint in $\2B(\2K,\2H)$ and in $Q$, respectively.)
If one restricts oneself to separable type III factors the corresponding
full sub-bicategories are equivalent in the sense of \cite{benab}: 
${\cal BIM}^{sep}_{III}\simeq{\cal MOR}^{sep}_{III}$, cf.\ \cite{connes,lo1}.
In the following discussion we will focus on the 2-category ${\cal MOR}^{sep}_{III}$.

The relevance of the concept of adjoints in 2-categories to subfactor theory is evident in
view of the following result which is due to Longo \cite{lo1}, though originally not
formulated in this way.

\btheor \label{th-lo1}
Let $N\subset M$ be an inclusion of separable type III factors. The embedding morphism
$\iota: N\hookrightarrow M$ has a two-sided adjoint $\ol{\iota}: M\hookrightarrow N$ iff
$[M:N]<\infty$. The dimensions of $\iota,\ol{\iota}$ (in the sense of \cite{lro}) are
related to the index by $d(\iota)=d(\ol{\iota})=[M:N]^{1/2}$. 
\etheor
Since separable type III factors are simple, every morphism $P\rarr Q$ is in fact an
embedding. We thus have the following more symmetric formulation.

\bcoro \label{c-lo2}
Let $\rho: P\rarr Q$ be a morphism of separable type III factors. Then $\rho$ has a
two-sided adjoint $\ol{\rho}: Q\rarr P$ iff $[Q:\rho(P)]<\infty$. Then
$d(\rho)=d(\ol{\rho})=[P:\ol{\rho}(Q)]^{1/2}=[Q:\rho(P)]^{1/2}$.
\ecoro
If these equivalent conditions are satisfied, an important object of study is the full
sub-bicategory of ${\cal MOR}^{sep}_{III}$ generated by $\rho$ and $\ol{\rho}$, which
consists of all morphisms between $P$ and $Q$ which are obtained by composition of $\rho,
\ol{\rho}$, the retracts and finite direct sums of such. (The type II$_1$ analogue is
Ocneanu's paragroup \cite{ocn1} associated with $N\subset M$.)

In this situation, $\gamma=\ol{\rho}\circ\rho: P\rarr P$ and 
$\sigma:\rho\circ\ol{\rho}: Q\rarr Q$ are endomorphisms of $P$ and $Q$, respectively, the
so-called canonical endomorphisms. By the considerations of the preceding section we know
that there is a canonical Frobenius algebra $(\gamma,v,v',w,w')$ in 
$\End(P):=\END_{{\cal MOR}}(P)$
and similarly for $\sigma$. Since we are in a $*$-categorical setting we have 
$v'=v^*, w'=w^*$. In \cite{lo3} such triples $(\gamma,v,w)$ were called Q-systems and it
was shown that for every Q-system where $\gamma\in\End(M)$ there exists a subfactor 
$N\subset M$ such that $(\gamma,v,w)$ arises as above. These results were clearly
motivated by the notion of conjugates (duals) in tensor categories, but the constructions
of Kleisli and Eilenberg/Moore do not seem to have played a r\^{o}le. 

Longo's results can be rephrased by saying that given a separable type III factor $M$
there is a bijection between subfactors $N\subset M$ with $[M:N]<\infty$ and Q-systems
(= canonical Frobenius algebras) in the tensor category $\End(M)$. Given such a
Frobenius algebra $\5Q$, we can on the one hand construct a subfactor $N$, an adjoint pair 
$\rho: M\rarr N,\, \ol{\rho}:N\rarr M$ and the 2-category generated by them. On the other
hand we can regard $\End(M)$ as an abstract $*$-category and apply to $\5Q$ the
construction announced in the preceding subsection. It will turn out that these two
procedures give rise to equivalent 2-categories.


\subsection{Organization of the paper}
In the following section we will discuss some preliminaries on (tensor) categories, mostly
concerning the notions of duality and dimension. The emphasis in this paper is on
spherical tensor categories, but we also consider $*$-categories. In Section 3 we
establish the connection between two-sided duals in bicategories and Frobenius algebras in
tensor categories. This is used in Section 4 to define the notion of weak monoidal Morita
equivalence. In Section 5 we consider categories that are linear over some field, in
particular spherical and $*$-categories. We show that the bicategory $\2E$ constructed
from a Frobenius algebra in a $*$-category or a (semisimple) spherical category is a
$*$-bicategory or (semisimple) spherical bicategory, and we prove the equality of certain
dimensions. Section 6 is devoted to several examples. We begin with classical Frobenius
algebras, i.e.\ Frobenius algebras in the category $\7F$-Vect, and specialize to finite
dimensional Hopf algebras. Our main result is the weak monoidal Morita equivalence
$H-\mod\approx\hat{H}-\mod$ whenever $H$ is semisimple and cosemisimple. We also discuss
in more detail the examples provided by subfactor theory. An application to the subject of
quantum invariants of 3-manifolds is given in Section 7. We outline a proof for the fact
that weakly monoidally Morita equivalent spherical categories define the same state sum
invariant for 3-manifolds. In the last section we briefly relate our results to previous
works and conclude with some announcements of further results and open problems.


\section{Categorical Preliminaries}\label{duals}
\subsection{Some basic notions and notations}
We assume known the standard definitions of tensor categories and bicategories, cf.\
\cite{cwm}. For definiteness all categories in this paper are supposed small. (In the
early stages essential smallness would suffice, whereas later on we will even
require the number of isoclasses to be finite.) We use `tensor category' and `monoidal
categories' interchangeably. Tensor categories will usually be assumed strict. (A tensor
category is strict if the tensor product satisfies associativity 
$X\otimes(Y\otimes Z)=(X\otimes Y)\otimes Z$ on the nose and the unit object $\11$
satisfies $X\otimes\11=\11\otimes X=X\ \forall X$. Similarly, in a strict bicategory (=
2-category) the composition of 1-morphisms is associative.) 
Since every tensor category is equivalent to a strict one \cite{cwm, js1} and every
bicategory to a 2-category this does not restrict the generality of our results.

Throughout, (2-)categories will be denoted by calligraphic letters
$\2A,\2B, \2C, \ldots$ and objects and morphisms in 1-categories by capital and lowercase
Latin letters, respectively. In 2-categories objects, 1- and 2-morphisms are denoted by 
Gothic ($\6A,\6B$,\ldots), capital and lowercase Latin letters, respectively. 
Unit objects and unit morphisms in tensor categories are denoted by $\11$ and 
$\id_X\in\Hom(X,X)$, respectively. Similarly, the unit 1- and 2-morphisms in 
2-categories are $\11_\6A\in\Hom(\6A,\6A)$ and $\id_X\in\Hom(X,X)$. If $\6A,\6B$ are
objects in a bicategory $\2E$ then $\Hom_\2E(\6A,\6B)$ and $\HOM_\2E(\6A,\6B)$ denote 
the corresponding 1-morphisms as a set and as a (1-)category (whose morphisms are the 
2-morphisms in $\2E$), respectively. As is well known,
$\END_\2E(\6A)\equiv\HOM_\2E(\6A,\6A)$ is a tensor category. (The composition $\circ$  
of $\6A-\6A$-morphisms in $\2E$ is the tensor product 
of objects in $\2A$ and the compositions $\circ,\otimes$ of 2-morphisms in $\2E$ are the
compositions $\circ,\otimes$ of morphisms in $\2A$.) Since the use of the composition
symbol $\circ$ for $\6A-\6A$-morphisms in $\2E$ as opposed to the tensor product $\otimes$
of objects in $\2A$ might lead to confusion (and there is only one composition for these
items anyway) we mostly omit the composition symbols $\otimes$ for objects in $\2A$ and
$\circ$ for 1-morphisms in $\2E$ altogether. 

We write $Y\prec X$ if $Y$ is a retract of $X$, i.e.\ there are morphisms 
$e: Y\rarr X, f: X\rarr Y$ such that $f\circ e=\id_Y$. We also say $Y$ is a subobject of
$X$, slightly incompatibly with common usage. In this situation $p=e\circ f\in\End(X)$ is
idempotent: $p^2=p$. Categories in which every idempotent arises in this way are called
`pseudo-abelian' (Karoubi), `Karoubian' (SGA) or `Cauchy-complete' \cite{bor}. Others say
that `idempotents split in  $\2A$' or `$\2A$ has subobjects' \cite{dr6}. We consider none
of these expressions particularly satisfactory. We will stick to the last alternative
since it goes best with `$\2A$ has direct sums'. Every category $\2A$ can be embedded as a
full subcategory into one which has subobjects. The latter can be defined as solution to a
universal problem, cf.\ \cite{bor}, which implies uniqueness up to equivalence. There
is, however, a well-known canonical solution, which we call $\ol{\2A}^p$, cf.\ e.g.\
\cite{gr}. Its objects are pairs $(X,p)$ where $X\in\obj\2A$ and $p=p^2\in\End(X)$. The
morphisms are given by 
\[ \Hom_{\ol{\2A}^p}((X,p),(Y,q))=\{ s: X\rarr Y\ | \ s=q\circ s\circ p\}. \]
If $\2A$ is a tensor category then so is $\ol{\2A}^p$ with 
$(X,p)\otimes (Y,q)=(X\otimes Y, p\otimes q)$.

A preadditive category (or Ab-category) is a category where all hom-sets are abelian
groups and the composition is additive w.r.t.\ both arguments. A preadditive category
$\2A$ is said to have direct sums (or biproducts) if for every $X_1, X_2\in\2A$ there are
$Y$ and morphisms $v_i\in\Hom(X_i,Y), v_i'\in\Hom(Y,X_i)$ such that 
$v_1\circ v_1'+v_2\circ v_2'=\id_Y$ and $v_i'\circ v_j=\delta_{i,j}\id_{X_i}$. We then
write $Y\cong X_1\oplus X_2$. Note that every $Y'\cong Y$ is a direct sum of $X_1,X_2$,
too. A preadditive category can always be embedded as a full subcategory into one with
(finite) direct sums. There is a canonical such category, cf.\ e.g.\ \cite{gr}, which we
call $\ol{\2A}^\oplus$. Again, this construction is compatible with a monoidal structure
on $\2A$. Constructions completely analogous to $\ol{\2A}^p, \ol{\2A}^\oplus$ exist for a 
bicategory $\2E$, cf.\ \cite[Appendix]{lro}. (Thus in $\ol{\2E}^p$ all idempotent
2-morphisms split and $\ol{\2E}^\oplus$ has direct sums for parallel 1-morphisms.) 
A pre-additive category is additive if it has direct sums and a zero object. If $\2A$ is
preadditive then $\ol{\2A}^p$ is additive: any $(X,0)$ is a zero object. Furthermore,
$\2A^{p\oplus}\simeq\2A^{\oplus p}$. 

Given a commutative ring $k$, a (monoidal) category is $k$-linear if all hom-sets are
finitely generated $k$-modules and the composition $\circ$ (and tensor product $\otimes$)
of morphisms are bilinear. Mostly, $k$ will be a field $\7F$. An object in a $k$-linear
category is called simple if $\End\,X=k\id_X$. (This property is often called absolute
simplicity or irreducibility. We drop the attribute `absolute', see the remarks below.) 
A $k$-linear category $\2C$ is called semisimple if it has direct sums and subobjects and
there are simple objects $X_i$ labeled by a set $I$ which are mutually disjoint 
($i\ne j\impl\Hom(X_i,X_j)=\{0\}$) such that the obvious map 
\[ \bigoplus_{i\in I} \Hom(Y,X_i)\otimes_\7F\Hom(X_i,Z)\longrightarrow \Hom(Y,Z) \]
is an isomorphism. Then every object $X$ is a finite direct sum of objects in 
$\{ X_i, i\in I\}$ and is determined up to isomorphism by the function
\[ I\rarr\7Z_+, \quad i\mapsto N^X_i =\dim\Hom(X,X_i)=\dim\Hom(X_i,X). \]
If $\2C$ is monoidal we also require $\11$ to be simple. 

In this paper we will not use the language of abelian categories, since we will not need
the notions of kernels and cokernels. (However, when applied to an abelian category the 
constructions given below give rise to abelian categories. See Remark \ref{r-alt}.) 
We briefly relate our definition of semisimplicity to more conventional terminology. An
abelian category is semisimple if it satisfies the following equivalent properties: (i)
all short exact sequences split, (ii) every monic is a retraction, (iii) every epi is a
section. An object $X$ in an abelian category is indecomposable if it is not a direct sum
of two non-zero objects, equivalently, if $\End\,X$ contains no idempotents besides $0$
and $\id_X$. If $\2C$ is semisimple and $X\in\2C$ is indecomposable then $\End\,X$ is a
skew field. Thus if $\7F$ is algebraically closed and $\2C$ is $\7F$-linear then
$X\in\2C$ is indecomposable iff it is (absolutely) simple. Since an abelian category has
direct sums and subobjects, $\2C$ is semisimple in our sense. Conversely, if $\2C$ is
semisimple in our sense and has a zero object then it is semisimple abelian, whether the
field $\7F$ is algebraically closed or not.)

A subcategory $\2S\subset\2C$ is full iff 
$\Hom_\2S(X,Y)=\Hom_\2C(X,Y) \ \ \forall X,Y\in\2S$, thus it is  determined by
$\obj\,\2S$. A subcategory is replete iff $X\in\obj\2S$ implies $Y\in\obj\2S$ for all
$Y\cong X$. (In the literature replete full subcategories are also
called strictly full.) Most subcategories we consider will be replete full. Isomorphism of
categories is denoted by $\cong$ and equivalence by $\simeq$. For preadditive (k-linear)
categories the functors establishing the equivalence/isomorphism are required to be
additive (k-linear) and for tensor categories they must be monoidal. In principle one
should use qualified symbols like 
$\cong_+, \simeq_k, \stackrel{\otimes}{\simeq},\stackrel{\otimes}{\cong}_k$ etc., where 
$+, k$ and $\otimes$ stand for additivity, k-linearity and monoidality, respectively, of
the equivalence/isomorphism. We will drop the $+, k$, hoping that they are obvious from
the context, but will write the $\otimes$.

The following definition is somewhat less standard. (Recall that we require the categories
to be small.)
\bdefin Two categories $\2A, \2B$ are Morita equivalent, denoted $\2A\approxeq\2B$, iff 
the categories 
$\mbox{Fun}(\2A^\op,\mbox{Sets})$, $\mbox{Fun}(\2B^\op,\mbox{Sets})$ are equivalent.
Preadditive categories are Morita equivalent iff the categories
$\mbox{Fun}_+(\2A^\op,\mbox{Ab})$, $\mbox{Fun}_+(\2B^\op,\mbox{Ab})$ of additive functors
are equivalent as preadditive categories.
$k$-linear categories are Morita equivalent iff the categories
$\mbox{Fun}_k(\2A^\op,\mbox{k-\mod})$, $\mbox{Fun}_k(\2B^\op,\mbox{k-\mod})$ of k-linear
functors are equivalent as k-linear categories.
\edefin

\bprop Let $\2A, \2B$ be categories. Then $\2A\approxeq\2B$ iff
$\ol{\2A}^p\simeq\ol{\2B}^p$. If $\2A, \2B$ are preadditive  (k-linear) categories then
$\2A\approxeq\2B$ iff $\ol{\2A}^{p\oplus}\simeq\ol{\2B}^{p\oplus}$ as preadditive  (k-linear)
categories.
\eprop
\prf For the first claim see \cite[Section 6.5]{bor}, for the others 
\cite[Chapter 2]{gr}. \qed 

This result motivates the following definition of (strong) Morita equivalence for
monoid\-al categories.

\bdefin Two tensor categories $\2A, \2B$ are strongly monoidally Morita equivalent 
($\2A\stackrel{\otimes}{\approxeq}\2B$) iff
$\ol{\2A}^p\stackrel{\otimes}{\simeq}\ol{\2B}^p$. If $\2A, \2B$ 
are preadditive  (k-linear) tensor categories then they are strongly Morita equivalent
iff $\ol{\2A}^{p\oplus}\stackrel{\otimes}{\simeq}\ol{\2B}^{p\oplus}$ as preadditive
(k-linear) tensor categories.
\edefin
While useful for certain purposes, this definition is unsatisfactory in that we cannot offer
an equivalent definition in terms of module-categories. In Section \ref{s-Morita} we will
define the notion of weak Morita equivalence $\2A\approx\2B$ of (preadditive, k-linear)
tensor categories. Is is genuine to tensor categories and satisfies
$\2A\stackrel{\otimes}{\approxeq}\2B\impl\2A\approx\2B$. We
speculate that $\2A\approx\2B$ iff suitably defined `representation categories' of $\2A,
\2B$ are equivalent, but we will leave this problem for future investigations.


\subsection{Duality in tensor categories and 2-categories}\label{ss-duality}
As explained in the Introduction, the notion of adjoint functors generalizes from
$\2C\2A\2T$ to an arbitrary 2-category $\2E$. Specializing to tensor categories, i.e.\
one-object 2-categories we obtain the following well-known notions. 
We recall that a tensor category $\2C$ is said to have left (right) duals if for every
$X\in\2C$ there is a $^*\!X$ ($X^*$) together with morphisms 
$e_X: \11\rarr X\otimes{}^*\!X$, $d_X: {}^*\!X\otimes X\rarr\11$ 
($\ve_X: \11\rarr X^*\otimes X$, $\eta_X: X\otimes X^*\rarr\11$) satisfying the usual
duality equations
\[ \id_X\otimes d_X\mcirc e_X\otimes\id_X = \id_X, \quad
  d_X\otimes\id_{\ol{X}}\mcirc\id_{\ol{X}}\otimes e_X=\id_{\ol{X}} \]
etc. Categories having left and right duals for every object
are called autonomous. Since duals, if they exist, are unique up to isomorphism the
conditions ${}^*\!X\cong X^*$ and $X^{**}\cong X$, which are easily seen to be equivalent,
do not involve any choices. If $X^*\cong {}^*\!X$ we speak of a two-sided dual of $X$. We
will exclusively consider categories with two-sided duals, for which we use the symmetric
notation $\ol{X}$. Assume now that $\2C$ is linear over the commutative ring
$k\equiv\End(\11)$. Then we have $\eta_X\circ e_X, d_X\circ\ve_X\in k$ and we would like
to consider them as dimensions of $X$ or $\ol{X}$. Yet, if $\lambda$ is a unit in $k$ then
replacing $e_X$ and $d_X$ by $\lambda e_X, \lambda^{-1}d_X$, respectively, the triangular
equations are still verified while $\eta_X\circ e_X$ and $d_X\circ\ve_X$ change. (The
product is invariant, though) We thus need a way to eliminate this indeterminacy. There
are three known solutions to this problem. If $\2C$ has a braiding $c_{-,-}$ and a twist
$\theta(-)$ we can determine the right duality in terms of the left duality by
$\ve_X=(\id_{\ol{X}}\otimes\theta_X)\circ c_{X,\ol{X}}\circ e_X$ and $\eta_X=d_X\circ
c_{X,\ol{X}}\circ(\theta_X\otimes\id_{\ol{X}})$, allowing to unambiguously define
$d(X)=\eta_X\circ e_X$. Since in this paper we do not require the existence of braidings
this approach is of no use to us. A fairly satisfactory way to eliminate the indeterminacy 
exists if $\2C$ has a $*$-operation, see Subsection \ref{stars}. The third solution, cf.\
the next subsection, is provided by the notion of spherical categories which was
introduced by Barrett and Westbury \cite{bw2,bw1}, elaborating on earlier work on pivotal
or sovereign categories. (For categories with braiding the latter approach is related to
the first one, in that there is a one-to-one correspondence between twists and spherical
structures, cf.\ \cite{y}.) We believe that this is the most general setting within which
results like Proposition \ref{eqofdims} and those of \cite{mue10} obtain without a
fundamental change of the methods. ($*$-categories can be turned into spherical
categories, cf.\ \cite{yamag2}, but this is not necessarily the most convenient thing to
do.) 

The following result shows that the {\it square} of the dimension of a simple object in a
linear category is well defined even in the absence of further structure.

\bprop Let $\2C$ be a $\7F$-linear tensor category with simple unit. If $X$ is simple and 
has a two-sided dual then $d^2(X)=(\eta_X\circ e_X)(d_X\circ\ve_X)\in\7F$ is a well
defined quantity. If $X, Y, XY$ are simple then $d^2(XY)=d^2(X)d^2(Y)$. Whenever $\2C$ has
a spherical or $*$- structure $d^2(X)$ coincides with $d(X)^2$ as defined using the latter.
\eprop
\prf Let $X$ be simple, $\ol{X}$ a two-sided dual and $e_X, d_X, \ve_X, \eta_X$ the
corresponding morphisms. By simplicity of $X$ we have
$\Hom(\11,X\ol{X})\cong\Hom(X,X)\cong\7F$, thus $\Hom(\11,X\ol{X})=\7F e_X$, 
$\Hom(\11,\ol{X}X)=\7F \ve_X$, etc. Therefore any other solution of the triangular
equations is given by 
\[ \tilde{e}_X=\alpha e_X,\ \  \tilde{d}_X=\alpha^{-1} d_X, \ \ 
   \tilde{\ve}_X=\beta\ve_X,\ \  \tilde{\eta}_X=\beta^{-1} \eta_X, \]
where $\alpha,\beta\in\7F^*$. Thus $(\eta_X\circ e_X)(d_X\circ\ve_X)$ does not depend on
the choice of the morphisms. The independence $d^2(X)$ of the choice of $\ol{X}$ and the
multiplicativity of $d^2$ (only if $XY$ is simple) are obvious. The final claim will follow
from the fact that in spherical and $*$- categories $d(X)$ is defined as $\eta_X\circ d_X$
for a certain choice of $d_X, \eta_X$ and the fact that $d(X)=d(\ol{X})$.
\qed

By this result the square of the dimension of a simple object is independent of the chosen 
spherical or $*$- structure and can in fact be defined without assuming the latter. Yet,
consistently choosing signs of the dimensions and extending $d(X)$ to an additive and
multiplicative function for all objects is a non-trivial cohomological problem to which
there does not seem to be a simple solution. A $*$-structure, when available, provides the
most natural way out, spherical structures being the second (but more general) choice.

The proposition implies that the following definition makes sense.
\bdefin Let $\2C$ be a semisimple $\7F$-linear tensor category with simple unit and two
sided duals. If $\2C$ has finitely many isomorphism classes of simple objects then we
define 
\[ \dim\2C=\sum_X d^2(X) \ \in\7F, \]
where the summation is over the isomorphism classes of simple objects and $d^2(X)$ is as 
in the proposition. If $\2C$ has infinitely many simple objects then we formally posit
$\dim\2C=\infty$.
\label{catdim}\edefin

We recommend the reader to skip the next two subsections until the structures introduced
there will be needed in Section \ref{lincat}.


\subsection{Spherical categories} \label{ss-spherical}
In spherical categories the problem mentioned above is solved by picking a two-sided dual 
$\ol{X}$ for every object $X$ and by specifying morphisms $\11\rarr X\otimes\ol{X}$, 
$X\otimes\ol{X}\rarr\11$ as part of the given data. This may look unnatural, but in
important cases such assignments are handy. (If, e.g., $H$ is an involutive Hopf algebra
then $\pi\mapsto \pi^t\circ S$ gives an involution on the representation category of $H$.)
We recall that we consider only strict tensor categories, and by the coherence theorem
of \cite{bw1} we may assume also strict duality. Since later we will have occasion to
construct spherical structures out of other data we give a redundancy-free definition of
spherical categories, which then will be proven equivalent to that of \cite{bw1,bw2}. 

\bdefin \label{pivotal}
A strict tensor category $\2C$ is a strict pivotal category if
there is a map $\obj\,\2C\rarr\obj\,\2C, X\mapsto\ol{X}$ such that
\begin{equation} \label{conjmap}
  \ol{\ol{X}}=X, \quad \ol{X\otimes Y}=\ol{Y} \otimes \ol{X},\quad \ol{\11}=\11 
\end{equation}
and there are morphisms $\ve(X): \11\rarr X\otimes\ol{X}$, 
$\ol{\ve}(X): X\otimes\ol{X}\rarr\11, X\in\2C$ satisfying the following conditions. 

\newarrow{Congruent} 33333
\begin{description}
\item[(1)] The composites
\begin{diagram}
 X\equiv X\otimes\11 & \rTo^{\id\otimes \ve(\ol{X})} & X\otimes\ol{X}\otimes X &
  \rTo^{\ol{\ve}(X)\otimes\id} & \11\otimes X\equiv X,
\end{diagram}
\begin{diagram}
 X\equiv \11\otimes X & \rTo^{\ve(X)\otimes\id} & X\otimes\ol{X}\otimes X &
  \rTo^{\id\otimes\ol{\ve}(\ol{X})} & X\otimes\11\equiv X
\end{diagram}
coincide with $\id_X$.
\item[(2)] The diagrams
\begin{diagram}
\11 & \rTo^{\ve(X)} & X\otimes\ol{X} \\
\dTo^{\ve(X\otimes Y)} & & \dTo_{\id\otimes\ve(Y)\otimes\id} \\
X\otimes Y\otimes \ol{X\otimes Y} & \rCongruent & X\otimes Y\otimes\ol{Y}\otimes \ol{X} 
\end{diagram}
and
\begin{diagram}
\11 & \lTo^{\ol{\ve}(X)} & X\otimes\ol{X} \\
\uTo^{\ol{\ve}(X\otimes Y)} && \uTo_{\id\otimes\ol{\ve}(Y)\otimes\id} \\
X\otimes Y \otimes\ol{X\otimes Y} & \rCongruent & X\otimes Y\otimes\ol{Y}\otimes \ol{X}
\end{diagram}
commute for all $X,Y$.
\item[(3)] For every $s: X\rarr Y$ the composites
\begin{diagram}
\ol{Y}\equiv\ol{Y}\otimes\11 & \rTo^{\id\otimes \ve(X)} & \ol{Y}\otimes X\otimes \ol{X} &
  \rTo^{\id\otimes s\otimes\id} & \ol{Y}\otimes Y\otimes\ol{X} 
  & \rTo^{\ol{\ve}(\ol{Y})\otimes\id} & \11\otimes\ol{X}\equiv\ol{X},
\end{diagram}
\begin{diagram}
\ol{Y}\equiv \11\otimes \ol{Y} & \rTo^{\ve(\ol{X})\otimes\id} & \ol{X}\otimes X\otimes\ol{Y} &
 \rTo^{\id\otimes s\otimes\id} & \ol{X}\otimes Y\otimes\ol{Y} &
  \rTo^{\id\otimes\ol{\ve}(Y)} & \ol{X}\otimes\11\equiv \ol{X}
\end{diagram}
coincide in $\Hom(\ol{Y},\ol{X})$.
\end{description}
For every $X$ and $s: X\rarr X$ we define morphisms in $\End(\11)$ by
\begin{diagram}
tr_L(s): &
\11 & \rTo^{\ve(X)} & X\otimes\ol{X} & \rTo^{s\otimes\id} & X\otimes\ol{X} & \rTo^{\ol{\ve}(X)} & \11,
\end{diagram}
\begin{diagram}
tr_R(s): &
\11 & \rTo^{\ve(\ol{X})} & \ol{X}\otimes X & \rTo^{\id\otimes s} & \ol{X}\otimes X & \rTo^{\ol{\ve}(\ol{X})} & \11.
\end{diagram}
$\2C$ is spherical iff $tr_L(s)=tr_R(s)$ for all $s$.
\edefin

\brem 1. If for $s\in\Hom(X,Y)$ we now define $\ol{s}\in\Hom(\ol{Y},\ol{X})$ by the formulas in
(3)\ we easily verify $\ol{\ol{s}}=s$ and $\ol{s\circ t}=\ol{t}\circ\ol{s}$. Thus the maps
$X\mapsto\ol{X}, s\mapsto\ol{s}$ constitute an involutive contravariant endofunctor of
$\2C$. Condition (1) can now be expressed as $\ol{\id_X}=\id_{\ol{X}}\ \forall X$. 

2. Using the conditions (1) one verifies that  $\ol{\ve}(X)=\ol{\ve(X)}$ and thus consistency
of our notation. 

3. The definition of the map $s\mapsto \ol{s}$ implies that the square 
\be
\begin{diagram}
\11 & \rTo^{\ve(X)} & X\otimes \ol{X} \\
\dTo^{\ve(Y)} && \dTo_{s\otimes\id} \\
Y\otimes\ol{Y} & \rTo^{\id\otimes \ol{s}} & Y\otimes\ol{X}
\end{diagram}
\label{dualm}\end{equation}
commutes for every $s: X\rarr Y$. 
Thus our axioms for strict pivotal categories imply those of \cite{bw1,bw2}, and the
converse is obvious if we put $\ol{\ve}(X):=\ol{\ve(X)}$. 

4. Even in those applications where the pivotal category $\2C$ under study is strict
monoidal, like in type III subfactor theory, the duals are rarely strict. In general one
has natural isomorphisms 
$\tau_X: X\rarr\ol{\ol{X}}$, $\gamma_{X,Y}: \ol{X}\otimes\ol{Y}\rarr \ol{Y\otimes X}$ and 
$\nu: \11\rarr\ol{\11}$ which must satisfy a number of compatibility conditions. See
\cite[Theorem 1.9]{bw2}, where it is proven that such a category is monoidally equivalent to
a strict pivotal category as defined above. To be sure, in practice one does not really
want to strictify the categories under consideration in order to work with them. Just as
the well known coherence results on (braided) tensor categories \cite{cwm,js1}, 
\cite[Theorem 1.9]{bw2} can be rephrased as follows: All computations in a strict pivotal or
spherical category remain valid in the non-strict case, provided one inserts the morphisms
$\tau_X, \gamma_{X,Y}, \nu$ wherever needed. If this is possible in different ways no
possible result depends on the choices one makes in this process. In order to travel with
slightly lighter luggage we may therefore calmly stick to the strict case. 

5. A pivotal category is called non-degenerate if for all $X,Y$ the pairing 
\[ \Hom(X,Y)\times\Hom(Y,X)\rarr \7F,\ \  s\times t\mapsto \langle s,t\rangle= 
    tr_X(s\circ t)=tr_Y(t\circ s) \]
between $\Hom(X,Y)$ and $\Hom(Y,X)$ is non-degenerate. In \cite[Lemma 3.1]{gk} it is
proved that semisimple pivotal categories are non-degenerate.
\erem

If $\2C$ is $\7F$-linear pivotal with simple unit $\11$ then we define dimensions by
\[ d(X)\,\id_\11=tr_L(X)=\ol{\ve}(X)\circ\ve(X). \]
If $\2C$ is spherical then 
\[ d(X)=d(\ol{X}) \quad \forall X. \]
We will show that in the semisimple case this is equivalent to sphericity. Let 
$X\cong\bigoplus_{j\in J}N_iX_i$, thus there are 
$u_i^\alpha: X_i\rarr X, u_i': X\rarr X_i^\alpha$ such that
\begin{equation} \sum_i\sum_{\alpha=1}^{N_i} u_i^\alpha\circ {u'}_i^\alpha=\id_X, \quad 
   {u'}_i^\alpha\circ u_j^\beta=\delta_{i,j}\delta_{\alpha,\beta}\id_{X_i}. 
\label{ctz}\end{equation}
Using the conjugation functor $\ol{{\ }}$ we define $v_i^\alpha:=\ol{{u'}_i^\alpha}:
X\rarr X_i$, ${v'}_i^\alpha:=\ol{u_i^\alpha}: X_i\rarr X$, and clearly
\[ \sum_i\sum_{\alpha=1}^{N_i} v_i^\alpha\circ {v'}_i^\alpha=\id_X, \quad 
   {v'}_i^\alpha\circ v_j^\beta=\delta_{i,j}\delta_{\alpha,\beta}\id_{X_i}. \]

\blemma Let $\2C$ be a semisimple pivotal tensor category with simple unit. Then $\2C$ is
spherical iff $d(X)=d(\ol{X})$ for all $X$.  \label{spher}\elemma
\prf We compute
\bean tr_L(s) &=& \ol{\ve}(X)\mcirc s\otimes \id_{\ol{X}} \mcirc\ve(X) \\ 
  &=& \sum_{i,\alpha}\sum_{j,\beta} \ol{\ve}(X)\mcirc 
  (u_i^\alpha\circ {u'}_i^\alpha\circ s\circ u_j^\beta\circ{u'}_j^\beta)\otimes\id_{\ol{X}} \mcirc\ve(X) \\
  &=& \sum_{i,\alpha}\sum_{j,\beta} \ol{\ve}(X)\mcirc 
  (u_i^\alpha\circ {u'}_i^\alpha\circ s\circ u_j^\beta)\otimes v_j^\beta \mcirc\ve(X_j) \\
  &=& \sum_{i,\alpha}\sum_{j,\beta} \ol{\ve}(X_j)\mcirc 
  ({u'}_j^\beta\circ u_i^\alpha\circ {u'}_i^\alpha\circ s\circ u_j^\beta)\otimes\id_{\ol{X_j}} \mcirc\ve(X_j) \\
  &=& \sum_{i,\alpha}\ol{\ve}(X_i)\mcirc ({u'}_i^\alpha\circ s\circ u_i^\alpha)\otimes\id_{\ol{X_j}} \mcirc\ve(X_i) \\
  &=& \sum_{i,\alpha} s_i^\alpha \,\ol{\ve}(X_i)\circ\ve(X_i) = \sum_{i,\alpha} s_i^\alpha d(X_i),
\eean
where $s_i^\alpha \id_{X_i}={u'}_i^\alpha\circ s\circ u_i^\alpha\in\Hom(X_i,X_i)$. We have
used (\ref{dualm}) and the fact that the $X_i$ are simple. In a similar way one computes
\[ tr_R(s)=\sum_{i,\alpha} s_i^\alpha\, \ol{\ve}(\ol{X_i})\circ\ve(\ol{X_i}) =
   \sum_{i,\alpha} s_i^\alpha d(\ol{X_i})
\]
and the result follows from the assumption. \qed

We will need the following facts concerning the behavior of sphericity under certain
categorical constructions: 
\blemma \label{herit}
Let $\2A, \2B$ be pivotal (spherical). Then $\2A^\op, \2A\otimes_\7F\2B$ (the product in
the sense of enriched category theory) and the completions $\ol{\2A}^\oplus, \ol{\2A}^p$
are pivotal (spherical).
\elemma
\prf We only sketch the proof, restricting us to the strict pivotal/spherical case. For
the opposite category the claimed facts are obvious, choose $\ve(X^\op)=\ol{\ve}(X)$ etc. 
As to $\ol{\2A}^p$, recall that its objects are given by pairs
$(X,p), X\in\obj\2A, p=p^2\in\End(X)$ and those of $\ol{\2A}^\oplus$ by finite sequences 
$(X_1,\ldots,X_l)$ of objects in $\2A$. We define the duality maps on $\2A^\op$,
$\2A\otimes_\7F\2B$, $\ol{\2A}^p$, $\ol{\2A}^\oplus$ by
\[ \ol{X^\op}=\ol{X}^\op, \quad \ol{X\boxtimes Y}=\ol{X}\boxtimes\ol{Y}, \quad \ol{(X,p)}=
  (\ol{X},\ol{p}), \quad \ol{(X_1,\ldots, X_l)}=(\ol{X_1},\ldots,\ol{X_l}), \]
respectively. The conditions (\ref{conjmap}) are clearly satisfied. We define further
\[ \ve(X^\op)=\ol{\ve}(X)^\op, \quad \ve(X\boxtimes Y)=\ve(X)\boxtimes\ve(Y), \quad\ve((X,p))
  = p\otimes\ol{p}\circ\ve(X)
\]
and
\[ \quad \ve((X_1,\ldots,X_l))= \sum_{i=1}^l    u_i\otimes\ol{u_i} \circ\ve(X_i)  , \]
where the $u_i: X_i\rarr (X_1,\ldots, X_l)$ are the injections which together with the
$\{u_i'\}$ satisfy (\ref{ctz}). The easy verification of the axioms is left to the
reader. It is clear that the new categories are spherical if $\2A, \2B$ are
spherical. \qed

The definition of strict pivotal tensor categories can be generalized to 2-categories
$\2E$ in a straightforward manner. Every 1-morphism $X:\6A\rarr\6B$ has a two-sided dual 
$\ol{X}:\6B\rarr\6A$. This map has properties which are obvious generalizations of 
the monoidal situation. In particular, $\ol{\11_\6X}=\11_\6X$ for all objects $\6X$.
For $X\in\Hom(\6A,\6B)$ there are 
$\ve(X): \11_\6B\rarr X\ol{X}$ and $\ol{\ve}(X): X\ol{X}\rarr\11_\6B$ satisfying 
conditions which are analogous to the those in pivotal categories. Again the conjugation
can be extended to a 2-functor $\ol{{\ }}: \2E^\op\rarr\2E$, where $\2E^\op$ has the same
objects as $\2E$ and 1- and 2-morphisms are reversed. (This functor acts trivially on the
objects.) All this can be obtained from \cite{mack} by ignoring the monoidal structure on
$\2E$ considered there. There is one difference, however, which requires attention. For
$X: \6A\rarr\6B$ and $s\in\End(X)$ the two traces 
\[ tr_L(s)=\ol{\ve}(X)\mcirc s\otimes\id_{\ol{X}}\mcirc\ve(X), \quad
   tr_L(s)=\ol{\ve}(\ol{X})\mcirc \id_X\otimes s\mcirc\ve(\ol{X}) \]
take values in different commutative monoids, viz.\ $\End(\11_\6B)$ and $\End(\11_\6A)$,
respectively. Thus defining spherical 2-categories as satisfying $tr_L(s)=tr_R(s)$ for all 
2-morphisms $s\in\End(X)$ -- not only those where $X$ has equal source and range --
makes sense only if there is a way to identify the $\End(\11_\6X)$ for different $\6X$. We
will restrict ourselves to 2-categories where $\11_\6X$ is absolutely simple for all
objects $\6X$, thus all $\End(\11_\6X)$ are canonically isomorphic to $\7F$ via 
$c\mapsto c\,\id_{\11_\6X}$. (The bicategories ${\cal BIM}$ and ${\cal MOR}$ of bimodules 
and $*$-homomorphisms of factors discussed in the introductions provide examples.) Under
this condition every 1-morphism has a $\7F$-valued dimension with the expected properties.


\subsection{Duality in $*$-Categories} \label{stars}
In this section we posit $\7F=\7C$. Many complex linear categories have an additional
piece of structure: a positive $*$-operation. See \cite{dr6, w} for two important classes
of examples. A $*$-operation on a $\7C$-linear category is map which assigns to every
morphism $s: X\rarr Y$ a morphism $s^*: Y\rarr X$. This map has to be antilinear,
involutive ($s^{**}=s$), contravariant ($(s\circ t)^*=t^*\circ s^*$)  
and monoidal ($(s\otimes t)^*=s^*\otimes t^*$) if the category is monoidal. A
$*$-operation is called positive if $s^*\circ s=0$ implies $s=0$. A {\it tensor
$*$-category} is a $\7C$-linear tensor category with a positive $*$-operation. (For a
braided tensor $*$-category one often requires unitarity of the braiding, but there are
examples where this is not satisfied by a naturally given braiding.) Some relevant
references, where infinite dimensional Hom-sets are admitted, are \cite{glr, dr6} and
\cite{lro}, the latter reference containing a very useful discussion of 2-*-categories. In 
\cite[Proposition 2.1]{mue6} we showed that $\7C$-linear categories with positive
$*$-operation and finite dimensional Hom-sets are $C^*$- and $W^*$-categories in the sense
of \cite{glr,dr6}. In $W^*$-categories one has a polar decomposition theorem for morphisms
\cite[Cor.\ 2.7]{glr}, which implies, e.g., that if $\Hom(Y,X)$ contains an invertible
morphism then it contains a unitary morphism ($u\circ u^*=u^*\circ u=\id$). 

The notion of duality in $*$-categories as considered in \cite{dr6,lro} has two peculiar
features. First, it is automatically two-sided. Secondly, there is no compelling reason to
fix a duality {\it map} $X\mapsto\ol{X}$ and to choose morphisms $\11\rarr X\ol{X}$ etc.
Rather it is sufficient to assume that all objects {\it have} a conjugate. We stick to the
term `conjugate' from \cite{dr6,lro} in order to underline the conceptual difference.

An object $\ol{X}$ is said to be a conjugate of $X$ if there are 
$r_X\in\Hom(\11,\ol{X}X), \ol{r}_X\in\Hom(\11,X\ol{X})$ satisfying the  
{\it conjugate equations}:
\begin{equation} \ol{r}_X^*\otimes\id_X\mcirc\id_X\otimes r_X=\id_X,
  \quad\quad  r_X^*\otimes\id_{\ol{X}}\mcirc\id_{\ol{X}}\otimes
  \ol{r}_X= \id_{\ol{X}}. \label{conj-eq0}\end{equation}
A category $\2C$ has conjugates if every object $X\in\2C$ has a conjugate 
$\ol{X}\in\2C$. The triple $(\ol{X},r_X,\ol{r}_X)$ is called a solution of the conjugate
equations. It is called normalized if
$r_X^*\circ r_X=\ol{r}_X^*\circ\ol{r}_X\in\Hom(\11,\11)$. (Since $\7C$ is algebraically
closed every solution can be normalized.) A solution of the conjugate equations is called
standard if $r_X=\sum_i \ol{w}_i\otimes w_i\mcirc r_i$ where $w_i, \ol{w}_i$ are
isometries effecting decompositions $X=\oplus_i X_i$, $\ol{X}=\oplus_i\ol{X_i}$ into
simple objects, and $(\ol{X_i},r_i, \ol{r}_i)$ are normalized solutions of
(\ref{conj-eq0}) for $X_i$. For any object $X$ we define a dimension by 
$d(X)=r_X^*\circ r_X$, where $(\ol{X},r_X,\ol{r}_X)$ is a normalized standard solution. 
Then $d(X)$ is well defined and satisfies $d(X)=d(\ol{X})$, $d(X\oplus Y)=d(X)+d(Y)$ and
$d(X\otimes Y)=d(X)d(Y)$ for all $X,Y\in\2C$. The dimension takes values in the set $\{
2\cos\frac{\pi}{n},\ n=3,4,\ldots\}\cup[2,\infty)$. If $(\ol{X},r,\ol{r})$ is any
normalized solution of the conjugate equations for $X$ then $d(X)\le r^*\circ r$ and
equality holds iff $(\ol{X},r,\ol{r})$ is standard. For the proofs we refer to \cite{lro}.

Furthermore, a braided $*$-category with conjugates automatically \cite{mue6} has a
canonical twist \cite{js1}. Together with the fact \cite{y} that in braided tensor
categories there is a one-to-one correspondence between sovereign structures and twists
this implies that every braided $*$-category has a canonical sovereign structure. It is
not unreasonable to believe that this is true even in the absence of a braiding. In fact,
in \cite{yamag1} Yamagami considered spherical structures (`$\ve$-structures') compatible
with a given $*$-structure, and in \cite{yamag2} he shows that every $*$-category can be
equipped with an essentially unique spherical structure. See also \cite{fgsv} for similar
considerations. In any case, in the $*$-case one does not need a spherical structure since
every scalar quantity (i.e.\ morphism $\11\rarr\11$) is unambiguously defined if one
sticks to the above normalization rules for the duality morphisms $r_X, \ol{r}_X$, cf.\
\cite{lro}.


\subsection{Graphical notation}
In this paper we will often represent computations with morphisms in a tensor category in
terms of tangle diagrams rather than by formulas or commutative diagrams. Since this
notation is well known, cf.\ e.g.\ \cite{ka}, we just explain our conventions.

Our diagrams are to be read upwards, thus with 
$a: X\rarr Y, b: Y\rarr Z, c: U\rarr V, d: V\rarr W$ we represent
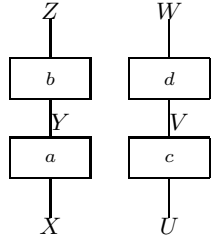
\[ b \otimes d \circ a\otimes c = (b\circ a)\otimes (d\circ c)\in\Hom(X\otimes U, Z\otimes W) 
\]
by
\[ 
\begin{tangle}
\hstep\object{Z}\step[3]\object{W}\\
\hh\hstep\id\step[3]\id\\
\frabox{b}\Step\frabox{d}\\
\hh\hstep\id\obj{Y}\step[3]\id\obj{V}\\
\frabox{a}\Step\frabox{c}\\
\hh\hstep\id\step[3]\id\\
\hstep\object{X}\step[3]\object{U}
\end{tangle}
\]
Representing $\ve(X)$ by
$\ \begin{tangle}\object{X}\step\object{\ol{X}}\\ \hh\hev \end{tangle}$
and $\ol{\ve}(X)$ by
$\ \begin{tangle} \hh\hcoev\\ \object{X}\step\object{\ol{X}}\end{tangle}$
we depict condition (1) of Definition \ref{pivotal} by
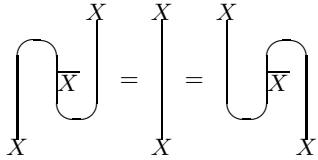
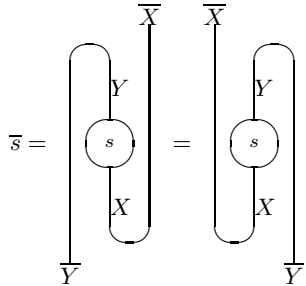
\[ 
\begin{tangle}
\Step\object{X}\\
\hh\hcoev\step\id\\
\hh\id\step\id\obj{\ol{X}}\step\id\\
\hh\id\step\hev\\
\object{X}
\end{tangle}
\ \ = \ \ 
\begin{tangle}
\object{X}\\
\id\\
\hh\id\\
\object{X}
\end{tangle}
\ \ = \ \
\begin{tangle}
\object{X}\\
\hh\id\step\hcoev\\
\hh\id\step\id\obj{\ol{X}}\step\id\\
\hh\hev\step\id\\
\Step\object{X}
\end{tangle}
\]
and (3) becomes
\[ 
\ol{s}= \ \ 
\begin{tangle}
\Step\object{\ol{X}}\\
\hh\hcoev\step\id\\
\hh\id\step\id\obj{Y}\step\id\\
\id\step\O s\step\id\\
\hh\id\step\id\obj{X}\step\id\\
\hh\id\step\hev\\
\object{\ol{Y}}
\end{tangle}
\ \ = \ \ 
\begin{tangle}
\object{\ol{X}}\\
\hh\id\step\hcoev\\
\hh\id\step\id\obj{Y}\step\id\\
\id\step\O s\step\id\\
\hh\id\step\id\obj{X}\step\id\\
\hh\hev\step\id\\
\Step\object{\ol{Y}}
\end{tangle}
\]


\section{Frobenius Algebras vs.\ Two-sided Duals}\label{catap}
\subsection{From two sided duals to Frobenius algebras}
Let $\2A$ be a tensor category which for notational simplicity we assume to be strict. As
is well known this does not restrict the generality of our considerations. Ultimately we
will be interested in linear categories over a field $\7F$, but a large part of our
considerations is independent of this and will be developed without assuming linearity.

\bdefin \label{d-Frob}
Let $\2A$ be a (strict) tensor category. A Frobenius algebra in $\2A$ is a
quintuple $(Q,v,v',w,w')$, where $Q$ is an object in $\2A$ and 
$v:\11\rightarrow Q, v': Q\rightarrow\11, w: Q\rightarrow Q^2, w':Q^2\rightarrow Q$ 
are morphisms satisfying the following conditions:
\begin{equation} \label{W1}
   w\otimes\id_Q\mcirc w=\id_Q\otimes w\mcirc w, 
\end{equation}
\begin{equation} \label{W1'}
   w' \mcirc w'\otimes\id_Q=w' \mcirc \id_Q\otimes w', 
\end{equation}
\begin{equation} \label{W2}
   v'\otimes\id_Q\mcirc w = \id_Q = \id_Q\otimes v'\mcirc w.
\end{equation}
\begin{equation} \label{W2'}
   w' \mcirc v\otimes\id_Q = \id_Q = w' \mcirc \id_Q\otimes v. 
\end{equation}
\begin{equation} \label{W3}
   w'\otimes\id_Q \mcirc \id_Q\otimes w = w\mcirc w' = \id_Q\otimes w'\mcirc w\otimes \id_Q. 
\end{equation} 
\edefin

\begin{figure}
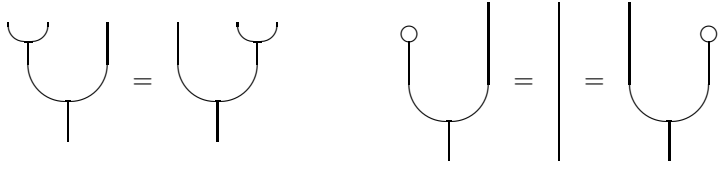

\bean
\begin{tangle}
\hh\hcu\step[1.5]\id \\
\hstep\cu \\
\end{tangle}
\ \ = \ \ 
\begin{tangle}
\hh\id\step[1.5]\hcu \\
\cu \\
\end{tangle}
\quad\quad & \quad\quad 
\begin{tangle}
\unit\Step\id \\
\cu \\
\end{tangle}
\ \ = \ \ 
\begin{tangle}
\id \\
\id \\
\end{tangle}
\ \ = \ \ 
\begin{tangle}
\id\Step\unit \\
\cu \\
\end{tangle}
\eean
\caption{Comonoids in a category}
\label{fig0}\end{figure}

\begin{figure}
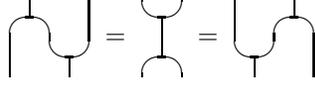

\[ 
\begin{tangle}
\hh\cd\step\id \\ 
\hh\id\step\cu
\end{tangle}
\ =\ 
\begin{tangle}
\hh\cu \\ 
\hh\cd
\end{tangle}
\ =\ 
\begin{tangle}
\hh\id\step\cd \\ 
\hh\cu\step\id
\end{tangle}
\]
\caption{The Frobenius condition}
\label{fig00}\end{figure}

\brem \label{frob-rem}
1. Throughout the paper we use the following symbols:
\[ \begin{tangle}\hh\hcu\end{tangle}=w, \quad \begin{tangle}\hh\hcd\end{tangle}=w', \quad
   \begin{tangle}\hh\counit\end{tangle}\ =v, \quad \begin{tangle}\hh\unit\end{tangle}\ =v'. \]
For the tangle diagrams corresponding to the above conditions see Figures\ \ref{fig0},
\ref{fig00}. 

2. Eqs.\ (\ref{W1}-\ref{W2'}) amount to requiring that $(Q,w',v)$ and $(Q,w,v')$ are
a monoid and a comonoid, respectively, in the category $\2C$. The  new ingredient is
the Frobenius condition (\ref{W3}), cf.\ Fig.\ \ref{fig00}, which can be interpreted as 
expressing that $w:Q\rarr Q^2$ is a map of $Q$-$Q$ bimodules. This must not be confused
with the more familiar bialgebra condition. The latter makes sense only if $\2C$ comes
with a braiding, which we do not assume. For this reason we avoid the usual symbols 
$m, \Delta, \ve, \eta$. 

3. To the best of our knowledge (\ref{W3}) makes its first appearance in 
\cite[Appendix A3]{q}, where it is part of an alternative characterization of symmetric
algebras (in Vect). This is a special case of the alternative characterization of
Frobenius algebras mentioned in the Introduction and discussed in more detail in
Subsection \ref{s-examples}. 

4. If $\2A$ is a $*$-category we will later on require $w'=w^*, v'=v^*$, which obviously
renders (\ref{W1'}), (\ref{W2'}) redundant.
\erem

\bdefin \label{d-isom}
Frobenius algebras $(Q,v,v',w,w')$,
$(\tilde{Q},\tilde{v},\tilde{v}',\tilde{w},\tilde{w}')$ in the (strict) tensor category
$\2A$ are isomorphic if there is an isomorphism $s':Q\rarr\tilde{Q}$ such that
\[ s\circ v=\tilde{v},\ \  v'=\tilde{v}'\circ s, \ \ 
   s\otimes s\circ w=\tilde{w}\circ s, \ \ s\circ w'=\tilde{w}'\circ s\otimes s. \]
\edefin

The above definitions are vindicated by the following.
\blemma \label{motiv}
Let $J: \6B\rarr\6A$ be a 1-morphism in a 2-category $\2E$ and let $\oj:\6A\rarr\6B$ be a
two-sided dual with duality 2-morphisms $d_J, e_J, \ve_J, \eta_J$. Positing
$Q=J\oj: \2A\rarr\2A$ there are $v,v',w,w'$ such that $(Q,v,v',w,w')$ is a Frobenius
algebra in the tensor category $\2A=\HOM_\2E(\6A,\6A)$. 
\elemma
\prf Since $J, \oj $ are mutually two-sided duals there are 
\[ e_J: \11_\6A\rarr J\oj , \ \ \eta_J: J\oj \rarr\11_\6A, \ \ 
   \ve_J: \11_\6B\rarr \oj J, \ \ d_J: \oj J\rarr\11_\6B \]
satisfying
\be \ba{ccccc} \displaystyle \id_J\otimes d_J\mcirc e_J\otimes\id_J & 
  = \id_J = & \eta_J\otimes\id_J\mcirc\id_J\otimes\eta_J, \\
 \displaystyle \id_\oj\otimes\eta_J\mcirc\ve_J\otimes\id_\oj & = 
  \id_\oj = & d_J\otimes\id_\oj\mcirc\id_\oj\otimes e_J.
\ea\label{cJ}\end{equation}
Defining  
\be \ba{ccccc} \displaystyle v &=& e_J &\in& \Hom_\2E(\11_\6A, J\oj ), \\
  \displaystyle v' &=& \eta_J &\in& \Hom_\2E(J\oj , \11_\6A), \\
  \displaystyle w &=& \id_J\otimes \ve_J\otimes\id_\oj &\in&  \Hom_\2E(J\oj , J\oj J\oj ), \\
  \displaystyle w' &=& \id_J\otimes d_J\otimes\id_{\oj } &\in& \Hom_\2E(J\oj J\oj , J\oj), 
\ea\end{equation}
we have $v\in\Hom_\2A(\11,Q),\ w\in\Hom_\2A(Q,Q^2)$, $v'\in\Hom_\2A(Q, \11)$ and
$w'\in\Hom_\2A(Q^2, Q)$. Now (\ref{W1}) follows simply from the interchange law
\[
\begin{tangle}
\object{J}\step\object{\oj }\step\object{J}\step\object{\oj }\step\object{J}\step\object{\oj }\\
\nw2\hev\step\id\step\id\step\id\\
\Step\d\hev\dd\\
\step[3]\object{J}\step\object{\oj }
\end{tangle}
\ \ \ \ = \ \ \ \ 
\begin{tangle}
\object{J}\step\object{\oj }\step\object{J}\step\object{\oj }\step\object{J}\step\object{\oj }\\
\id\step\id\step\id\step\hev\ne2\\
\d\hev\dd\\
\step\object{J}\step\object{\oj }
\end{tangle}
\]
and (\ref{W2}) from the duality equations (\ref{cJ}):
\[
\begin{tangle}
\hh\hcoev\step\id\step\id\\
\d\hev\dd\\
\step\object{J}\step\object{\oj }
\end{tangle}
\ \ \ \ = \ \ \ \ 
\begin{tangle}
\id\step\id\\
\id\step\id\\
\object{J}\step\object{\oj }
\end{tangle}
\ \ \ \ = \ \ \ \ 
\begin{tangle}
\hh\id\step\id\step\hcoev\\
\d\hev\dd\\
\step\object{J}\step\object{\oj }
\end{tangle}
\]
The conditions (\ref{W1'}), (\ref{W2'}) are analogous, and the Frobenius algebra condition
(\ref{W3}) follows from 
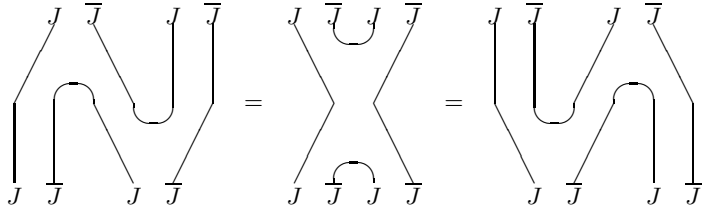
\[ 
\begin{tangle}
\step\object{J}\step\object{\oj }\Step\object{J}\step\object{\oj }\\
\dd\hcoev\d\step\id\step\id\\
\id\step\id\step\d\hev\dd\\
\object{J}\step\object{\oj }\Step\object{J}\step\object{\oj }
\end{tangle}
\quad = \quad
\begin{tangle}
\object{J}\step\object{\oj }\step\object{J}\step\object{\oj }\\
\d\hev\dd\\
\dd\hcoev\d\\
\object{J}\step\object{\oj }\step\object{J}\step\object{\oj }
\end{tangle}
\quad = \quad
\begin{tangle}
\object{J}\step\object{\oj }\Step\object{J}\step\object{\oj }\\
\id\step\id\step\dd\hcoev\d\\
\d\hev\dd\step\id\step\id\\
\step\object{J}\step\object{\oj }\Step\object{J}\step\object{\oj }
\end{tangle}
\]
\qed 

\brem 1. By duality we obviously also obtain a Frobenius algebra 
$\hat{\5Q}=(\hat{Q}=\oj J, \ldots)$ in $\2B=\HOM_\2E(\6B,\6B)$. 

2. Considering a tensor category $\2C$ as a 2-category with a single object $\6A$ we
obtain the special case of an object $X$in $\2C$ with a two-sided dual $\ol{X}$. 

3. Since duals are unique up to isomorphism a different choice of $\oj$ changes $Q$ only
within its isomorphism class. Yet it is in general not true that the Frobenius algebra 
$(Q,v,v',w,w')$ is well defined up to isomorphism.
\erem

\blemma \label{lemma-selfd}
The object $Q$ of a Frobenius algebra is two-sided self-dual, i.e.\
there are $r\in\Hom(\11,Q^2)$, $r'\in\Hom(Q^2,\11)$ such that the duality equations are
satisfied with $e_Q=\ve_Q=r$, $d_Q=\eta_Q=r'$. 
\elemma
\prf
\prf Set $r=w\circ v, r'=v'\circ w'$. Then the duality equations follow by 
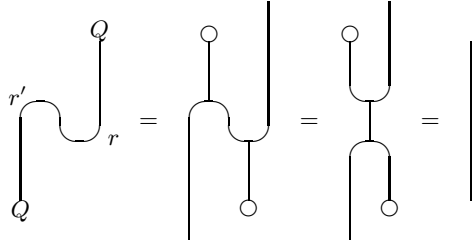
\[
\begin{tangle}
\Step\object{Q}\\
\hh\Step\id\\
\hh\step[-.3]\mobj{r'}\step[.3]\coev\step\id\\
\hh\id\step\ev\step[.2]\mobj{r} \\
\hh\id\\
\object{Q}
\end{tangle}
\quad = \quad
\begin{tangle}
\hstep\unit\step[1.5]\id\\
\hh\cd\step\id\\
\hh\id\step\cu\\
\id\step[1.5]\counit
\end{tangle}
\quad = \quad
\begin{tangle}
\unit\step\id\\
\hh\cu\\
\hh\cd\\
\id\step\counit
\end{tangle}
\quad =\quad
\begin{tangle}
\id\\
\id
\end{tangle}
\]
Here the first equality holds by definition, the second follows from (\ref{W3})
and the last from (\ref{W2}). 
\qed


\subsection{A universal construction}
Now we prove a converse of Lemma \ref{motiv} by embedding a given tensor category $\2A$
with a Frobenius algebra $(Q,v,v',w,w')$ into a suitably constructed 2-category such that
$Q=J\oj $.

\bdefin An almost-2-category is defined as a 2-category \cite{ks} except that we do not
require the existence of a unit 1-morphism $\11_{\6X}$ for every object $\6X$. \edefin

\bprop \label{constr}
Let $\2A$ be a strict tensor category and $\5Q=(Q,v,v',w,w')$ a Frobenius algebra in
$\2A$. Then there is an almost-2-category $\2E_0$ satisfying
\begin{enumerate}
\item $\displaystyle \obj\,\2E_0=\{\6A,\6B\}$.
\item There is an isomorphism $I: \2A\rarr\displaystyle\HOM_{\2E_0}(\6A,\6A)$ of tensor
categories. 
\item There are 1-morphisms $J: \6B\rarr \6A$ and $\oj : \6A\rarr \6B$ such that
$J\oj=I(Q)$.  
\end{enumerate}
If $\2A$ is $\7F$-linear then so is $\2E_0$. Isomorphic Frobenius algebras give rise to
isomorphic almost 2-categories.
\eprop
\prf The proof is constructive, the definition of the {\it objects} obviously being forced
upon us: $\ \displaystyle \obj\,\2E_0=\{\6A, \6B\}$. \\  

\noindent{\it 1-morphisms.} We define formally
\bea \Hom_{\2E_0}(\6A,\6A) &=& \obj\,\2A, \nn\\
  \Hom_{\2E_0}(\6B,\6A) &=& \{ ``X J",\ X\in\obj\,\2A \}, \nn\\
  \Hom_{\2E_0}(\6A,\6B) &=& \{``\oj  X",\ X\in\obj\,\2A \}, \label{1morphisms}\\
  \Hom_{\2E_0}(\6B,\6B) &=& \{ ``\oj  X J",\ X\in\obj\,\2A \}. \nn
\eea
For the moment this simply means that $\Hom_{\2E_0}(\6A,\6B), \Hom_{\2E_0}(\6B,\6A)$ and 
$\Hom_{\2E_0}(\6B,\6B)$ are isomorphic to $\Hom_\2A(\6A,\6A)$ as sets. In 
particular, with $X=\11\in\obj\2A$ we obtain the following distinguished 1-morphisms:
\bean J &=& ``\11 J" \ \in\Hom_{\2E_0}(\6B,\6A), \\
  \oj  &=& ``\oj  \11" \ \in\Hom_{\2E_0}(\6A,\6B). 
\eean
\\
{\it Composition of 1-morphisms}, wherever legal, is defined by juxtaposition, followed
by replacing a possibly occurring composite $J\oj $ by the underlying object $Q$ of the
Frobenius algebra. (The latter is the case whenever one considers $X\circ Y$ where
$Y\in\Hom_{\2E_0}(?,\6B), X\in\Hom_{\2E_0}(\6B,?)$.) 
We refrain from tabulating all the possibilities and give instead a few examples:
\[\ba{lcc} ``\oj  X" \mcirc ``Y J" &:=& ``\oj (XY) J"  \in\Hom_{\2E_0}(\6B,\6B), \\
   ``X J"\mcirc ``\oj  Y J" &:=& ``(XQY) J" \in\Hom_{\2E_0}(\6B,\6A).
\ea\]
$\2A$ being strict by assumption, the composition of 1-morphisms is obviously strictly
associative. With this definition the set of 1-morphisms is the free semigroupoid (=
small category) generated by $\obj\2A\cup\{J, \oj \}$ modulo $J\oj =Q$ and the
relations in $\2A$.
If we had to consider only 1-morphisms we could drop the quotes in (\ref{1morphisms})
since now $J, \oj $ are legal 1-morphisms and, e.g., $``XJ"$ is just the composition 
$X\circ J$ of $X$ and $J$. But in order to define the 2-morphisms and their compositions
and to verify that we obtain a 2-category we must continue to distinguish between 
$``X" ``J"$ and $``XJ"$ for a while. \\  

\noindent{\it 2-morphisms as sets.}
Since we want $\HOM(\6A,\6A)\cong\2A$, we clearly have to set
\[ \Hom_{\2E_0}(X,Y)=\Hom_\2A(X,Y), \quad X,Y\in\Hom_{\2E_0}(\6A,\6A)=
  \obj\,\2A. \]
In order to identify the remaining 2-morphisms we appeal to duality which
in the end should hold in $\2E$. Applied to $\Hom_{\2E_0}(\6B,\6A)$ this means
\[ \Hom_{\2E_0}(``X\circ J", ``Y\circ J") \cong 
  \Hom_{\2E_0}(``X\circ J\circ\oj ", Y) = 
  \Hom_{\2E_0}(XQ,Y)= \Hom_\2A(XQ,Y). \]
This means that the elements of $\Hom_{\2E_0}(``X J", ``Y J")$
\[
\begin{tangle}
\object{Y}\step\object{J}\\
\hh\id\step\id\\
\hh\frabox{s}\\
\hh\id\step\id\\
\object{X}\step\object{J}
\end{tangle}
\]
where $X,Y\in\obj\2A$, are represented by those of $\Hom_\2A(XQ,Y)$, etc:
\[
\begin{tangle}
\object{Y}\\
\hh\id\step\hcoev\mobj{e_J}\\
\hh\frabox{s}\step\id\\
\hh\id\step\id\step\id\\
\object{X}\step\underbrace{\object{J}\step\object{\oj }}_{Q}\\
\step[1.5]
\end{tangle}
\]
Thus we \ul{define}
\bean \Hom_{\2E_0}(``X J", ``Y J") &=& 
   \Hom_\2A(XQ,Y), \\
\Hom_{\2E_0}(``\oj  X", ``\oj  Y") &=& 
   \Hom_\2A(X,Q Y),\\
\Hom_{\2E_0}(``\oj  X J", ``\oj  Y J") &=& 
  \Hom_\2A(XQ,Q Y), \\
\eean
as sets. 
Now we must define the vertical and horizontal compositions of the 2-morphisms in 
$\2E_0$, which we denote $\bullet$ and $\times$, respectively, in order to avoid
confusion with the compositions in $\2A$. \\ \\
\noindent{\it Vertical ($\bullet$-) Composition of 2-morphisms.}
\begin{figure}
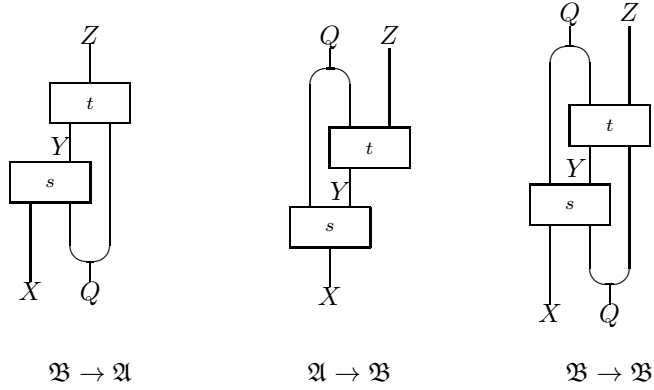

\[\ba{ccccc}
\begin{tangle}
\step[1.5]\object{Z}\\
\hh\step[1.5]\id\\
\step\frabox{t}\\
\hh\hstep\obj{Y}\hstep\id\step\id\step\\
\hh\frabox{s}\step\id\\
\hh\id\step\id\step\id\\
\hh\id\step\hcu \\
\object{X}\step[1.5]\object{Q}
\end{tangle}
&\quad\quad\quad\quad&
\begin{tangle}
\hstep\object{Q}\step[1.5]\object{Z}\\
\hh\hcd\step\id\\
\hh\id\step\id\step\id\\
\hh\id\step\frabox{t}\\
\hh\id\hstep\obj{Y}\hstep\id\\
\hh\frabox{s}\\
\hh\hstep\id\\
\hstep\object{X}
\end{tangle}
&\quad\quad\quad\quad&
\begin{tangle}
\hstep\object{Q}\step[1.5]\object{Z}\\
\hh\hcd\step\id\\
\hh\id\step\id\step\id\\
\hh\id\step\frabox{t}\\
\hh\id\step[.4]\obj{Y}\step[.6]\id\step\id\step\\
\hh\frabox{s}\step\id\\
\hh\id\step\id\step\id\\
\hh\id\step\hcu \\
\object{X}\step[1.5]\object{Q}
\end{tangle} \\ \\
\6B\rightarrow\6A && \6A\rightarrow\6B && \6B\rightarrow\6B
\ea\]
\caption{Vertical ($\bullet$-)composition of 2-morphisms in $\2E_0$}
\label{fig1}\end{figure}
Let $X,Y,Z\in\obj\2A$. For $s\in\Hom_{\2E_0}(X,Y)\equiv\Hom_{\2A}(X,Y)$, 
$t\in\Hom_{\2E_0}(Y,Z)\equiv\Hom_{\2A}(Y,Z)$ it is clear that $t\bullet s=t\circ s$. 
\[
\6A%
\cthree{X}{Y}{Z}{s}{t}%
\6A
\]

Let now $s\in\Hom_{\2E_0}(``X\circ J", ``Y\circ J")\equiv\Hom_\2A(XQ,Y)$,  
$t\in\Hom_{\2E_0}(``Y\circ J", ``Z\circ J")\equiv\Hom_\2A(YQ,Z)$. Then we define
$t\bullet s\in \Hom_{\2E_0}(``X\circ J", ``Z\circ J")\equiv\Hom_\2A(XQ,Z)$ by
\[ t\bullet s= t\mcirc s\otimes\id_Q\mcirc\id_X\otimes w. \]

Similarly, for $s\in\Hom_{\2E_0}(``\oj \circ X", ``\oj \circ Y")\equiv\Hom_\2A(X,QY)$, 
$t\in\Hom_{\2E_0}(``\oj \circ Y", ``\oj \circ Z")\equiv\Hom_\2A(Y,Q Z)$ we define
$t\bullet s \in\Hom_{\2E_0}(``\oj \circ X",``\oj \circ Z")\equiv\Hom_\2A(X,QZ)$ by
\[ t\bullet s= w'\otimes\id_Z\mcirc\id_Q\otimes t\mcirc s. \]

Finally, for 
$s\in\Hom_{\2E_0}(``\oj \circ X\circ J", ``\oj \circ Y\circ J")\equiv\Hom_\2A(XQ,Q Y)$, 
$t\in\Hom_{\2E_0}(``\oj \circ Y\circ J", ``\oj \circ Z\circ J")\equiv\Hom_\2A(YQ,Q Z)$
we define 
$t\bullet s\in\Hom_{\2E_0}(``\oj \circ X\circ J",``\oj \circ Z\circ J")\equiv\Hom_\2A(XQ,Q Z)$ by
\begin{equation} t\bullet s= w'\otimes\id_Z\mcirc\id_Q\otimes t\mcirc 
  s\otimes\id_Q\mcirc\id_X\otimes w. \label{bb}\end{equation}
See Fig.\ \ref{fig1} for diagrams corresponding to these definitions. Associativity of the
$\bullet$-composit\-ions is easily verified using the coassociativity (\ref{W1}) and
associativity (\ref{W1'}) of the Frobenius algebra $(Q,v,v',w,w')$. \\

\begin{figure}
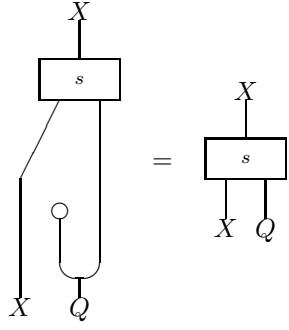

\bean
\begin{tangle}
\step[1.5]\object{X}\\
\hh\step[1.5]\id\\
\step\frabox{s}\\
\dd\step\id \\
\id\step\unit\step\id\\
\hh\id\step\hcu \\
\object{X}\step[1.5]\object{Q}
\end{tangle}
&\quad=\quad&
\begin{tangle}
\hstep\object{X}\\
\hh\hstep\id\\
\frabox{s}\\
\hh\id\step\id\\
\object{X}\step\object{Q}
\end{tangle}
\eean
\caption{Unit 2-morphism $\id_{``XJ"}$}
\label{fig4}\end{figure}

\noindent{\it Unit 2-arrows.} It is clear $\id_X\in\Hom_\2A(X,X)$ for $X\in\obj\2A$ 
is the unit 2-arrow for the 1-morphism $X: \6A\rightarrow \6A$. Furthermore, using
the above rules for the $\bullet$-composition of 2-morphisms and the equations
(\ref{W2}), (\ref{W2'}), it is easily verified that
\[ \id_{``XJ"}=\id_X\otimes v'\in \Hom_\2A(XQ,X)\equiv\Hom_{\2E_0}(``XJ",``XJ") \]
is in fact the unit 2-arrow $\id_{``XJ"}$. Diagrammatically, the equation
$s\bullet\id_{``XJ"}=s$ with 
$s\in\Hom_\2A(XQ,X)\equiv\Hom_{\2E_0}(``XJ",``XJ")$ looks as in
Fig.\ \ref{fig4}. 
Similarly, we have 
\bean \id_{``\oj X"} &=& v\otimes\id_X\in \Hom_\2A(X,QX) \equiv
   \Hom_{\2E_0}(``\oj X",``\oj X"), \\ 
 \id_{``\oj XJ"} &=& v\otimes\id_X\otimes v'
  \in \Hom_\2A(XQ,QX) \equiv \Hom_{\2E_0}(``\oj XJ",``\oj XJ").
\eean
\\
\noindent{\it Horizontal ($\times$-) Composition of 2-morphisms.}
Let $\6E, \6F, \6G\in \{\6A,\6B\}, Z_1, Z'_1\in\Hom_{\2E_0}(\6E,\6F)$, 
$Z_2, Z'_2\in\Hom_{\2E_0}(\6F,\6G)$ and $s_i\in\Hom_{\2E_0}(Z_i,Z'_i), i=1,2$. 

\[
\6E\ctwo{Z_1}{Z_1'}{s_1}%
\6F\ctwo{Z_2}{Z_2'}{s_2}%
\6G
\]

Then we must define 
$s_2\times s_1\in\Hom_{\2E_0}(Z_2\circ Z_1,Z'_2\circ Z'_1)$. To this purpose we
consider, as before, $s_i$  as morphisms in $\2A$. Thus for $\6E=\6F=\6A$ we have 
$Z_1, Z_1'\in\obj\2A$ and $s_1\in\Hom_\2A(Z_1,Z_1')$, and for $\6E=\6B, \6F=\2A$ we have
$Z_1=X_1J, Z_1'=X_1'J$ with $X_1, X_1'\in\obj\2A$ and 
$s_1\in\Hom_\2E(X_1j,X'_1J)\equiv \Hom_\2A(X_1Q,X_1')$ etc. We define 
\bea s_2\times s_1= s_2\otimes s_1 & \mbox{if} & \6F=\6A, \nn\\
   s_2\times s_1= \id_? \otimes w'\otimes\id_? \mcirc s_2\otimes\id_Q\otimes s_1 \mcirc 
  \id_? \otimes w\otimes\id_? & \mbox{if} & \6F=\6B, \label{t2}\eea
To illustrate the second equation, consider the case $\6E=\6F=\6G=\6B$, thus
$Z^{(')}_i=\oj \circ X^{(')}_i\circ J, i=1,2$, with $X^{(')}_i\in\2A$ and
$s_i\in\Hom_\2A(X_iQ,QX'_i)$. Then (\ref{t2}) looks like:
\[ 
\begin{tangle}
\object{Q}\step\object{X_2'}\step[1.5]\object{Q}\step[1.5]\object{X_1'}\\
\id\step\id\step\hcd\step\id\\
\hh\id\step\id\step\id\step\id\step\id\\
\hh\frabox{s_2}\step\id\step\frabox{s_1}\\
\hh\id\step\id\step\id\step\id\step\id\\
\id\step\hcu\step\id\step\id\\
\object{X_2}\step[1.5]\object{Q}\step[1.5]\object{X_1}\step\object{Q}
\end{tangle}
\]

Our task is now two-fold. On the one hand we must convince the reader that this is the
`right' definition. We remark that the definition (\ref{t2}) is motivated by the
interpretation of the 2-morphisms of $\2E_0$ in terms of morphisms in $\2A$. For the
horizontal composition of 2-morphisms where the intermediate object $\6F$ is $\6B$ this
looks as follows: 
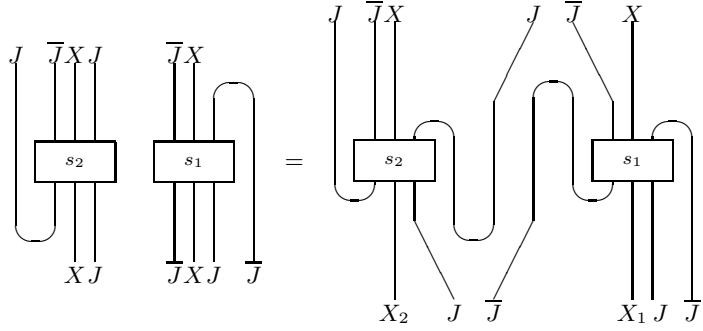
\[ 
\begin{tangle}
\object{J}\step\object{\oj }\hstep\object{X}\hstep\object{J}\Step\object{\oj }\hstep\object{X}\\
\hh\id\step\id\hstep\id\hstep\id\Step\id\hstep\id\hstep\hcoev\\
\hh\id\step\id\hstep\id\hstep\id\Step\id\hstep\id\hstep\id\step\id\\
\hh\id\step\frabox{s_2}\Step\frabox{s_1}\step\id\\
\hh\id\step\id\hstep\id\hstep\id\Step\id\hstep\id\hstep\id\step\id\\
\hh\hev\hstep\id\hstep\id\Step\id\hstep\id\hstep\id\step\id\\
\step[1.5]\object{X}\hstep\object{J}\Step\object{\oj }\hstep\object{X}\hstep\object{J}\step\object{\oj }
\end{tangle}
\quad = \quad
\begin{tangle}
\object{J}\step\object{\oj }\hstep\object{X}\step[3.5]\object{J}\step\object{\oj }\step[1.5]\object{X}\\
\id\step\id\hstep\id\step[2.5]\dd\hcoev\d\hstep\id\\
\hh\id\step\id\hstep\id\hstep\hcoev\step\id\step\id\step\id\step\id\hstep\id\hstep\hcoev\\
\hh\id\step\frabox{s_2}\step\id\step\id\step\id\step\id\step\frabox{s_1}\step\id\\
\hh\hev\hstep\id\hstep\id\step\id\step\id\step\id\step\hev\hstep\id\hstep\id\step\id\\
\step[1.5]\id\hstep\d\hev\dd\step[2.5]\id\hstep\id\step\id\\
\step[1.5]\object{X_2}\step[1.5]\object{J}\step\object{\oj }\step[3.5]\object{X_1}\step[.7]\object{J}\step[.8]\object{\oj }
\end{tangle}
\]
It should be clear that for $\6F=\6A$ the first formula in (\ref{t2}) is the correct
definition. 

The second task is, of course, to prove that our definition renders $\2E_0$ an
almost-2-category. This means the horizontal composition of three 2-morphisms must be
associative for all legal compositions. Luckily, this is quite obvious from associativity
of the $\otimes$-composition in the tensor category $\2A$ and we refrain from
formalizing this. It remains to show the interchange law, which again we do only for
$\6E=\6F=\6G=\6B$, all other cases being easier. Let thus $Z_i=JX_i\oj $, 
$s_1\in\Hom_\2E(Z_1,Z_1')\equiv\Hom_\2A(X_1Q,QX_1')$ etc.  
We compute the composition
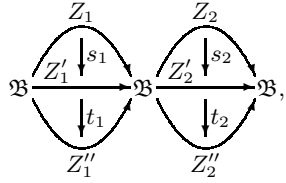
\[
\6B%
\cthree{Z_1}{Z_1'}{Z_1''}{s_1}{t_1}%
\6B
\cthree{Z_2}{Z_2'}{Z_2''}{s_2}{t_2}%
\6B,
\]
which is in $\Hom_\2E(Z_ZZ_1,Z_2''Z_1'')\equiv\Hom_\2A(X_2QX_1Q,QX_2''QX_1'')$, in two
different ways. We can first do the horizontal compositions and obtain

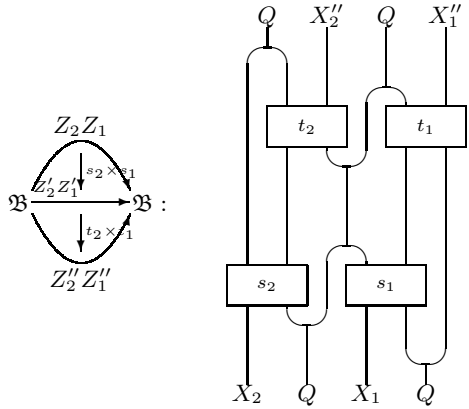
\[
\6B%
\cthree{Z_2Z_1}{{\scriptstyle Z_2'Z_1'}}{Z_2''Z_1''}{{\scriptscriptstyle s_2\times s_1}}{{\scriptscriptstyle t_2\times t_1}}%
\6B
: \quad \quad \quad 
\begin{tangle}
\hstep\object{Q}\step[1.5]\object{X_2''}\step[1.5]\object{Q}\step[1.5]\object{X_1''}\\
\hh\hcd\step\id\step[1.5]\id\step[1.5]\id\\
\hh\id\step\id\step\id\step\hcd\step\id\\
\hh\id\step\frabox{t_2}\step\id\step\frabox{t_1}\\
\id\step\id\step\hcu\step\id\step\id\\
\hh\id\step\id\step\hcd\step\id\step\id\\
\hh\frabox{s_2}\step\id\step\frabox{s_1}\step\id\\
\hh\id\step\hcu\step\id\step\id\step\id\\
\hh\id\step[1.5]\id\step[1.5]\id\step\hcu\\
\object{X_2}\step[1.5]\object{Q}\step[1.5]\object{X_1}\step[1.5]\object{Q}
\end{tangle}
\]

Beginning with the vertical compositions we arrive at
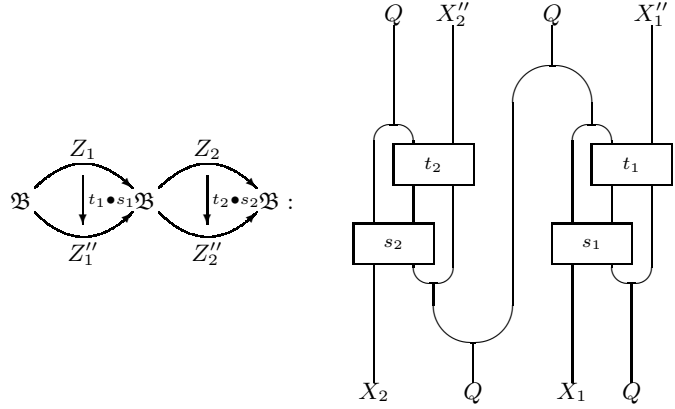
\[
\6B\ctwo{Z_1}{Z_1''}{{\scriptstyle t_1\bullet s_1}}%
\6B\ctwo{Z_2}{Z_2''}{{\scriptstyle t_2\bullet s_2}}%
\6B
: \quad \quad \quad 
\begin{tangle}
\hstep\object{Q}\step[1.5]\object{X_2''}\step[2.5]\object{Q}\step[2.5]\object{X_1''}\\
\hstep\id\step[1.5]\id\step[1.5]\cd\step[1.5]\id\\
\hh\hcd\step\id\step[1.5]\id\step[1.5]\hcd\step\id\\
\hh\id\step\frabox{t_2}\step[1.5]\id\step[1.5]\id\step\frabox{t_1}\\
\hh\id\step\id\step\id\step[1.5]\id\step[1.5]\id\step\id\step\id\\
\hh\frabox{s_2}\step\id\step[1.5]\id\step[1.5]\frabox{s_1}\step\id\\
\hh\id\step\hcu\step[1.5]\id\step[1.5]\id\step\hcu\\
\id\step[1.5]\cu\step[1.5]\id\step[1.5]\id\\
\object{X_2}\step[2.5]\object{Q}\step[2.5]\object{X_1}\step[1.5]\object{Q}
\end{tangle}
\]
That the two expressions coincide is again verified easily using (\ref{W1}-\ref{W3}).

Assume now that $\2A$ is $\7F$-linear. Clearly, the spaces of 2-morphisms in $\2E_0$
are $\7F$-vector spaces and the compositions $\bullet, \times$ are bilinear. Thus $\2E_0$
is $\7F$-linear. Finally, let $s: Q\rarr \tilde{Q}$ be an isomorphism between the
Frobenius algebras $\5Q, \tilde{\5Q}$ and consider the almost 2-categories $\2E_0,
\tilde{\2E}_0$ constructed from $\5Q, \tilde{\5Q}$. The objects and 
1-morphisms (as well as the composition of the latter) of $\2E_0, \tilde{\2E}_0$ not
depending on the Frobenius algebra, there are obvious bijections. Furthermore, there is a
bijection between, e.g., $\Hom_{\2E_0}(``XJ",``YJ")\equiv\Hom_\2A(XQ,Y)$ and 
$\Hom_{\tilde{\2E}_0}(``XJ",``YJ")\equiv\Hom_\2A(X\tilde{Q},Y)$ given by
\[ \Hom_\2A(XQ,Y)\ni t \mapsto t\circ s^{-1}\otimes\id_Y   \in\Hom_\2A(X\tilde{Q},Y). \]
Since these isomorphisms commute with the compositions $\bullet, \times$ of 2-morphisms
in $\2E_0, \tilde{\2E}_0$ we have an isomorphism of almost 2-categories.
\qed

\brem 1. It is obvious how to modify the proposition if $\2A$ is non-strict: The
definition of the objects and 1-morphisms and the composition of the latter are
unchanged. Since we may still require $J\oj=Q$, the associativity constraint of $\2A$
gives rise to that of $\2E_0$. As to the 2-morphisms, the only change are appropriate
insertions of associativity morphisms in the definitions of $\bullet, \times$.

2. If $\5Q=(Q,\ldots)$ is a Frobenius algebra in a tensor category $\2A$ then the functor 
$F=Q\otimes-$ is part of a Frobenius algebra in $\End\,\2A$, thus in particular of a
monoid in $\End\,\2A$, equivalently a monad
$(Q\otimes-,\{w'\otimes\id_X\},\{v\otimes\id_X\})$ in $\2A$. It is easy to verify that
our construction of the category $\HOM_{\2E_0}(\6A,\6B)$ is precisely the Kleisli 
construction \cite[Section VI.5]{cwm} starting from $\2A$ and the monad $(F,\ldots)$. 
(Alternatively one may invoke the Kleisli type construction for monoids in 2-categories,
cf.\ e.g.\ \cite{gray}.) Similar statements hold for the categories $\HOM_\2E(\6B,\6A)$
and $\HOM_\2E(\6B,\6B)$. But our way of pasting everything together in order to obtain a 
bicategory seems to be new. 

3. Just as the Kleisli category is the smallest solution to a certain problem (i.e.\ an
initial object in the category of all adjunctions producing the given monad), it is
intuitively clear that $\2E_0$ is an initial object in the category of all solutions of
our problem. Consider a 2-category $\2E'$ with $\{ \6A', \6B'\} \subset\obj\2E'$ such that
$\HOM_{\2E'}(\6A',\6A')\cong\2A$ (with given isomorphism) and with mutually dual
$J'\in\Hom_{\2E'}(\6B',\6A'), \oj'\in\Hom_{\2E'}(\6A',\6B')$, such that $Q=J'\oj'$. Then
there is a unique functor $K$ of almost-2-categories $K:\2E_0\rarr \2F$ such that
$F(\6A)=\6A', F(\6B)=\6B', F(J)=J', F(\oj )=\oj'$. We omit the details, since in Theorem
\ref{t-uniq} we will prove a more useful uniqueness result.
\erem


With the preceding constructions $\2E_0$ is a 2-category up to one defect: there is no
unit-$\6B-\6B$-morphism. We could, of course, try to add one by hand but that would be
difficult to do in a consistent manner. Fortunately, it turns out that taking the closure
of $\2E_0$ in which idempotent 2-morphisms split automatically provides us with a
(non-strict) $\6B-\6B$-unit $\11_\6B$ provided the Frobenius algebra $(Q,v,v',w,w')$
satisfies an additional condition.

In all applications we are going to discuss $\2A$ is linear over a field $\7F$ and
$\End(\11)\cong\7F$. Yet we wish to emphasize the generality of our basic construction.
This is why we give the following definition,  motivated by considerations
in \cite[p.\ 72-3]{t}.

\begin{defprop}
A strict (not necessarily preadditive) tensor category $\2A$ is $\End(\11)$-linear
if $\lambda\otimes s=s\otimes\lambda=:\lambda s$ for all $\lambda\in\End(\11)$ and 
$s: X\rarr Y$. It then follows that 
$\lambda(s\circ t)=(\lambda s)\circ t=s\circ(\lambda t)$ and
$\lambda(s\otimes t)=(\lambda s)\otimes t=s\otimes(\lambda t)$. (All this generalizes to
non-strict categories.)
\end{defprop}

\btheor \label{main0}
Let $\2A$ be a strict tensor category and $\5Q=(Q,v,v',w,w')$ a Frobenius
algebra in $\2A$. Assume that one of the following conditions is satisfied:
\begin{itemize}
\item[(a)] $\displaystyle w'\circ w = \id_Q$.
\item[(b)] $\2A$ is $\End(\11)$-linear and 
\[ w'\circ w = \lambda_1 \, \id_Q, \]
where $\lambda_1$ is an invertible element of the commutative monoid $\End(\11)$. 
\end{itemize}
Then the completion $\2E=\ol{\2E_0}^p$ of the $\2E_0$ defined in Proposition \ref{constr} is a
bicategory such that 
\begin{enumerate}
\item $\displaystyle \obj\,\2E=\{\6A,\6B\}$.
\item There is a fully faithful tensor functor $I: \2A\rarr\HOM_\2E(\6A,\6A)$ such that
for every $Y\in\HOM_\2E(\6A,\6A)$ there is $X\in\2A$ such that $Y$ is a retract of $I(X)$.
(Thus $I$ is an equivalence if $\2A$ has subobjects.)
\item There are 1-morphisms $J: \6B\rarr \6A$ and $\oj : \6A\rarr \6B$ such that
$Q=J\oj $.  
\item $J$ and $\oj $ are mutual two-sided duals, i.e.\ there are 2-morphisms
\[ e_J: \11_\6A\rarr J\oj,\ \ \ve_J: \11_\6B\rarr\oj J,\ \ d_J: \oj J\rarr\11_\6B,\ \ 
   \eta_J: J\oj\rarr\11_\6A \]
satisfying the usual relations.
\item We have the identity 
\[ I(Q,v,v',w,w')=(J\oj,e_J,\eta_J,\id_J\otimes\ve_J\otimes\id_{\oj},\id_J\otimes
   d_J\otimes\id_{\oj }) \]
of Frobenius algebras in $\END_\2E(\6A)$.
(In particular, $d_J\circ\ve_J=\lambda_1\id_{\11_\6B}$.)
\item If $\2A$ is a preadditive ($\7F$-linear) category then $\2E$ is a preadditive
($\7F$-linear) 2-category.
\item If $\2A$ has direct sums then $\2E$ has direct sums of 1-morphisms.
\end{enumerate}
Isomorphic Frobenius algebras $\5Q, \tilde{\5Q}$ give rise to isomorphic bicategories
$\2E, \tilde{\2E}$.
\etheor
\prf If we are in case (a) put $\lambda_1=1$ in the sequel. Then $\End(\11)$-linearity will
not be needed. We define the bicategory $\2E$ as the completion $\ol{\2E_0}^p$. Thus 
$\obj\,\2E=\{\6A, \6B\}$ and for $\6X, \6Y\in\{\6A, \6B\}$ the 1-morphisms are
\[ \Hom_{\2E}(\6X,\6Y)= \{ (X,p) \ \vert\ X\in\Hom_{\2E_0}(\6X,\6Y),\
  p=p\bullet p\in\Hom_{\2E_0}(X,X) \}.  \]
Furthermore,
\bean \Hom_{\2E}((X,p),(Y,q)) & = & \{ s\in\Hom_{\2E_0}(X,Y) \ \vert\ s\bullet p=q\bullet s=s
   \} \\   &=& q\bullet\Hom_{\2E_0}(X,Y)\bullet p. \eean
In order to alleviate the notation we allow $X$ to denote also $(X,\id_X)$. With this
definition it is clear that $(X,p)\prec X\equiv(X,\id_X)$. ($p\in\Hom_{\2E}(X,X)$ is an
invertible  morphism from $(X,p)$ to $X$, since also $p\in\Hom_{\2E}(X,(X,p))$ and
$p\bullet p=p=id_{(X,p)}$.) 
\\ \\
{\it Exhibiting a unit $\6B-\6B$-morphism.}
Recall that $\End_\2E(\oj J)\equiv\End_\2A(Q)$ as a vector space and consider the morphism 
$p_1\in\End_\2E(\oj J)$ represented by $\lambda_1^{-1}\id_Q\in\End_\2A(Q)$. (Note that
$\id_Q$ is the unit of the monoid $\End_\2A(Q)$, but not of $\End_{\2E}(\oj J)$, whose
unit is $v\circ v'$!) As a consequence of condition (b) we see that 
$p_1\bullet p_1=\lambda_1^{-2}w'\circ w=\lambda_1^{-1}\id_Q=p_1$. Thus $p_1$ is idempotent
and $(\oj J, p_1)$ is a $\6B-\6B$-morphism in $\2E$. We claim that $(\oj J, p_1)$
is a (non-strict) unit $\6B-\6B$-morphism. In order to see that
$\11_\6B$ is a right unit we have to show that there are isomorphisms 
\bean r((XJ,p)) &:&  (XJ,p)\11_\6B=(XQJ,p\times p_1)\rarr(XJ,p), \\
   r((\oj XJ,p)) &:& (\oj XJ,p)\11_\6B=(\oj XQJ,p\times p_1)\rarr(\oj XJ,p) \eean
for all $(XJ, p): \6B\rarr\6A$ and $(\oj XJ,p): \6B\rarr\6B$, respectively. We consider
only $r((\oj XJ,p))$ and leave the other case to the reader. (For $r((XJ,p))$ the only
change is that the upper left Q-leg of $p$ disappears.)

By definition of the horizontal composition of 2-morphisms in $\2E$ we have
\[ p\times p_1 = \lambda_1^{-1}\quad
\begin{tangle}
\object{Q}\step\object{X}\step[1.5]\object{Q}\\
\hh\id\step\id\step\hcd\\
\hh\frabox{p}\step\id\step\id\\
\hh\id\step\id\step\id\step\id\\
\hh\id\step\hcu\step\id\\
\object{X}\step[1.5]\object{Q}\step[1.5]\object{Q}
\end{tangle}
\quad = \lambda_1^{-1}\quad
\begin{tangle}
\object{Q}\step\object{X}\step\object{Q}\\
\hh\id\step\id\step\id\\
\hh\frabox{p}\step\id\\
\hh\id\step\id\step\id\\
\hh\id\step\hcu\\
\hh\id\step\hcd\\
\object{X}\step\object{Q}\step\object{Q}
\end{tangle}
\] 

Consider
\[ r((\oj XJ,p))= \quad
\begin{tangle}
\object{Q}\step\object{X}\\
\id\step\id\step[1.5]\unit\\
\hh\id\step\id\step\hcd\\
\hh\frabox{p}\step\id\step\id\\
\hh\id\step\id\step\id\step\id\\
\hh\id\step\hcu\step\id\\
\object{X}\step[1.5]\object{Q}\step[1.5]\object{Q}
\end{tangle}
\quad = \quad
\begin{tangle}
\object{Q}\step\object{X}\\
\hh\id\step\id\\
\hh\frabox{p}\\
\hh\id\step\id\\
\hh\id\hstep\hcd\\
\object{X}\hstep\object{Q}\step\object{Q}
\end{tangle}
\quad \in\Hom_\2A(XQQ,QX)\equiv\Hom_\2E(\oj XQJ,\oj XJ)
\] 
\[ r'((\oj XJ,p))= \quad
\begin{tangle}
\object{Q}\step\object{X}\step[1.5]\object{Q}\\
\hh\id\step\id\step\hcd\\
\hh\frabox{p}\step\id\step\id\\
\hh\id\step\id\step\id\step\id\\
\hh\id\step\hcu\step\id\\
\id\step[1.5]\counit\step[1.5]\id\\
\object{X}\step[3]\object{Q}
\end{tangle}
\quad = \quad
\begin{tangle}
\object{Q}\step\object{X}\step\object{Q}\\
\hh\id\step\id\step\id\\
\hh\frabox{p}\step\id\\
\hh\id\step\id\step\id\\
\id\step\hcu\\
\object{X}\step[1.5]\object{Q}
\end{tangle}
\quad \in\Hom_\2A(XQ,QXQ)\equiv\Hom_\2E(\oj XJ,\oj XQJ)
\] 
Using $p\bullet p=p$ and the rules of computation in $\2E$ it is easy to verify that 
$r\bullet r'=\lambda_1 p=\lambda_1\id_{(\oj XJ,p)}$ and 
$r'\bullet r=\lambda_1 p\times p_1=\lambda_1\id_{(\oj XQJ,p\times p_1)}$, such that 
$r((\oj XJ,p))$ is an isomorphism from $(\oj XQJ,p\times p_1)$ to $(\oj XJ,p)$. We leave
this computation to the reader as an exercise. That $r((\oj XJ,p))$ is natural w.r.t.\
$(\oj XJ,p)$ is easy, and the coherence law connecting the unit constraint with the tensor
product becomes almost obvious since $\2E$ is strict except for the unit morphism under
study. That $\11_\6B$ is a left unit is shown by a similar argument defining 
$l((\oj XJ,p))\in\Hom_\2E(\oj QXJ,\oj XJ)\equiv\Hom_\2A(QXQ,QX)$ analogously.
Finally, one sees that for $(\oj XJ,p)=\11_\6B$ the left and right unit constraints
coincide: $l(\11_\6B)=r(\11_\6B)$. \\ \\
\noindent{\it Duality of $J$ and $\oj $.} 
We refer to \cite[Section I.6]{gray} for a discussion of adjoint 1-morphisms in
(non-strict) bicategories. In order to show that $J, \oj $ are two-sided dual 1-morphisms
we must exhibit morphisms 
\[ e_J: \11_\6A\rarr J\oj ,\ \  d_J: \oj J\rarr\11_\6B, \ \ 
  \ve_J: \11_\6B\rarr \oj J,\ \ \eta_J: J\oj \rarr\11_\6A \]
satisfying the usual triangular equations.
Motivated by $Q=J\oj $ and by Lemma \ref{motiv} we set 
\bea e_J &=& v \quad \in\Hom_\2A(\11,Q)\equiv\Hom_\2E(\11_\6A,J\oj ), \label{vebj}\\
   \eta_J &=& v'  \quad \in\Hom_\2A(Q,\11)\equiv\Hom_\2E(J\oj ,\11_\6A). \nn\eea
Now we observe that 
\[ \Hom_\2E(\11_\6B,\oj J)\subset\Hom_{\2E_0}(\oj J,\oj J)\equiv\Hom_\2A(Q,Q), \]
\[ \Hom_\2E(\oj J,\11_\6B)\subset\Hom_{\2E_0}(\oj J,\oj J)\equiv\Hom_\2A(Q,Q). \]
(By construction of $\2E$, $\11_\6B$ is the retract of $\oj J: \6B\rarr\6B$ corresponding
to the idempotent $\lambda_1^{-1}\id_Q\in\End_\2A(Q)\equiv\End_{\2E_0}(\oj J)$.)
Thus it is reasonable to consider the following candidates for $d_J$ and
$\ve_J$ (both of which live in $\Hom_\2A(Q,Q)$):
\begin{equation} \label{vej}
   d_J=\id_Q,\quad \lambda_1^{-1}\ve_J= \id_Q.
\end{equation}
(Whether $d_J$ or $\ve_J$ contains the factor $\lambda_1^{-1}$ is immaterial, but the
normalization of the left and right unit constraints $l$ and $r$ depends on this choice.)
With this definition we have 
\[ d_J\bullet\ve_J=\lambda_1^{-1} w'\circ w=\id_Q=\lambda_1 p_1=\lambda_1\id_{\11_\6B} \]
as desired. In the verification of the triangular equations we must be aware that $\2E$ is
only a bicategory since there are non-trivial unit constraints for $\11_\6B$. The
computations tend to be somewhat confusing. We prove only one of the four equations,
namely that
\be
\begin{diagram}
J \equiv \11_\6A J  & \rTo^{e_J\otimes\id_J} & J\oj J & 
   \rTo^{\id_J\otimes d_J} & J\11_\6B & \rTo^{r(J)} & J
\end{diagram}
\label{ww}\end{equation}
is the identity 2-morphism $\id_J$, the computation being completely analogous in the other
cases. With $\id_J=v'\in\Hom_\2A(Q,\11)\equiv\Hom_\2E(J,J)$, (\ref{vebj}), and (\ref{vej})
we compute
\bean e_J\times\id_J &=& \ \ 
\begin{tangle} \counit\step\unit \end{tangle}\ \  \in\Hom_\2A(Q,Q)\equiv\Hom_\2E(J,J\oj J)\\
   \id_J\times d_J &=& \
\begin{tangle}
\unit\step\hcd\\
\hcu\step\id
\end{tangle} 
\ \ = \
\begin{tangle} \hcd \end{tangle} 
\in\Hom_\2A(Q^2,Q)\equiv\Hom_\2E(J\oj J,J\oj J)
\eean
$\bullet$-Composing these 2-morphisms between $\6B\rarr\6A$-morphisms according to the
rules of Proposition \ref{constr} we obtain
\[ 
\begin{tangle}
\cd\\
\counit\step\unit\step\id\\
\hh\step\hcu
\end{tangle}\quad = \quad
\begin{tangle}
\id\\
\object{Q}
\end{tangle}
\]
This is precisely the isomorphism $r'(J): J\11_\6B\rarr J$ given by 
$r'(J)=\id_Q\in\Hom_\2A(Q,Q)\equiv\Hom_\2E(J\oj J,J)$ provided by Theorem \ref{main0}.
Now $r(J), r'(J)$ are mutually inverse, which proves that (\ref{ww}) gives the unit
morphism $\id_J$. 
The last statement is obvious since isomorphic (almost) bicategories 
$\2E_0, \tilde{\2E}_0$ have isomorphic completions $\ol{\2E_0}^p, \ol{\tilde{\2E}_0}^p$. 
\qed

\brem 1. The bicategory $\2E$ fails to be strict (thus a 2-category) only due to the
presence of non-trivial unit constraints for $\11_\6B$. This defect could be repaired by
adding a strict unit 1-morphism for $\6B$ which is isomorphic to $(\oj J, p_1)$. There
will, however, be no compelling reason to do so. 

2. The condition (a/b) in Theorem \ref{main0} was crucial for identifying the unit
1-morphism $\11_\6B$ as a retract of $\oj J$. Furthermore, we obtained a distinguished  
retraction/section $\11_\6B\leftrightarrow \oj J$.
So far, our assumptions are not symmetric in that they do not imply $\11_\6A\prec J\oj$,
let alone provide a canonical retraction and section. This is achieved by the following
definition. 
\erem

\bdefin \label{d-canon}
Let $\2A$ be an $\End(\11)$-linear (but not necessarily a preadditive) category.
A Frobenius algebra $\5Q=(Q,v,v',w,w')$ in $\2A$ is canonical iff 
\bea w'\circ w &=& \lambda_1 \id_Q, \label{W4} \\
   v'\mcirc v &=& \lambda_2, \label{W5}\eea
where $\lambda_1, \lambda_2\in\End(\11)$ are invertible. If $\lambda_1=\lambda_2$ then 
$\5Q$ is called normalized. 
\edefin

\brem 1. The notion `canonical Frobenius algebra is motivated by subfactor theory, where
canonical endomorphisms are the objects of canonical Frobenius algebras in certain tensor 
categories, see Section \ref{ss-subfact}. 

2. If $\alpha,\beta\in\End(\11)^*$ and $\5Q=(Q,v,v',w,w')$ is a canonical Frobenius
algebra then clearly also $\tilde{\5Q}=(Q,\alpha v,\beta^{-1}v',\beta w,\alpha^{-1}w')$ is
one. $\tilde{\5Q}$ is isomorphic to $\5Q$ (in the sense of Definition \ref{d-isom}) iff
$\alpha=\beta$, in which case an isomorphism is given by$s=\alpha\,\id_Q:Q\rarr\tilde{Q}=Q$. 
Yet we consider Frobenius algebras related by this renormalization as equivalent. Note
that $\tilde{\lambda}_1\tilde{\lambda}_2=\lambda_1\lambda_2$, thus 
$v'\circ w'\circ w\circ v\in\End(\11)$ is invariant under isomorphism and renormalization.
\erem

From now on all tensor categories are assumed to be $\End(\11)$-linear. Thus if we state
this explicitly it is only for emphasis.

\bprop \label{connect}
Let $\2E$ be a bicategory and $J:\6B\rarr\6A$ a 1-morphism with two-sided dual
$\oj:\6A\rarr\6B$. Assume that the corresponding 2-morphisms $d_J, e_J, \ve_J, \eta_J$ can
be chosen such that $\eta_J\circ e_J$ and $d_J\circ \ve_J$ are invertible in the monoids 
$\End(\11_\6A), \End(\11_\6B)$, respectively. Then the functor 
$F=-\otimes J:\ \HOM_\2E(\6A,\6A)\rarr\HOM_\2E(\6B,\6A)$ is faithful and do\-mi\-nant in the
sense that every $X: \6B\rarr\6A$ is a retract of $Y\circ J$ for some $Y:\6A\rarr\6A$.
The same holds for the other seven functors given by composition with $J$ or $\oj$ from
the left or right.  
\eprop
\prf Our conditions obviously imply that $e_X: \11_\6A\rarr J\oj$ and 
$\ve_X: \11_\6A\rarr\oj J$ are retractions, viz.\ have left inverses. Thus
$\11_\6A\prec\oj J$ and therefore $X\prec X\circ(\oj\circ J)\cong (X\circ\oj)\circ J$ for
any $X:\6B\rarr\6A$. Since $X\circ\oj$ is a $\6A-\6A$ morphism this implies the dominance
of $F$. Faithfulness can be proved using \cite[Theorem IV.3.1]{cwm}, but we prefer to give a 
direct argument. Let $X,Y: \6A\rarr\6A$ and $s\in\Hom_\2E(X,Y)$. If $s\otimes\id_J=0$ then
also $s\otimes\id_J\otimes\id_\oj=s\otimes\id_Q=0$. Sandwiching between $\id_Y\otimes v'$
and $\id_X\otimes v$ gives $s\otimes(v'\circ v)=\lambda_2 s=0$ and thus $s=0$ by
invertibility of $\lambda_2$. 
\qed

\bcoro \label{coro1}
Let $\2A$ be $\End(\11)$-linear and $(Q,v,v',w,w')$ a canonical Frobenius algebra
in $\2A$. Then the bicategory $\2E$ defined above is such that 
$\eta_J\circ e_J$ and $d_J\circ \ve_J$ are invertible in the monoids $\End(\11_\6A),
\End(\11_\6B)$, respectively. Conversely, if $J: \6B\rarr\6A$ has a two-sided dual $\oj$
such that $e_J, d_J, \ve_J, \eta_J$ satisfy these conditions then the Frobenius algebras
$(J\oj,\ldots)$ and $(\oj J,\ldots)$ in $\END_\2E(\6A), \END_\2E(\6B)$, respectively,  are
canonical. 
\ecoro
\prf By Theorem \ref{main0}, $d_J\circ\ve_J=\lambda_1\id_Q$ with
$\lambda_1\in\End(\11)^*$. On the other hand, $\eta_J\circ e_J=v'\circ v$ which is
invertible in $\End(\11)$ since $\5Q$ is canonical.
The converse is obvious in view of Lemma \ref{motiv}.
\qed

Now we are in a position to consider the uniqueness of our bicategory $\2E$.
\btheor \label{t-uniq}
Let $\2A$ be $\End(\11)$-linear and $(Q,v,v',w,w')$ a canonical Frobenius algebra
in $\2A$. Let $\2E$ be as constructed in Theorem \ref{main0} and let $\tilde{\2E}$ be any
bicategory such that
\begin{enumerate}
\item $\displaystyle\obj\,{\tilde{\2E}}=\{\6A,\6B\}$.
\item Idempotent 2-morphisms in $\tilde{\2E}$ split.
\item There is a fully faithful tensor functor $\tilde{I}: \2A\rarr\END_{\tilde{\2E}}(\6A)$ 
such that every object of $\END_{\tilde{\2E}}(\6A)$ is a retract of $\tilde{I}(X)$ for some
$X\in\2A$.
\item There are mutually two-sided dual 1-morphisms $\tilde{J}: \6B\rarr\6A$,
$\tilde{\oj}: \6A\rarr\6B$ and an isomorphism $\tilde{s}: I(Q)\rarr\tilde{J}\tilde{\ol{J}}$
between the Frobenius algebras $I(Q,v,v',w,w')$ and 
$(\tilde{J}\tilde{\ol{J}}, \tilde{e}_{\tilde{J}},\ldots)$ in $\End_{\tilde{\2E}}(\6A)$.
\end{enumerate}
Then there is an equivalence $E: \2E\rarr\tilde{\2E}$ of bicategories such that there is a
monoidal natural isomorphism between the tensor functors 
$\tilde{I}$ and $(E\restr\END_\2E(\6A))\circ I$.
\etheor
\prf Replacing $\2E, \tilde{\2E}$ by equivalent bicategories, we may assume that the
functors $I: \2A\rarr\END_\2E(\6A), \tilde{I}: \2A\rarr\END_{\tilde{\2E}}(\6A)$ are
injective on the objects. Thus $\END_\2E(\6A)$ and $\END_{\tilde{\2E}}(\6A)$ contain $\2A$
as a full subcategory. In view of the coherence theorem for bicategories we may replace 
$\2E,\tilde{\2E}$ by equivalent strict bicategories or 2-categories and we suppress the
symbols $I,\tilde{I}$. (As a consequence of these replacements, we will no more have the
identity $I(Q)=J\oj$ but only an isomorphism $s: Q\equiv I(Q)\rarr J\oj$ compatible with
the Frobenius algebra structures.) 
In view of Proposition \ref{connect}, every $Y\in\Hom_\2E(\6B,\6A)$ is a retract of $Y\oj J$
and therefore of $XJ$ for $X=Y\oj\in\2A$. Similarly, every $Z\in\Hom_\2E(\6A,\6B)$ 
($Z\in\Hom_\2E(\6B,\6B)$) is a retract of $\oj X$ ($\oj XJ$) for some $X\in\2A$, and
similarly for $\tilde{\2E}$. Let $\2E_0$ be the full sub 2-category of $\2E$ with objects
$\{\6A,\6B\}$ and 1-morphisms $X, XJ,\oj X,\oj XJ$ with $X\in\2A$, and similarly for
$\tilde{\2E}$. Now we can define $E: \2E_0\rarr\tilde{\2E}_0$ as the identity on objects
and 1-morphisms. Composing the obvious isomorphisms 
$\Hom_{\2E_0}(XJ,YJ)\cong\Hom_\2A(XQ,Y)\cong\Hom_{\tilde{\2E}_0}(X\tilde{J},Y\tilde{J})$,
etc., provided by the duality of $J,\oj$ and $\tilde{J},\tilde{\oj}$ we can define the
functor $E_0$ on the 2-morphisms. That $E$ commutes with the horizontal and vertical
compositions is obvious by the isomorphism 
$(J\oj,e_J,\ldots)\cong(Q,v,\ldots)\cong(\tilde{J}\tilde{\oj},\tilde{e}_{\tilde{J}},\ldots)$
of Frobenius algebras. In order to obtain a (non-strict) isomorphism $E$ of (strict)
bicategories we need to define invertible 2-cells $\phi_{gf}: Eg\circ Ef\rarr E(g\circ f)$
satisfying the usual conditions \cite{benab}. When $\mbox{Ran}\,f=\mbox{Src}\,g=\6A$ we
choose them to be identities and for $\mbox{Ran}\,f=\mbox{Src}\,g=\6B$ we use the isomorphism 
$\tilde{s}\circ s^{-1}: J\oj\rarr\tilde{J}\tilde{\oj}$. The verification of the coherence
conditions is straightforward but very tedious to write down, and therefore omitted.
In view of $\2E\simeq\ol{\2E_0}^p\cong\ol{\tilde{\2E}_0}^p\simeq\tilde{\2E}$
the isomorphism $E_0: \2E_0\rarr\tilde{\2E}_0$ extends to an equivalence
$E: \2E\rarr\tilde{\2E}$ which has all desired properties.
\qed

\brem \label{r-alt}
1. The construction of the bicategory $\2E$ given in this section reflects the author's 
understanding as of 1999. More recently, it has become clear that there exists an
alternative construction which can be stated quite succinctly. Namely, define a 
bicategory $\tilde{\2E}$ with $\obj\,\tilde{\2E}=\{\6A,\6B\}$ by positing
\begin{eqnarray*} \HOM_{\tilde{\2E}}(\6A,\6A) &=& \2C, \\
\HOM_{\tilde{\2E}}(\6A,\6B) &=& Q-\mod, \\
\HOM_{\tilde{\2E}}(\6B,\6A) &=& \mod-Q, \\
\HOM_{\tilde{\2E}}(\6B,\6B) &=& Q-\mod-Q, 
\end{eqnarray*}
where $Q-\mod,\mod-Q,Q-\mod-Q$ are the categories of left and right $Q$-modules and of
$Q-Q$-bimodules in $\2C$, respectively. It is very easy to prove that these categories are
abelian if $\2C$ is abelian, in particular they `have subobjects', i.e.\ idempotents split.
The compositions of 1-morphisms are defined as suitable quotients like in ring theory. The
two different constructions are related to each other as the constructions of Kleisli and
Eilenberg/Moore. That they lead to equivalent bicategories is, at least morally, due to
the completion w.r.t.\ subobjects which we apply in our construction. (Module categories
automatically have subobjects.) To prove this rigorously it is, in view of Theorem
\ref{t-uniq}, sufficient to prove that $\tilde{\2E}$ satisfies the requirements of the
latter. Inspired by \cite{fs} (which in turn was influenced by ideas of the author), the
alternative definition of $\tilde{\2E}$ was proposed also by S. Yamagami in \cite{yamag4},
which he kindly sent me.

2. By the above remark, our category $\2B=\END_\2E(\6B)$ is equivalent to the bimodule
category $Q-\mod-Q$. Under certain technical assumptions, which are satisfied if $\2A$ is
semisimple spherical, \cite[Theorem 3.3]{sch} then implies the braided monoidal
equivalence $\2Z(\2A)\stackrel{\otimes}{\simeq}_{br}\2Z(\2B)$ claimed in the abstract.
We hope to discuss these matters in more detail in \cite{mue15}.
\erem


\section{Weak Monoidal Morita Equivalence `$\approx$'}\label{s-Morita}
\bdefin \label{Morita0}
A Morita context is a bicategory $\2F$ satisfying 
\begin{enumerate}
\item $\displaystyle\obj\,\2F=\{\6A,\6B\}$.
\item Idempotent 2-morphisms in $\2F$ split.
\item There are mutually two-sided dual 1-morphisms $J: \6B\rarr\6A$,
$\oj: \6A\rarr\6B$ such that the compositions
$\eta_J\circ e_J\in\id_{\11_\6A}$ and $d_J\circ\ve_J\in\id_{\11_\6B}$ are invertible.
\end{enumerate}
\edefin

\bdefin \label{Morita}
Two (preadditive, k-linear) tensor categories $\2A, \2B$ are weakly
mo\-no\-id\-ally Morita equivalent, denoted $\2A\approx\2B$, iff there exists a
(preadditive, k-linear) Morita context $\2F$ such that 
$\2A\stackrel{\otimes}{\approxeq}\END_\2F(\6A)$ and 
$\2B\stackrel{\otimes}{\approxeq}\END_\2F(\6B)$. We recall that in the non-additive 
case this means that there are monoidal equivalences 
$\ol{\2A}^p\stackrel{\otimes}{\simeq}\END_\2F(\6A)$ and  
$\ol{\2B}^p\stackrel{\otimes}{\simeq}\END_\2F(\6B)$, whereas for preadditive and k-linear
categories we require $\ol{\2A}^{p\oplus}\stackrel{\otimes}{\simeq}\END_\2F(\6A)$ and  
$\ol{\2B}^{p\oplus}\stackrel{\otimes}{\simeq}\END_\2F(\6B)$. In this situation $\2F$ is
called a Morita context for $\2A, \2B$. 
\edefin

\brem 1. If in Definition \ref{Morita0} we admit $\6A=\6B$ the implication 
$\2A\stackrel{\otimes}{\approxeq}\2B\ \impl\ \2A\approx\2B$ is obvious. 

2. In \cite{par} a notion of Morita equivalence for module categories of Hopf
algebras was considered, which has some similarities with the ours. Furthermore, it was
shown that Hopf algebras with Morita equivalent module categories (in the sense of 
\cite{par}) have the same dimension. This is reminiscent of our Proposition \ref{eqofdims}. 

3. If the structure morphisms $v, w$ of the Frobenius algebra are isomorphisms with 
$v^{-1}=v', w^{-1}=w'$ it is easy to see that the functor $X\mapsto \oj XJ$ is faithful,
full, essentially surjective and monoidal. Thus 
$\ol{\2A}^{p}\stackrel{\otimes}{\simeq}\ol{\2B}^{p}$, viz.\ $\2A$ and $\2B$ are
{\it strongly} monoidally Morita equivalent $\2A\stackrel{\otimes}{\approxeq}\2B$. 

4. Additional restrictions on the Morita context will be required if $\2A, \2B$ are
spherical or $*$-categories.
\erem

By Lemma \ref{motiv} a Morita context $\2F$ for tensor categories $\2A, \2B$ provides us
with canonical Frobenius algebras $(Q=J\oj, \ldots)$ and $(\hat{Q}=\oj J, \ldots)$ in $\2A$
and $\2B$, respectively. Conversely, the construction of the preceding subsection provides
us with a means of constructing tensor categories which are weakly Morita equivalent to a
given one, together with a Morita context: 

\blemma Let $\2A$ be a strict tensor category, let $Q$ be a Frobenius algebra satisfying
(\ref{W4}) and let $\2E$ be as constructed in the preceding subsection.
Then $\2B=\HOM_\2E(\6B,\6B)$ is weakly Morita equivalent to $\2A$, a Morita context 
being given by $\2E$. 
\elemma 
\prf Obvious by Theorems \ref{main0} and Definition \ref{Morita}. \qed

\bprop Let $\2A\approx\2B$ with Morita context $\2F$. Let $(Q,v,v',w,w')$ be the
Frobenius algebra in $\2A$ arising as in Lemma \ref{motiv} and $\2E$ as in Theorem
\ref{main0}. Then there is an equivalence of bicategories $\2E\simeq\2F$. In particular,
we have an equivalence $\HOM_\2E(\6B,\6B)\simeq\HOM_\2F(\6B,\6B)$ of tensor categories
under which the canonical Frobenius algebras $\hat{Q}=\oj J$ of $\HOM_\2E(\6B,\6B)$ and 
$\HOM_\2F(\6B,\6B)$ go into each other. 
\label{main2}\eprop
\prf Obvious in view of Theorem \ref{t-uniq}. \qed 

In view of this result it is intuitively clear that given a tensor category $\2A$ there is
a `bijective  correspondence' between (1) tensor categories which are weakly Morita
equivalent to $\2A$ and (2) canonical Frobenius algebras in $\2A$. This relation certainly
deserves to be made precise.  

That weak Morita equivalence is in fact an equivalence relation is not entirely obvious.
\bprop 
Weak monoidal Morita equivalence is an equivalence relation.
\eprop
\prf Symmetry and reflexivity of the relation $\approx$ are obvious. Assume
$\2A\approx\2B$ and $\2B\approx\2C$ with respective Morita contexts $\2E_1, \2E_2$ whose
objects we call $\6A, \6B_1$ and $\6B_2, \6C_2$, respectively. In order to prove transitivity
we must find a Morita context for $\2A$ and $\2C$. Since the definition of weak monoidal
Morita equivalence involves only the subobject-completions, we may assume without
restriction of generality that $\2A,\2B,\2C$ have subobjects. We identify $\2A$ and
$\END_{\2E_1}(\6A)$. By definition of a Morita context we have 
\[ \END_{\2E_1}(\6B_1) \, \stackrel{\otimes}{\simeq} \, \2B \,
   \stackrel{\otimes}{\simeq}\END_{\2E_2}(\6B_2), \]
and replacing $\2E_1, \2E_2$ by equivalent bicategories we may assume  
$\END_{\2E_1}(\6B_1)\stackrel{\otimes}{\cong}\END_{\2E_2}(\6B_2)$. The bicategories 
$\2E_1, \2E_2$ come with 1-morphisms $J_1\in\Hom_{\2E_1}(\6B_1,\6A)$,
$J_2\in\Hom_{\2E_2}(\6C_2,\6B_2)$ and their two-sided duals. We thus have a
Frobenius algebra $\5Q_2=(J_2\oj_2, v_2, v_2', w_2, w_2')$ in
$\END_{\2E_2}(\6B_2)$. Using 
the isomorphism $\END_{\2E_1}(\6B_1)\stackrel{\otimes}{\cong}\END_{\2E_2}(\6B_2)$ we obtain
the Frobenius algebra $(\tilde{Q}_2, \tilde{v}_2, \tilde{v}_2', \tilde{w}_2, \tilde{w}_2')$
in $\End_{\2E_1}(\6B_1)$. We define $Q_3=J_1\tilde{Q}_2\ol{J}_1\in\End_{\2E_1}(\6A)$ and
claim that this is part of a Frobenius algebra $\5Q_3$ in $\2A=\END_{\2E_1}(\6A)$ if we
put  
\bean v_3 &=& \id_{J_1}\otimes \tilde{v}_2\otimes\id_{\ol{J}_1}\mcirc v_1, \\
 v_3' &=& v_1' \mcirc \id_{J_1}\otimes \tilde{v}'_2\otimes\id_{\ol{J}_1}, \\
 w_3 &=& \id_{J_1\tilde{Q}_2}\otimes e_{J_1}\otimes\id_{\tilde{Q}_2\ol{J}_1} \mcirc 
    \id_{J_1}\otimes \tilde{w}_2\otimes\id_{\ol{J}_1},\\
 w_3' &=&  \id_{J_1}\otimes \tilde{w}_2\otimes\id_{\ol{J}_1} \mcirc 
   \id_{J_1\tilde{Q}_2}\otimes \eta_{J_1}\otimes\id_{\tilde{Q}_2\ol{J}_1}.
\eean
The verification of (\ref{W1}-\ref{W3}) is straightforward and therefore omitted. 
(This is quite similar to \cite[Theorem IV.8.1]{cwm} on the composition of adjoints.)
Let $\2E_3$ be the bicategory obtained from $\2A$ and $\5Q_3$ and $\{\6A_3, \6C_3\}$ its
objects. We denote $\2C_3=\END_{\2E_3}(\6C_3)$. There are functors 
$F_1: X\mapsto\ol{J}_2\ol{J}_1XJ_1J_2$ from $\2A$ to $\END_{\2E_2}(\6C_1)$ (the
composition of $X\mapsto\ol{J}_1XJ_1, \2A\rarr\END_{\2E_1}(\6B_1)$ and 
$X\mapsto\ol{J}_2XJ_2, \END_{\2E_2}(\6B_2)\rarr\END_{\2E_2}(\6C_2)$)
and $F_2: X\mapsto\ol{J}_3XJ_3$ from $\2A$ to $\END_{\2E_3}(\6C_3)$.
In view of the definitions of $\2E_1, \2E_2, \2E_3$ and of $\5Q_3$ it is clear that the
images of $F_1, F_2$ as full subcategories of $\END_{\2E_2}(\6C_2)$ and 
$\END_{\2E_3}(\6C_3)$, respectively, are equivalent as tensor categories. Since the
tensor categories $\END_{\2E_2}(\6C_2)$ and $\END_{\2E_3}(\6C_3)$ are equivalent to the
subobject closures of the respective full subcategories they are themselves equivalent: 
$\END_{\2E_2}(\6C_2)\stackrel{\otimes}{\simeq}\END_{\2E_3}(\6C_3)$.
Together with $\2C\stackrel{\otimes}{\simeq}\END_{\2E_2}(\6C_2)$ and
$\END_{\2E_3}(\6C_3)\approx\2A$ this implies $\2C\approx\2A$.
\qed

\brem 1. Comparing our notion of weak Morita equivalence with the one for rings we see
that Definition \ref{Morita} is similar (but not quite, see below) to the property of two
rings $R,S$ of admitting an invertible $A-B$ bimodule. Now, it is known \cite{bass} that
this is the case iff one has either of the equivalences $R-\mod\simeq S-\mod$,
$\mod-R\simeq\mod-S$. Since there is a notion of representation bicategory of a tensor
category, cf.\ \cite{par0} and \cite[Chapter 4]{neuchl}, is is very natural to conjecture
that two tensor categories are weakly monoidally Morita equivalent iff their
representation categories are equivalent bicategories. (This would make the transitivity
of the relation $\approx$ obvious.) We hope to go into this question elsewhere. There is,
however, one caveat, viz.\ in Definition \ref{Morita} we do not require the 1-morphisms to be
mutually inverse (in the sense $J\oj \cong\11_\6A$, $\oj J\cong\11_\6B$) but only to be
adjoint (conjugate). Already as applied to rings this yields a weaker equivalence
relation. 

2. It is interesting to note that the usual Morita equivalence of (non-monoidal)
categories can be expressed via the existence of a pair of mutually inverse 1-morphisms in
the 2-category of small categories, distributors and their morphisms, cf.\ 
\cite[Section 7.9]{bor}. One might ask whether a useful generalization is obtained by
requiring only the existence of a two-sided adjoint pair of distributors.

3. Let $\2A$ be a $\End(\11)$-linear tensor category. By the definitions and results of the 
preceding and the present section, a canonical Frobenius algebra $\5Q$ in $\2A$ gives rise
to a tensor category $\2B\approx\2A$ together with a Morita context $\2E$. Conversely, a
Morita context $\2E$ for $\2B\approx\2A$ gives rise to a canonical Frobenius algebra $\5Q$
in $\2A$. It is clear that this correspondence can be formalized as a one-to-one
correspondence, modulo appropriate equivalence relations on both sides, between canonical
Frobenius algebras in $\2A$ and tensor categories $\2B\approx\2A$ together with a Morita
context. Here we do not pursue this further for lack of space.
\erem


\section{Linear Categories}\label{lincat}
\subsection{Linear categories}
From now on all categories are linear over a field $\7F$ and all hom-sets are finite
dimensional over $\7F$. If $\11$ is not simple, i.e.\ $\End(\11)\not\cong\7F$, we require
canonical Frobenius algebras to satisfy $\lambda_1,\lambda_2\in\7F^*$, not just
$\lambda_1,\lambda_2\in\End(\11)^*$.

\bprop \label{p-simple}
Let $\2A$ be $\7F$-linear with possibly non-simple unit. Let $\5Q$ be a canonical
Frobenius algebra. Then the following holds for the bicategory $\2E$ of Theorem \ref{main0}. 
\begin{description}
\item[(i)] $J$ is simple iff $\oj$ is simple iff $\dim\Hom_\2A(Q, \11)=1$ iff
$\dim\Hom_\2A(\11, Q)=1$.
\item[(ii)] $\11_\6B$ is simple iff 
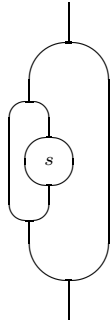
\begin{equation}
\begin{tangle}
\hstep\cd\\
\hh\hcd\step[1.5]\id\\
\id\step\O s\step[1.5]\id\\
\hh\hcu\step[1.5]\id\\
\hstep\cu
\end{tangle}
\label{sandwich}\end{equation}
is a multiple $F(s)$ of $\id_Q$ for every $s\in\End(Q)$. 
\end{description}
Furthermore, (i) implies simplicity of $\11_\2A=\11_\6A$ and $\11_\6B$.
\eprop
\prf $J$ is simple iff $\End_\2E(J)\cong\7F$. By definition of $\2E$ this is the case iff 
$\dim\Hom_\2A(Q,\11)=1$. Similarly, $\oj$ is simple iff $\dim\Hom_\2A(\11,Q)=1$. The
remaining equivalence in (i) follows from
$\Hom(J,J)\cong\Hom(J\oj,\11)\cong\Hom(\oj,\oj)$, which is a trivial consequence of
duality of $J, \oj$ in $\2E$. By Corollary \ref{coro1} all functors $J\otimes -$, 
$\oj\otimes -$, $-\otimes J$, $-\otimes\oj$ are faithful. Thus both $\End(\11_\6A)$ and
$\End(\11_\6B)$ embed as subalgebras into $\End(J)$. Thus $\End(J)\cong\7F$ implies   
$\End(\11_\6A)\cong\7F$ and $\End(\11_\6B)\cong\7F$.
\qed

\brem 1. The implication `$J$ simple $\Leftrightarrow \oj $ simple $\impl\ \11_\6B$ 
simple' is reminiscent of the situation for an inclusion $B\subset A$ of von Neumann
algebras, where we trivially have 
\[ A\cap B'=\7C\11 \quad\Leftrightarrow\quad B'\cap (A')'=\7C\11 
  \quad \impl\quad B'\cap B=\7C\11 \ \mbox{and}\ A'\cap A=\7C\11.
\]

2. Note that we did not need a duality between $\Hom(X,Y)$ and $\Hom(Y,X)$ as it exists in
$*$-categories and non-degenerate spherical categories in order to conclude
$\dim\Hom(\11,Q)=\dim\Hom(Q,\11)$.
\erem

\bdefin Let $\2A$ be a $\7F$-linear tensor category. A canonical Frobenius algebra $\5Q$
in $\2C$ is irreducible if the equivalent conditions of (i) above hold. (By the above this
is possible only if $\11_\2A$ is simple.)
\edefin


\subsection{$*$-Categories}
We recall that a 1-morphism $\ol{X}$ of a $*$-2-category (or object in a 
$\otimes$-$*$-category) is conjugate to $X$ in the sense of \cite{lro} if there are  
$r_X: \11\rarr X\ol{X}, \ol{r}_X: \11\rarr\ol{X}X$ satisfying (\ref{conj-eq0}). In the
generic notation of Subsection \ref{ss-duality} this amounts to 
$e_X=\ol{r}_X$, $\ve_X=r_X$, $d_X=r^*_X$, $\eta_X=\ol{r}^*_X$, thus 
$d_X=(\ve_X)^*, \eta_X=(e_X)^*$. This implies that the Frobenius algebra
$\5Q=(\ol{X}X,\ldots)$ of Lemma \ref{motiv} satisfies the conditions $v'=v^*, w'=w^*$
which we mentioned in Remark \ref{frob-rem}. Therefore the canonicity conditions
(\ref{W4}-\ref{W5}) amount to saying that $v, w$ are (non-zero) multiples of isometries
thus up to renormalization $(Q,v,w)$ is an `abstract Q-system' in the sense of
\cite{lro}. In \cite{lro} is is shown that, quite remarkably, in this situation (\ref{W3})
holds automatically. Therefore a canonical Frobenius algebra in a $*$-category is
the same as an algebra $(Q, v, w^*)$ where $v,w$ are multiples of isometries. 

\bdefin A Q-system is a canonical Frobenius algebra in a $*$-category satisfying 
$v'=v^*, w'=w^*$. It is normalized if the Frobenius algebra $(Q,v,v^*,w,w^*$) is
normalized, i.e.\ if $v^*\circ v=w^*\circ w$.
\edefin

\bprop Let $\2A$ be a tensor $*$-category and $(Q,v,v'=v^*,w,w'=w^*)$ a Q-system in
$\2A$. Then $\2E_0$ has a positive $*$-operation $\#$ which extends the given one on
$\2A$. Let $\2E_*$ be the full sub-bicategory of $\2E$ whose 1-morphisms are $(X, p)$
where $p=p\bullet p=p^\#$. Then $\2E_*$ is equivalent to $\2E$ and has a positive
$*$-operation $\#$. 
\eprop
\prf Since the 2-morphisms of $\2E$ are given in terms of 1-morphisms in $\2A$, we must 
denote the $*$-operation of $\2E_0$ by $\#$ in order to avoid confusion. For morphisms in
$\END_\2E(\6A)\equiv\2A$ we obviously define $s^\#=s^*$. Let $X,Y\in\obj\2A$
throughout. For $s\in\Hom_\2E(XJ,YJ)\equiv\Hom_\2A(XQ,Y)$ we define
\[ s^{\#}= \id_X\otimes r^*\mcirc s^*\otimes\id_Q 
   \quad\in\Hom_\2A(YQ,X)\equiv\Hom_\2E(YJ,XJ),  \]
for $s\in\Hom_\2E(\oj X,\oj Y)\equiv\Hom_\2A(X,QY)$ we posit
\[ s^{\#}= \id_Q\otimes s^* \mcirc r\otimes\id_Y 
   \quad\in\Hom_\2A(Y,QX)\equiv\Hom_\2E(\oj Y,\oj X) \] 
and or $s\in\Hom_\2E(\oj XJ,\oj YJ)\equiv\Hom_\2A(XQ,QY)$ we put
\[ s^{\#}= \id_{QX}\otimes r^* \mcirc \id_Q\otimes s^*\otimes\id_Q \mcirc r\otimes\id_{YQ}
   \quad\in\Hom_\2A(YQ,QX)\equiv\Hom_\2E(\oj YJ,\oj XJ). \]
(Recall that $r=w\circ v$.) Graphically the last definition looks like
\[ 
{\left( \quad
\begin{tangle}
\object{Q}\step\object{Y}\\
\hh\id\step\id\\
\hh\frabox{s}\\
\hh\id\step\id\\
\object{X}\step\object{Q}
\end{tangle}
\quad \right)}^\# = \quad
\begin{tangle}
\object{Q}\step\object{X}\\
\hh\id\step\id\step\hcoev\\
\hh\id\step\frabox{s^*}\step\id\\
\hh\hev\step\id\step\id\\
\Step\object{Y}\step\object{Q}
\end{tangle}
\]
Antilinearity of these operations is obvious and involutivity follows from the duality
equation satisfied by $r, r^*$. The easy verification of contravariance 
($(s\circ t)^{\#}=t^{\#}\circ s^{\#}$) and monoidality 
($(s\times t)^{\#}=s^{\#}\times t^{\#}$) are left to the reader. We limit ourselves to
showing that $\#$ is positive. We consider only the case of morphisms between
$\6B-\6B$-morphisms, the others being similar. With 
$s\in\Hom_\2E(\oj XJ,\oj YJ)\equiv\Hom_\2A(XQ,QY)$ we compute 
\[ s^{\#}\bullet s= \quad
\begin{tangle}
\hstep\object{Q}\step[1.5]\object{X}\\
\hcd\step\id\step\hcoev\\
\hh\id\step\id\step\frabox{s^*}\step\id\\
\d\hev\dd\step\id\\
\hh\step\frabox{s}\Step\id\\
\step\id\step\cu\\
\step\object{X}\Step\object{Q}
\end{tangle}
\quad =\quad
\begin{tangle}
\object{Q}\step\object{X}\\
\hh\id\step\id\step\hcoev\\
\hh\id\step\frabox{s^*}\step\id\\
\hdcu\step\mobj{Y}\id\step\id\\
\hh\step\frabox{s}\step\id\\
\step\id\step\hcu\\
\step\object{X}\step[1.5]\object{Q}
\end{tangle}
\]
If this vanishes then (by sandwiching between $v'\otimes\id_X$ and $\id_X\otimes v$) also 
\[
\begin{tangle}
\step[1.5]\object{r^*}\\
\hh\id\step\hcoev\\
\hh\frabox{s^*}\step\id\\
\hh\id\step\id\step\id\\
\hh\frabox{s}\step\id\\
\hh\id\step\hev\\
\step[1.5]\object{r}
\end{tangle}
\quad\quad\mbox{and}\quad\quad
\begin{tangle}
\hh\id\step\id\step\id\\
\hh\frabox{s}\step\id\\
\hh\id\step\hev\\
\end{tangle}
\]
vanish, the latter by positivity of the $*$-operation in $\2A$. Now duality implies $s=0$,
thus $\#$ is a positive $*$-operation. 

We now turn to the bicategory $\2E$. Let $(X,p), (Y,q)$ be parallel 1-morphisms and 
$s: X\rarr Y$. By definition, $s$ is a morphism $(X, p)\rarr (Y,q)$ iff 
$s=s\bullet p=q\bullet s$, which is equivalent to $s^\#=p^\#\bullet s^\#=s^\#\bullet q^\#$.
Thus $s^\#: Y\rarr X$ is in fact in $\Hom_\2E((Y,q^\#),(X,p^\#))$. In the full sub-bicategory
$\2E_*$ we have $p^\#=p, q^\#=q$, thus $s^\#\in\Hom_\2E((Y,q),(X,p))$ as it should.
Finally, in a finite dimensional $*$-algebra (like $\End_\2E(X)$) every projection is
similar to an orthogonal projection. Thus every $(X,p)$ is isomorphic to $(X,q)$ where $q$
is an orthogonal projection. This proves $\2E_*\simeq\2E$.
\qed

\brem \label{semi1}
1. Let $\2A$ be a $*$-category which has subobjects and finite dimensional
hom-sets. Then positivity of the $*$-operation implies $\End(X)$ to be a multi matrix 
algebra for every $X$ and therefore semisimplicity of $\2A$. Since $\2E_*$ by construction
has retracts of 1-morphisms, finite dimensional hom-sets and a positive $*$-operation, we
conclude that $\2E_*$ is semisimple (in the sense that all categories $\HOM_\2E(?,!)$ are
semisimple). 

2. As explained in Section \ref{stars}, the notion of (two-sided) duals in $*$-categories
is local in that it does not necessitate a conjugation {\it map} (or functor)
$X\rarr\ol{X}$ together with {\it chosen} morphisms $\11\rarr X\otimes\ol{X}$ etc. If
$X\in\2A$ has a conjugate $\ol{X}$ it is easy to see that $XJ: \6B\rarr\6A$ has
$\oj\ol{X}$ as conjugate etc. Thus if $\2A$ has conjugates for all objects then $\2E_*$
has conjugates for all 1-morphisms. Therefore the above construction of a $*$-structure on
$\2E$ completes the discussion of $*$-categories.

3. A self-conjugate object $X$ in a $*$-category is called real (or orthogonal) if there
exists a solution $(X,r_X,\ol{r}_X)$ of the conjugate equations where $r_X=\ol{r}_X$ and
pseudo-real (or symplectic) if we can put $r_X=-\ol{r}_X$. (Every simple self-conjugate
object is either real or pseudo-real, cf.\ \cite{lro}.) As already observed in \cite{lro}
the object of a `Q-system' is real since $(Q,r,\ol{r})$ with $\ol{r}=r=w\circ v$ is a
solution of the conjugate equations. By minimality of the intrinsic dimension $d(Q)$, this
solution of the conjugate equations is standard iff $r^*\circ r=v^*\circ w^*\circ w\circ v$
equals $d(Q)\id_\11$. This is automatic when $(Q,v,w)$ is irreducible, i.e.\
$\dim\Hom(\11,Q)=1$, as was shown in \cite{lro} using the construction of a subfactor from
a Q-system. Our construction of the bicategory $\2E$ allows to give a simple purely
categorical argument. If the Frobenius algebra $\5Q$ is irreducible then Proposition
\ref{p-simple} implies irreducibility of $J: \6B\rarr\6A$ and $\oj: \6B\rarr\6A$ in $\2E$.
Then $\Hom_\2E(\11_\6A,J\oj)$ and $\Hom_\2E(\11_\6B,\oj J)$ are one dimensional, which
implies $v^*\circ v=d(J)\id_\11$ and $w^*\circ w=d(J)\id_Q$. Thus 
$v^*\circ w^*\circ w\circ v=d(J)^2\id_\11=d(Q)\id_\11$ and $(Q,r,\ol{r})$ is standard.
\erem


\subsection{Spherical categories}
Just as Frobenius algebras in $*$-categories we required the compatibility condition
$v'=v^*, w'=w^*$, we need a compatibility of $(Q,v,v',w,w')$ with the spherical structure
of $\2A$. Let $Q$ be a 1-morphism in a strict spherical 2-category $\2F$. Then $Q=J\oj$ is
strictly selfdual: $\ol{Q}=\ol{J\oj}=J\oj=Q$. If we consider the Frobenius algebra
obtained from Lemma \ref{motiv} with $e_J=\ve(J)$, $\ve_J=\ve(\oj)$, $d_J=\ol{\ve}(\oj)$,
$\eta_J=\ol{\ve}(J)$ then obviously 
\[ w\circ v=\id_J\otimes \ve(\oj)\otimes\id_\oj \mcirc \ve(J)=\ve(J\oj)=\ve(Q) \]
and $v'\circ w'=\ol{\ve}(Q)$. Conversely, we have the following

\blemma \label{lemma-piv}
Let $\2A$ be a strict pivotal $\7F$-linear tensor category. Let $(Q, \ldots)$ be a 
canonical Frobenius algebra such that $\ol{Q}=Q$ and $w\circ v=\ve(Q)$, 
$v'\circ w'=\ol{\ve}(Q)$. Then the bicategory $\2E$ of Theorem \ref{main0} has a strict 
pivotal structure which restricts to the one of $\2A$. 
\elemma
\prf We extend the conjugation map to $\2E_0$ as follows:
\[ \ol{``XJ"} := ``\oj \,\ol{X}", \quad \ol{``\oj X"}:=``\ol{X}J", \quad 
   \ol{``\oj XJ"}:=``\oj \,\ol{X}J". \]
Thus $\ol{\ol{``XJ"}}=``\ol{\ol{X}}J"=``XJ"$ etc., and the conditions (\ref{conjmap}) are
obvious consequences of those for $\2A$. Using the notations of Theorem \ref{main0} we
define $\ve(J)=e_J, \ve(\oj)=\ve_J, \ol{\ve}(J)=\eta_J, \ol{\ve}(\oj)=d_J$. Since we know
$\ve(X)$ for $X\in\2A$, condition (2) of Definition \ref{pivotal} enforces
$\ve(XJ)=\id_X\otimes\ve(J)\otimes\id_{\ol{X}} \circ\ve(X)$ and analogously for
$\ve(\oj X), \ve(\oj XJ)$. With this definition conditions (i) and (ii) of Definition
\ref{pivotal} are clearly satisfied. (The conditions on $w\circ v, v'\circ w'$ are 
necessary and sufficient for $\ve, \ol{\ve}$ being well defined for all 1-morphisms since
they guarantee $\ve(J\oj)=\ve(Q)$.) It remains to verify (iii). Let, e.g.,
$s\in\Hom_\2E(XJ,YJ)$, represented by $\tilde{s}\in\Hom_\2A(XQ,Y)$. Now an easy
computation shows that
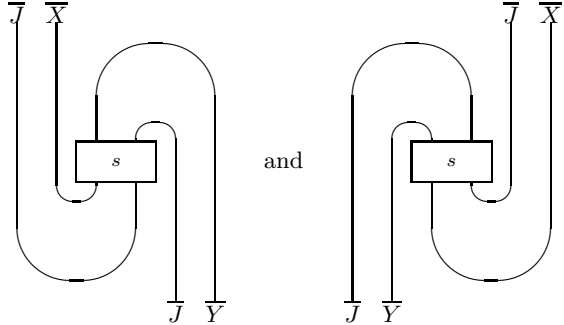
\[
\begin{tangle}
\object{\oj}\step\object{\ol{X}}\\
\id\step\id\step\mcoev\\
\hh\id\step\id\step\id\step\hcoev\step\id\\
\hh\id\step\id\step\frabox{s}\step\id\step\id\\
\hh\id\step\hev\step\id\step\id\step\id\\
\mev\step\id\step\id\\
\step[4]\object{\oj}\step\object{\ol{Y}}
\end{tangle}
\quad\quad \mbox{and} \quad\quad
\begin{tangle}
\step[4]\object{\oj}\step\object{\ol{X}}\\
\mcoev\step\id\step\id\\
\hh\id\step\hcoev\step\id\step\id\step\id\\
\hh\id\step\id\step\frabox{s}\step\id\step\id\\
\hh\id\step\id\step\id\step\hev\step\id\\
\id\step\id\step\mev\\
\object{\oj}\step\object{\ol{Y}}
\end{tangle}
\]
are represented in $\2A$ by
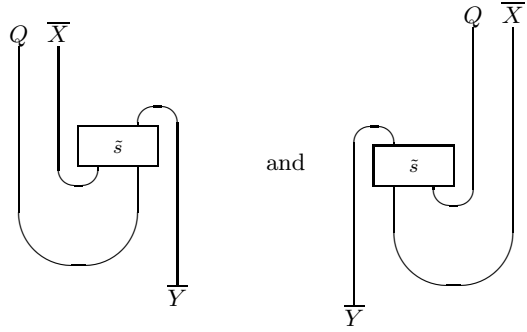
\[
\begin{tangle}
\object{Q}\step\object{\ol{X}}\\
\id\step\id\Step\hcoev\step\\
\hh\id\step\id\step\frabox{\tilde{s}}\step\id\\
\hh\id\step\hev\step\id\step\id\\
\mev\step\id\\
\step[4]\object{\ol{Y}}
\end{tangle}
\quad\quad \mbox{and} \quad\quad
\begin{tangle}
\step[3]\object{Q}\step\object{\ol{X}}\\
\step[3]\id\step\id\\
\hh\hcoev\Step\id\step\id\\
\hh\id\step\frabox{\tilde{s}}\step\id\step\id\\
\hh\id\step\id\step\hev\step\id\\
\id\step\mev\\
\object{\ol{Y}}
\end{tangle}
\]
respectively. These two expressions coincide since $\2A$ is pivotal, showing that the dual
of the 2-morphism $s\in\Hom_{\2E_0}(XJ,YJ)\equiv\Hom_\2A(XQ,Y)$ is precisely given by 
$\ol{s}\in\Hom_\2A(\ol{Y},Q\ol{X})\equiv\Hom_{\2E_0}(\oj Y,\oj X)$, where $\ol{s}$ is
computed in $\2A$. The same holds for the other types of 2-morphisms. 
The completion w.r.t.\ subobjects in a spherical tensor category
behaves nicely w.r.t.\ duality ($\ol{(X,p)}:=(\ol{X},\ol{p})$), and the same holds for
1-morphisms in a 2-category. Thus $\2E$ has strict duals and is pivotal. 
\qed

As explained in Subsection \ref{ss-spherical}, the notion of sphericity of a (non-monoidal)
2-category seems to make sense only if all identity 1-morphisms are simple, which is why
we need more stringent conditions in the following

\bprop \label{prop-sph}
Let $\2A$ be a strict pivotal $\7F$-linear tensor category with simple unit. Let
$(Q, \ldots)$ be a canonical Frobenius algebra such that 
\begin{itemize}
\item[(i)] $\displaystyle\ol{Q}=Q$.
\item[(ii)] $\displaystyle w\circ v=\ve(Q)$ and $\displaystyle v'\circ w'=\ol{\ve}(Q)$.
\item[(iii)] The Frobenius algebra satisfies the equivalent conditions (ii) of Proposition
\ref{p-simple}. 
\end{itemize}
Then $\2E$ is spherical iff $\2A$ is spherical and $\5Q$ is normalized. If furthermore 
the trace on $\2A$ is non-degenerate then also its natural extension to $\2E$ is
non-degenerate.
\eprop
\prf By the preceding result and  Proposition \ref{p-simple}, $\2E$ is a bicategory with strict
pivotal structure and simple $\11_\6A, \11_\6B$. Since $\2A$ sits in $\2E$ as the corner
$\END_\2E(\6A)$, sphericity of $\2A$ is clearly necessary for $\2E$ to be
spherical. Condition (ii) implies 
\[ d(Q)=\ol{\ve}(Q)\circ\ve(Q)=v'\circ w'\circ w\circ v=\lambda_1\lambda_2. \]
Next we observe that the $F(s)$ in Proposition \ref{p-simple} is given by
\begin{equation} \label{Fs} F(s)=\frac{\lambda_1}{\lambda_2} tr_Q(s), \end{equation}
where $tr_Q(s)$ is the trace in $\2A$ of $s\in\End(Q)$. (Using (ii), sphericity of $\2A$
and $w'\circ w=\lambda_1\id_Q$ the trace (in $\End_\2A(Q)$) of (\ref{sandwich}) is seen to 
equal $\lambda_1^2 tr_Q(s)$. Since this trace is also equal to
$F(s)tr_Q(\id_Q)=F(s)d(Q)=F(s)\lambda_1\lambda_2$, the claimed equation follows.)

In order to check sphericity of $\2E$ we need to consider the traces on $\End_\2E(X)$
where $X$ is an $\6A-\6B$-, $\6B-\6A$- or $\6B-\6B$-morphism. 
Let $s\in\End_\2E(XJ,XJ)$ represented by $\tilde{s}\in\Hom_\2A(XQ,X)$. In view of our
definition of $\2E$ we have
\[
\begin{tangle}
\object{X}\\
\hh\id\\
\frabox{\tilde{s}}\\
\hh\id\step\id\\
\object{X}\step\object{Q}
\end{tangle}
\quad=\quad
\begin{tangle}
\object{X}\step[1.3]\object{\ol{\ve}(J)}\\
\hh\id\step\hcoev\\
\hh\frabox{s}\step\id\\
\hh\id\step\id\step\id\\
\object{X}\step\object{J}\step\object{\oj}
\end{tangle}
\]
and with $\ve(J)=e_J=v$ we have
\[
tr_L(s)=\quad
\begin{tangle}
\step[1.5]\object{\ol{\ve}(X)}\\
\mcoev\step[-2]\mobj{\ol{\ve}(J)}\\
\hh\id\step\hcoev\step\id\\
\hh\frabox{s}\step\id\obj{\oj}\step\id\obj{\ol{X}}\\
\hh\id\step\hev\step\id\\
\mev\\
\step[1.5]\object{\ve(X)}\\
\end{tangle}
\quad=\quad
\begin{tangle}
\step\object{\ol{\ve}(X)}\\
\coev\\
\hh\frabox{\tilde{s}}\step\id\\
\id\step\counit\step[.2]\mobj{v}\step[.8]\id\obj{\ol{X}}\\
\ev\\
\step\object{\ve(X)}\\
\end{tangle}
\]
which expresses $tr_L(s)$ in terms of a formula in $\2A$. The computation of
\[ tr_R(s)=\quad
\begin{tangle}
\step[1.5]\object{\ol{\ve}(\oj)}\\
\mcoev\step[-2]\mobj{\ol{\ve}(\ol{X})}\\
\hh\id\step\hcoev\step\id\\
\hh\id\step\id\step\frabox{s}\\
\hh\id\step\hev\step\id\\
\mev\\
\step[1.5]\object{\ve(\oj)}\\
\end{tangle}
\]
is slightly more involved. In order to simplify matters we pretend that $X=\11$ which
eliminates the inner trace over $X$. Thus
\[ tr_R(s)=\quad
\begin{tangle}
\step\object{\ol{\ve}(\oj)}\\
\coev\\
\hh\id\step\frabox{s}\\
\ev\\
\step\object{\ve(\oj)}\\
\end{tangle}
\quad=\quad
\begin{tangle}
\step[1.5]\object{\ol{\ve}(\oj)}\\
\mcoev\\
\hh\id\step\frabox{\tilde{s}}\step\id\\
\hh\id\step\id\step\id\step\id\\
\hh\hev\step\hev\\
\hstep\object{\ve(\oj)}\Step\object{\ve(\oj)}
\end{tangle}
\]
Now, $tr_R(s)\in\End_\2E(\11_\6B)\subset\End_\2E(\oj J)\equiv\End_\2A(Q)$ is represented by
\[ \lambda_1^{-1}\quad
\begin{tangle}
\hstep\cd\\
\hh\hcd\step[1.5]\id\\
\id\step\counit\step[.2]\mobj{v}\step[1.3]\id\\
\hh\id\step\mfrabox{\tilde{s}}\hstep\id\\
\hh\hcu\step[1.5]\id\\
\hstep\cu
\end{tangle}
\]
which by assumption (iii) and (\ref{Fs}) equals 
\[ \lambda_1^{-1}\frac{\lambda_1}{\lambda_2}tr_Q(v\circ\tilde{s})\id_Q=
  \lambda_1^{-1}\frac{\lambda_1}{\lambda_2}(\tilde{s}\circ v)\,\id_Q=
   \frac{\lambda_1}{\lambda_2}(\tilde{s}\circ v)\,\id_{\11_\6B}. \]
Thus, reintroducing the trace over $X$ we have
\[ tr_R(s)=\frac{\lambda_1}{\lambda_2}\quad
\begin{tangle}
\hh\hcoev\\
\hh\id\step\frabox{\tilde{s}}\\
\hev\step\counit\step[.2]\mobj{v}
\end{tangle}
\]
which, using the sphericity of $\2A$ coincides with $tr_L(s)$ iff $\5Q$ is normalized.
Completely analogous considerations hold for the traces on $\End(\oj X)$ and $\End(\oj XJ)$.

Assume now, that the trace on $\2A$ is non-degenerate and let
$s\in\Hom_{\2E_0}(\oj X,\oj Y)\equiv\Hom_\2A(X,QY)$. By non-degeneracy of the trace on
$\2A$ there is $u\in\Hom_\2A(QY,X)$ such that $tr_X(u\circ s)\ne 0$. Now defining
\[ t= \ 
\begin{tangle}
\object{Q}\step[1.5]\object{X}\\
\hh\id\step[1.5]\id\\
\hh\id\step\frabox{u}\\
\hh\ev\step\id\\
\Step\object{Y}
\end{tangle}
\quad \in\Hom_\2A(Y,QX)\equiv\Hom_{\2E_0}(\oj Y,\oj X), \]
one easily verifies $tr_{\oj X}(t\bullet s)=tr_X(u\circ s)\ne 0$. Thus the pairing
$\Hom_{\2E_0}(\oj X,\oj Y)\times\Hom_{\2E_0}(\oj Y,\oj X)\rarr\7F$ provided by the
trace is non-degenerate. The other cases are verified similarly. 

That the trace on $\2E_0$ extends to a non-degenerate trace on the completion
$\2E=\ol{\2E_0}^p$ follows from the following simple argument. Let $tr$ be a
non-degenerate trace, $e,f $ idempotents and $exf\ne 0$. Then there is $y$ such that
$tr(exfy)\ne 0$. By cyclicity of the trace $tr(exfy)=tr((exf)(fye))$, thus $y$ can chosen
such that $y=fye$. \qed 

\brem 1. If $\7F$ is quadratically closed every canonical Frobenius algebra can be turned
into a normalized one by renormalization. The sign of $d(J)$ depends on the choice of the
renormalization, but in the $*$-case one can achieve $d(J)=+\sqrt{d(Q)}>0$.

2. Every simple self-dual object in a spherical category is either orthogonal or
symplectic, cf.\ \cite[p.\ 4018]{bw1}. A simple object in a spherical category is
orthogonal iff we can obtain $\ol{X}=X$ in a suitable strictification of the category.
In the latter sense the object of a Frobenius algebra is orthogonal.
\erem

As mentioned in Remark \ref{semi1}, $*$-categories are automatically semisimple and
therefore semisimplicity of $\2A$ entails semisimplicity of $\2E_*$. In order to prove an
analogous result for $\2A$ semisimple spherical we need the following facts, which we
include since we are not aware of a convenient reference.

A trace on a finite dimensional $\7F$-algebra $A$ is a linear map $A\rarr\7F$ such that
$tr(ab)=tr(ba)$. It is non-degenerate if for every $a\ne 0$ there is $b$ such that
$tr(ab)\ne 0$. 

\blemma Let $A$ be a finite dimensional $\7F$-algebra and $tr: A\rarr \7F$ a
non-degenerate trace. If $\,tr$ vanishes on nilpotent elements then $A$ is
semisimple. Conversely, every trace (not necessarily non-degenerate) on a semisimple
algebra vanishes on nilpotent elements. \label{semi}\elemma 
\prf Let $R$ be the radical and $0\ne x\in R$. By non-degeneracy there is $y\in A$ such
that $tr(xy)\ne 0$. On the other hand $xy\in R$ and $tr(xy)=0$ since $R$ is
nilpotent. Thus $R=\{0\}$.  
As to the second statement, observe that every trace on a matrix algebra coincides up to a
normalization with the usual trace. (Thus $tr(e_{i,j})=\alpha\delta_{i,j}$ with
$\alpha\in\7F$.) The latter vanishes for nilpotent matrices. The trace of a multi matrix
algebra is just a linear combination of such matrix traces on the simple  subalgebras, and
for general semisimple algebras the result follows by passing to an algebraic closure of
$\7F$.  
\qed

\bprop \label{th-semisim} Let $\7F$ be algebraically closed and $\2A$ strict spherical
$\7F$-linear and semi\-sim\-ple. Let $\5Q$ as in Proposition \ref{prop-sph}, including
$\lambda_1=\lambda_2$. Then $\2E$ is spherical and semisimple.
\eprop
\prf By construction, idempotent 2-morphisms in $\2E$ split. Since $\2A$ is semisimple,
thus has direct sums, the same holds for the 1-morphisms in $\2E$. In order to prove that 
$\2E$ it remains to show that $\End_\2E(X)$ is a multi matrix algebra for every 1-morphism
$X$. 

We begin by proving that $\End_{\2E_0}(XJ)=\End_{\2E}(XJ)$ is semisimple.
By Proposition \ref{lemma-piv} the trace $tr_{XJ}$ on $\End_{\2E_0}(XJ)$ is non-degenerate.
By Corollary \ref{coro1} the algebra homomorphism
\[ -\otimes\id_{\oj }:\ \End_{\2E_0}(XJ)\rarr \End_{\2E_0}(XJ\oj )=\End_\2A(XQ) \]
is injective, such that we can consider $\End_{\2E_0}(XJ)$ as a subalgebra of
$\End_{\2E_0}(XJ\oj )$. Furthermore, by sphericity of $\2E$ we have for
$s\in\End_{\2E_0}(XJ)$ 
\[ tr_{XJ\oj }(s\otimes\id_{\oj })=d(J)\,tr_{XJ}(s). \]
If $s$ is nilpotent then also $s\otimes\id_{\oj }\in\End_{\2E_0}(XQ)$ is nilpotent. 
Thus $tr_{XQ}(s\otimes\id_{\oj })=0$ by Lemma \ref{semi} and thus $tr_{XJ}(s)=0$ since
$d(J)\ne 0$. Therefore $\End_{\2E_0}(XJ)$ is semisimple by Lemma \ref{semi} and a multi
matrix algebra by algebraic closedness of $\7F$. If $A$ is a matrix algebra and $p=p^2\in
A$ then also $pAp\subset A$ is a matrix algebra. Thus also the endomorphism algebras
$\End_\2E((XJ, p))$ in the completion $\2E=\ol{\2E_0}^p$ are multi matrix
algebras. Perfectly similar arguments apply to $\End_{\2E}((\oj X,p))$ and 
$\End_{\2E}((\oj XJ,p))$ for all $X$. 
\qed

The conditions (i) and (ii) in Proposition \ref{prop-sph} on the Frobenius algebra are fairly
rigid and probably not satisfied in many applications. Furthermore, the above results 
should be generalized to the situation where neither the tensor product nor the duality of
$\2A$ are strict. In the following result we limit ourselves to the degree of generality
which will be needed for the application in \cite{mue10}. It is fairly obvious that also
the strictness conditions on $\2A$ can be dropped by inserting the appropriate
isomorphisms wherever needed.

\btheor \label{t-spherical2}
Let $\7F$ be algebraically closed and $\2A$ be $\7F$-linear, strict monoidal, strict
spherical and semisimple. Let $\5Q=(Q,v,v',w,w')$ be a normalized canonical Frobenius
algebra in $\2A$ satisfying condition (ii) of Proposition \ref{p-simple}.
Then the bicategory $\2E$ of Theorem \ref{main0} has simple $\6B$-unit $\11_\6B$, is
semisimple and has a (non-strict) spherical structure extending that of  
$\2A\stackrel{\otimes}{\simeq}\END_\2E(\6A)$.
\etheor
\prf The $\7F$-linear bicategory $\2E$ is defined as in Theorem \ref{main0}. We define a
conjugation map on the 1-morphisms as in Lemma.\ \ref{lemma-piv}. Thus we still have
$\ol{\ol{X}}=X$ for all 1-morphisms. By Lemma \ref{lemma-selfd}, $Q$ is self-dual, and
since duals are unique up to isomorphism the conjugation map $X\mapsto\ol{X}$ of $\2A$
satisfies $\ol{Q}\cong Q$. In fact, there is a unique isomorphism $s: Q\rarr\ol{Q}$ such
that  
\begin{equation} \label{bla}
 \id_Q\otimes s\mcirc r=\id_Q\otimes s\mcirc w\mcirc v=\ve(Q): \ \11\rarr Q\otimes\ol{Q}. 
\end{equation}
This is seen as follows. If $s:Q\rarr\ol{Q}$ satisfies (\ref{bla}) then
\[
\begin{tangle}
\Step\object{\ol{Q}}\\
\hh\Step\id\\
\hh\Step\O s\\
\hh\step[-.2]\mobj{r'}\step[.2]\coev\step\id\\
\hh\id\step\ev\step[.2]\mobj{r}\\
\object{Q}
\end{tangle}
\quad = \quad\quad
\begin{tangle}
\object{\ol{Q}}\\
\hh\id\\
\hh\O s\\
\hh\id\\
\object{Q}
\end{tangle}
\quad\quad = \quad
\begin{tangle}
\Step\object{\ol{Q}}\\
\hh\Step\id\\
\hh\step[-.3]\mobj{r'}\step[.3]\coev\step\id\\
\hh\id\step\ev\step[.2]\mobj{\ve(Q)}\\
\object{Q}
\end{tangle}
\]
On the other hand, it is equally easy to see that $s: Q\rarr\ol{Q}$ as defined by the
second half of this equation satisfies (\ref{bla}). In view of $Q=J\circ\oj$ (which is
true by construction of $\2E$) we have
$\ol{``XJ"\circ``\oj Y"}=\ol{XQY}=\ol{Y}\,\ol{Q}\,\ol{X}$. This coincides with
$\ol{``\oj Y"}\circ\ol{``XJ"}=``\ol{Y}J"\circ``\oj\,\ol{X}"=\ol{Y}Q\ol{X}$ only
if $Q=\ol{Q}$, which we do not assume. In any case there is an isomorphism
\[ \gamma_{``XJ",``\oj Y"}=\id_{\ol{Y}}\otimes s\otimes\id_{\ol{X}}: \ 
   \ol{``\oj Y"}\circ\ol{``XJ"} \rarr \ol{``XJ"\circ``\oj Y"} \]
and similar for all other pairs of composable 1-morphisms. This makes $\2E_0$ and $\2E$
bicategories with dual data in the sense of an obvious generalization of
\cite{bw2}. The definition of $\ve(J), \ol{\ve}(J), \ve(\oj), \ol{\ve}(\oj)$ and therefore
of $\ve$ and $\ol{\ve}$ for all 1-morphisms is as in Proposition \ref{prop-sph}. Yet the
verification of the conditions in Definition \ref{pivotal} is slightly more involved since we
must insert appropriate isomorphisms in the lower lines of the commutative diagrams in
condition (2). To illustrate this we consider the diagram
\[
\begin{diagram}
\11_\6A & \rTo^{\ve(XJ)} & XJ\otimes\ol{XJ} \\
\dTo^{\ve(XJ\otimes\oj Y)} & & \dTo_{\id_{XJ}\otimes\ve(\oj Y)\otimes\id_{\ol{XJ}}} \\
XJ\otimes\oj Y\otimes \ol{XJ\otimes\oj Y} & \rTo_{\id_{XJ\otimes\oj Y}\otimes\gamma^{-1}_{XJ,\oj Y}} & XJ\otimes\oj Y\otimes\ol{\oj Y}\otimes \ol{XJ} 
\end{diagram}
\]
In terms of the category $\2A$, where all this ultimately takes place, this is
\[
\begin{diagram}[midshaft]
\11 & \rTo^{\ve(X)} & X\ol{X} & \rTo^{\id_X\otimes v\otimes\id_{\ol{X}}} & XQ\ol{X} \\
 && & & \dTo_{\id_X\otimes w\otimes\id_{\ol{X}}} \\
\dTo^{\ve(XQY)} && & & XQQ\ol{X} \\
 && & & \dTo_{\id_{XQ}\otimes\ve(Y)\otimes\id_{Q\ol{X}}} \\
XQY\ol{Y}\,\ol{Q}\, \ol{X} & \rTo_{\id_{XQY\ol{Y}}\otimes s^{-1}\otimes\id_{\ol{X}}} & &&
   XQY\ol{Y}Q\ol{X}
\end{diagram}
\]
which commutes in view of (\ref{bla}) and the assumption that $\2A$ is strict pivotal.
The other conditions on $\ve, \ol{\ve}$ are verified similarly, concluding that $\2E$ is a 
spherical bicategory. The details are omitted. The proof of semisimplicity is unchanged. 
\qed

We conclude our discussion of spherical categories by showing that a Frobenius algebra in
such a category is determined by almost as little data as in the case of $*$-categories.

\bprop \label{alt}
Let $\2A$ be a non-degenerate spherical category with simple unit. Let $(Q,v,w')$ be an
algebra in $\2A$ such that
\begin{itemize}
\item[(i)] $\displaystyle \dim\Hom(\11,Q)=1$.
\item[(ii)] There is an isomorphism $s:Q\rarr\ol{Q}$ such that 
\begin{equation}\label{f-ii}
  \ol{\ve}(Q)\mcirc\id_Q\otimes s=\ol{\ve}(\ol{Q})\mcirc s\otimes\id_Q=v'\circ w'
\end{equation}
with some non-zero $v': Q\rarr\11$. 
\end{itemize}
Then there is also $w: Q\rarr Q^2$ such that $(Q,v,v',w,w')$ is a canonical
Frobenius algebra.
\eprop
\prf We first remark that by (i) a non-zero $v': Q\rarr\11$ exists and is unique up to a
scalar. If there is a $s:Q\rarr\ol{Q}$ satisfying (\ref{f-ii}) with some $v'$ then this
obviously is the case for every choice of $v'$. Since $\Hom(\11,Q), \Hom(Q,\11)$ are
one-dimensional and in duality we have $v'\circ v=\lambda_2\id_\11$ with $\lambda_2\ne 0$.
We write $r'=v'\circ w': Q^2\rarr\11$. Using the fact that $v$ is the unit for the
multiplication $w'$ we find
\begin{equation} \label{vw}
   r'\mcirc \id_Q\otimes v=v'\mcirc w'\mcirc\id_Q\otimes v=v'. 
\end{equation}
Using (\ref{f-ii}), the duality equations for
$\ve, \ol{\ve}$ and property (3) in Definition \ref{pivotal} one easily shows 
$\id_Q\otimes s^{-1}\mcirc\ve(Q)=s^{-1}\otimes\id_Q\mcirc\ve(\ol{Q})$. We take this as the
definition of $r:\11\rarr Q^2$. One readily verifies that $r, r'$ satisfy the triangular
equations. Using the latter and (\ref{vw}) we compute
\begin{equation} \label{vw'}
   v= \id_Q\otimes r'\mcirc r\otimes\id_Q\mcirc v= \id_Q\otimes v'\mcirc r. 
\end{equation}
and similarly $v'\otimes\id_Q\mcirc r=v$. In the following computation the first and last
equalities hold by definition of $r, r'$ and the middle by sphericity, viz.\ property (3)
in Definition \ref{pivotal}.

\[ w=\quad
\begin{tangle}
\step[3.5]\object{Q}\step\object{Q}\\
\step[-.5]\mobj{r'}\step[.5]\coev\step[1.5]\id\step\id\\
\hh\id\step[1.5]\hcd\mobj{w'}\step\id\step\id\\
\hh\id\step[1.5]\id\step\hev\mobj{r}\step\id\\
\id\step[1.5]\mev\mobj{r}\\
\object{Q}
\end{tangle}
\quad=\quad
\begin{tangle}
\hstep\object{\ol{\ve}(\ol{Q})}\step[3]\object{Q}\step\object{Q}\\
\coev\step[1.5]\O{s^-1}\step\id\\
\id\step[1.5]\hcd\mobj{w'}\step\id\step\O{s^-1}\\
\O{s}\step[1.5]\id\step\hev\step[-1.5]\mobj{\ve(Q)}\step[2.5]\id\\
\id\step[1.5]\mev\mobj{\ve(Q)}\\
\object{Q}
\end{tangle}
\quad=\quad
\begin{tangle}
\object{Q}\step\object{Q}\step[3]\object{\ol{\ve}(Q)}\\
\id\step\O{s^-1}\step[1.5]\coev\\
\O{s^-1}\step\id\step\hcd\mobj{w'}\step[1.5]\id\\
\id\hstep\mobj{\ve(\ol{Q})}\hstep\hev\step\id\step[1.5]\O{s}\\
\mev\step[-.2]\mobj{\ve(\ol{Q})}\step[1.7]\id\\
\step[4.5]\object{Q}
\end{tangle}
\quad=\quad
\begin{tangle}
\object{Q}\step\object{Q}\step[3]\\
\id\step\id\step[1.5]\coev\mobj{r'}\\
\hh\id\step\id\step\hcd\mobj{w'}\step[1.5]\id\\
\hh\id\hstep\mobj{r}\hstep\hev\step\id\step[1.5]\id\\
\mev\step[-.2]\mobj{r}\step[1.7]\id\\
\step[4.5]\object{Q}
\end{tangle}
\]
This defines a comultiplication $w:Q\rarr Q^2$ whose coassociativity is obvious. Together
with $r=\id_Q\otimes s^{-1}\mcirc\ve(Q)=s^{-1}\otimes\id_Q\mcirc\ve(\ol{Q})$  
and (\ref{vw}) one shows $w\circ v=r$, and (\ref{vw'}) implies 
$v'\otimes\id_Q\mcirc w=\id_Q\otimes v'\mcirc w=\id_Q$, thus $(Q,v',w)$ is a comonoid.
Furthermore,
\begin{equation} \label{ww1}
\begin{tangle}
\cu\mobj{w}
\end{tangle}
\quad\equiv\quad
\begin{tangle}
\step[3.5]\object{Q}\step\object{Q}\\
\step[-.5]\mobj{r'}\step[.5]\coev\step[1.5]\id\step\id\\
\hh\id\step[1.5]\hcd\mobj{w'}\step\id\step\id\\
\hh\id\step[1.5]\id\step\hev\mobj{r}\step\id\\
\id\step[1.5]\mev\mobj{r}\\
\object{Q}
\end{tangle}
\quad=\quad
\begin{tangle}
\step[3.5]\object{Q}\hstep\object{Q}\\
\mobj{r'}\step[.5]\coev\step\id\hstep\id\\
\hh\hcd\step[1.5]\id\mobj{w'}\step\id\hstep\id\\
\hh\id\step\id\step[1.5]\hev\mobj{r}\hstep\id\\
\id\step\mev\mobj{r}\\
\object{Q}
\end{tangle}
\quad=\quad
\begin{tangle}
\hcd\mobj{w'}\step\id\\
\id\step\hev\mobj{r}\\
\end{tangle}
\quad,
\end{equation}
where we have used $r'=v'\circ w'$ and associativity of the multiplication. Similarly,
\begin{equation} \label{ww2}
\begin{tangle}
\cu\mobj{w}
\end{tangle}
\quad=\quad
\begin{tangle}
\id\step\hcd\mobj{w'}\\
\hev\mobj{r}\step\id
\end{tangle}
\quad,\quad\quad\quad
\begin{tangle}
\cd\mobj{w'}
\end{tangle}
\quad=\quad
\begin{tangle}
\id\step\hcoev\mobj{r'}\\
\hh\hcu\mobj{w}\step\id\\
\hh\hstep\id\step[1.5]\id
\end{tangle}
\quad=\quad
\begin{tangle}
\hcoev\mobj{r'}\step\id\\
\hh\id\step\hcu\mobj{w}\\
\hh\id\step[1.5]\id
\end{tangle}
\quad.
\end{equation}
Therefore
\[ 
\begin{tangle}
\hh\hcu\\
\hh\hcd
\end{tangle}
\quad=\quad
\begin{tangle}
\hh\step\id\step[1.5]\id\\
\hh\hstep\hcd\step\id\\
\hh\hstep\id\step\hev\\
\hh\hcd
\end{tangle}
\quad=\quad
\begin{tangle}
\hcd\step[1.5]\id\\
\hh\id\hstep\hcd\step\id\\
\hh\id\hstep\id\step\hev
\end{tangle}
\quad=\quad
\begin{tangle}
\hcd\step\id\\
\id\step\hcu
\end{tangle}
\quad,
\]
and the other part of condition (\ref{W3}) is proved analogously. Thus $(Q,v,v',w,w')$ is a
Frobenius algebra. Using the relations proved so far we compute further
\[
\begin{tangle}
\hh\hcd\mobj{w'}\\
\hh\hcu\mobj{w}
\end{tangle}
\quad=\quad
\begin{tangle}
\hh\hcd\hstep\hcoev\\
\hh\hcu\hstep\id\step\id\\
\hh\hstep\hev\step\id
\end{tangle}
\quad=\quad
\begin{tangle}
\hh\hcd\step\hcoev\\
\id\step\hddcu\step\id\\
\hh\hev\Step\id
\end{tangle}
\quad=\quad
\begin{tangle}
\hh\id\step\hcoev\\
\hh\hcu\step\id\\
\hh\hcd\step\id\\
\hh\hev\step\id
\end{tangle}
\quad=\quad
\begin{tangle}
\hstep\cd\\
\hh\hcd\mobj{w'}\step[1.5]\id\\
\hh\hev\step[.2]\mobj{r}\step[1.3]\id\\
\hh\step[2.5]\id
\end{tangle}
\]
Now, by assumption (i), $w'\circ r\in\Hom(\11,Q)$ satisfies $w'\circ r=\lambda_1 v$,
implying $w'\circ w=\lambda_1\id_Q$. Furthermore, 
$v'\circ w'\circ r=r'\circ r=\ol{\ve}(Q)\circ\ve(Q)=d(Q)$,  
thus $\lambda_1\lambda_2=d(Q)\ne 0$. Therefore $(Q,v,v',w,w')$ is a canonical Frobenius
algebra in $\2A$ and we are done. 
\qed

\brem 1. This result is in perfect accord with the classical definition according to which 
an algebra over a field $\7F$ is Frobenius iff it is isomorphic to $\hat{A}$ as a left 
(equivalently, right) $A$-module. Further support for our terminology will be supplied in 
Subsection \ref{ss-Frob}.

2. Instead of the existence of $s: Q\rarr\ol{Q}$ one might assume the existence of 
$r: \11\rarr Q^2$ satisfying the duality equations together with $r'=v'\circ w'$.
Unfortunately, this approach meets a problem. One easily shows the existence of
isomorphisms $s_1, s_2: Q\rarr\ol{Q}$ such that $\id_Q\otimes s_1\mcirc r=\ve(Q)$ and 
$s_2\otimes\id_Q\mcirc r=\ve(\ol{Q})$. Yet, it is unclear whether $s_1=s_2$ as
required. This condition can in fact be shown if $\2C$ is braided, $r'\circ c(Q,Q)=r'$
and $\theta(Q)=\id$. Here the twist $\theta$ is defined using the spherical structure 
\cite{y}. These are precisely the defining properties of a `rigid algebra' in the sense of
\cite{ko}. In view of the preceding remark, we find the terminology `Frobenius algebra'
more appropriate. 
\erem


\subsection{More on weak monoidal Morita equivalence}
In order to maintain the correspondence between canonical Frobenius algebras in $\2A$ and
tensor categories $\2B\approx\2A$ it is clear that we need the following

\bdefin $*$-categories $\2A, \2B$ are Morita equivalent if there is a Morita context $\2F$
which is a $*$-bicategory such that the equivalences 
$\ol{\2A}^{p\oplus}\stackrel{\otimes}{\simeq}\END_\2F(\6A)$, 
$\ol{\2B}^{p\oplus}\stackrel{\otimes}{\simeq}\END_\2F(\6B)$,
are equivalences of $*$-categories.
If $\2A, \2B$ are spherical categories they are weakly monoidally Morita equivalent if
there is a Morita $\2F$ which is spherical such that the above equivalences are
equivalences of spherical categories. 
\edefin

We summarize our findings on $*$- and spherical categories.
\btheor If $\2A$ is a $*$-category and $(Q,v,v^*,w,w^*)$ is a canonical Frobenius algebra
then $\2E_*$ is a $*$-bicategory. If $\2A$ is spherical (and non-degenerate (and
semisimple)) and $(Q,v,v',w,w')$ is a canonical Frobenius algebra then $\2E$ is spherical
(and non-degenerate (and semisimple)). In both cases $\2B=\END_{\2E_*}(\6B)$ is weakly
monoidally Morita equivalent to $\2A$.
\etheor

As a first application of weak monoidal Morita equivalence for spherical or $*$-categories 
we prove the analogue of a well known result in subfactor theory, cf.\ e.g.,
\cite{vfr-su,ek2}. The proof extends to the present setting without any change. 
\bprop \label{eqofdims}
Let $\2A, \2B$ be a finite dimensional semisimple spherical tensor categories over $\7F$
(or $*$-categories) with simple units. If $\2A, \2B$ are weakly monoidally Morita
equivalent ($\2A\approx\2B$) then they have the same dimension in the sense of
(\ref{catdim}).
\eprop
\prf Let $\2E$ be a Morita context for $\2A\approx\2B$. Let $I, K$ be (finite) sets
labeling the isomorphism classes of simple $\6A-\6A$-morphisms and $\6B-\6A$-morphisms,
respectively, and let $\{X_i, i\in I\},\ \{Y_k, k\in K\}$ be representers of the
latter. The integers 
\[ N_i^k=\dim\Hom_\2E(Y_k,X_i J) \]
do not depend on the chosen representers, and by duality we have
\[ N_i^k=\dim\Hom_\2E(Y_k\oj ,X_i). \]
Thus $X_i J\cong\bigoplus_k N_i^k Y_k$ and 
$Y_k\oj \cong\bigoplus_i N_i^k X_i$. Using additivity and multiplicativity of the dimension
function we compute 
\bean  d(J) \sum_{i\in I} d(X_i)^2 &=& \sum_{i\in I} d(X_i) d(X_i J) =
  \sum_{i\in I \atop k\in K} N_i^k d(X_i)d(Y_k)  \\
  &=& \sum_{k\in K} d(Y_k) d(Y_k\oj )\ =\ d(\oj ) \sum_{k\in K} d(Y_k)^2. 
\eean
Since $d(J)=d(\oj )\ne 0$ we conclude
$\dim\HOM_\2E(\6A,\6A)=\dim\HOM_\2E(\6B,\6A)$. Entirely analogous arguments yield
$\dim\HOM_\2E(\6B,\6A)=\dim\HOM_\2E(\6B,\6B)$ and therefore $\dim\2A=\dim\2B$. Of course, also
$\dim\HOM_\2E(\6A,\6B)$ has the same dimension.
\qed

\brem 1. Note that the categories $\HOM_\2E(\6A,\6B)$ and $\HOM_\2E(\6B,\6A)$ are not
tensor categories, thus the intrinsic notion of dimension of \cite{bw1,lro} does not
apply and a priori it does not make sense to speak of their dimensions. But every object
of $\HOM_\2E(\6A,\6B)$ or $\HOM_\2E(\6B,\6A)$ is a 1-morphism in $\2E$ and as such has a 
dimension. It is this dimension which is intended in the above statement. 

2. Given a linear Morita context $\2E$, the common dimension of the four categories
$\HOM_\2E(\cdot,\cdot)$ is also called the dimension of $\2E$. 

3. A less elementary application of Morita equivalence is the fact that weakly Morita
equivalent spherical categories define the same state sum invariant (in the sense of
\cite{bw1, gk}) for 3-manifolds. The proof is sketched in Section \ref{s-last}. 

4. Let $\2A$ be a tensor category and $(Q,\ldots)$ a canonical Frobenius algebra in
$\2A$. It should be obvious that the tensor category $\2B=\Hom_\2E(\6B,\6B)$ can be
defined directly, avoiding the construction of the entire bicategory $\2E$. When we are 
interested only in $\2B$ we may suppress the $J,\oj$ in $(\oj XJ, p)$. Thus the objects of
$\2B$ are pairs $(X, p)$, where $X\in\obj\2A$ and $p\in\Hom_\2A(XQ,QX)$ satisfies
$p\bullet p=p$. The morphisms are given by 
\[ \Hom_\2B((X,p), (Y,q))=\{ s\in\Hom_\2A(XQ,QY)\ | \ s=q\bullet s\bullet p \}, \]
the tensor product of objects by $(X,p)\otimes (Y,q)=(XQY, p\times q)$, etc. (Here 
$\bullet, \times$ are defined as in Proposition \ref{constr}.) For many purposes, like the
study of the categorical version \cite{mue15} of `$\alpha$-induction' 
\cite{lre, xu1, bek}, this is sufficient. But proceeding in this way the weak monoidal
Morita equivalence $\2A\approx\2B$ become obscure and even the proof of Proposition
\ref{eqofdims} (which is an instance of the `2x2-matrix trick') seems very difficult
without the Morita context $\2E$. 
\erem


\section{Examples}\label{s-examples}
In this section we will consider two classes of examples: Classical Frobenius algebras
over a field, in particular Hopf algebras, and subfactors with finite index. Both examples
are essential for obtaining a deep understanding of categorical Frobenius theory. Whereas
most of of our discussion essentially amounts to reformulating known facts, in the final
subsection we will obtain a new result, viz.\ the weak monoidal Morita equivalence of
$H-\mod$ and $\hat{H}-\mod$ for certain Hopf algebras. This result relies on input from
both the classical Frobenius theory and subfactor theory. (That the latter are related is
not new and has been discussed, e.g., in \cite{kad}.)


\subsection{Frobenius and Hopf Algebras over Fields} \label{ss-Frob}
Here we briefly review a beautiful recent result of L.\ Abrams \cite{abrams} which implies
that in the case $\2A=\7F$-Vect (which we treat as strict, following common usage)
our notion of Frobenius algebras is equivalent to the classical one. This justifies the
terminology. Note, however, that this was not our main motivation for Definition
\ref{d-Frob}. 

Let $A$ be a finite dimensional (associative, with unit) algebra over a field $\7F$.
The dual vector space $\hat{A}$ comes with two natural coalgebra structures
\[ \langle \hat{\Delta}_1(\alpha), x\otimes y\rangle=\langle\alpha,xy\rangle, \quad\quad
   \langle \hat{\Delta}_2(\alpha), x\otimes y\rangle=\langle\alpha,yx\rangle, \]
both of which have the counit $\hat{\ve}(\alpha)=\langle\alpha,1\rangle$. 
Given an isomorphism $\Phi: A\rarr\hat{A}$ of vector spaces we can provide $A$ with a
coalgebra structure by
\[ \Delta= \Phi^{-1}\otimes\Phi^{-1}\circ \hat{\Delta}\circ\Phi, 
   \quad \ve= \hat{\ve}\circ\Phi, \]
where $\hat{\Delta}=\hat{\Delta}_1$ or $\hat{\Delta}=\hat{\Delta}_2$.

Whenever $A$ admits an isomorphism  $\Phi: A\rarr\hat{A}$ of left (equivalently right)
$A$-modules (with the natural left or right $A$-module structures) $A$ is called a
Frobenius algebra. We prefer the following equivalent definition, see \cite{kad}.

\bdefin A finite dimensional algebra (associative, with unit) over a field $\7F$ is a
Frobenius algebra if it admits a linear form $\phi: A\rarr\7F$ which is non-degenerate (in
the sense that the bilinear form $b(x,y)=\phi(xy)$ is non-degenerate). 
\edefin
The linear form $\phi$ gives rise to two isomorphisms between $A$ and its dual $\hat{A}$
via
\[ \Phi_1: x\mapsto \phi(x\cdot),\quad\quad \Phi_2: x\mapsto \phi(\cdot x). \]
Clearly, $\Phi_1=\Phi_2$ iff $\phi$ is a trace. By the preceding discussion we thus have
four canonical ways of providing $A$ with a coalgebra structure, depending on which
combination of $\Phi_1/\Phi_2, \hat{\Delta}_1/\hat{\Delta}_2$ we use. (If $\phi$ is a
trace these possibilities reduce to two and if $A$ is commutative we are left with a
unique one. The commutative case is discussed in \cite{abrams0}.) In any case, the counit
is given by $\ve=\phi$.

\btheor\cite{abrams}\label{t-abrams}
Let $A$ be a Frobenius algebra with given $\phi\in\hat{A}$. Let
$\Delta_1, \Delta_2$ be the coproducts defined as above using the combinations
$(\hat{\Delta}_2, \Phi_1)$ and $(\hat{\Delta}_2, \Phi_2)$, respectively. Then
$\Delta_1=\Delta_2$ and with $\Delta=\Delta_1$ the following diagrams commute:
\[
\begin{diagram}
  A\otimes A & \rTo^{m} & A \\
  \dTo^{\id\otimes\Delta} && \dTo_{\Delta} \\
  A\otimes A\otimes A & \rTo_{m\otimes\id} & A\otimes A
\end{diagram}
\quad\quad\quad
\begin{diagram}
  A\otimes A & \rTo^{m} & A \\
  \dTo^{\Delta\otimes\id} && \dTo_{\Delta} \\
  A\otimes A\otimes A & \rTo_{\id\otimes m} & A\otimes A
\end{diagram}
\]
Thus $(A,\11,\ve,\Delta,m)$ is a Frobenius algebra in the sense of Definition \ref{d-Frob}.
Every Frobenius algebra in $\7F$-Vect arises in this way.
\etheor
\brem \label{r-frob}
1. In the proof, one first shows that the first diagram commutes with $\Delta=\Delta_1$
and the second with $\Delta=\Delta_2$.
Using these facts one proves $\Delta_1=\Delta_2$, thus $(A,\11,\ve,\Delta,m)$ is a
Frobenius algebra in Vect. Our version of the converse statement differs slightly from the
one in \cite{abrams}, but it is easily seen to be equivalent. 

2. An obvious consequence of the alternative characterization of Frobenius algebras is
that the dual vector space of a Frobenius algebra is again a Frobenius algebra. 

3. A special case of this had been shown earlier by Quinn in the little noticed appendix
of \cite{q}. He defines an `ambialgebra' (in the category Vect) as an algebra and
coalgebra satisfying commutativity of the above diagrams plus symmetry conditions on
$\Delta(1)$ and $\ve\circ m$. He states that these ambialgebras are the same as
symmetric algebras (i.e.\ algebras admitting a non-degenerate trace). This result is
intermediate in generality between those of \cite{abrams0} and \cite{abrams}. 

4. Let $X$ be a finite dimensional $\7F$-vector space with dual vector space
$\ol{X}$. Then the Frobenius algebra $(X\ol{X},\ldots)$ defined as in Lemma \ref{motiv} is
well known to be just the matrix algebra $\End\,X$. Since there are Frobenius algebras
which are not isomorphic to some $M_n(\7F)$, already the category $\7F$-Vect provides an
example of a tensor category where the Frobenius algebras are not exhausted by those of
the form $(X\ol{X},\ldots)$. 

5. In view of $\dim\Hom_{\mbox{Vect}}(\11,Q)=\dim H$, non-trivial Frobenius algebras 
in Vect are not irreducible.
\erem

In order to understand when a Frobenius algebra in $\7F$-Vect is canonical we need the
more general notion of a Frobenius extension \cite{kad}.
\bdefin A ring extension $A/S$ is a Frobenius extension iff there exists a Frobenius
system $(E,x_i,y_i)$, where $E\in\Hom_{S-S}(A,S)$ (i.e.\ 
$E(abc)=aE(b)c\ \forall b\in A, a,c\in S$), $|I|<\infty$ and $x_i, y_i\in A,\ i\in I$ such
that
\begin{equation} \label{e-fr}
  \sum_{i\in I} x_i E(y_ia) = a = \sum_{i\in I} E(ax_i)y_i \quad \forall a\in A. 
\end{equation}
We call $\sum_i x_i\otimes y_i\in A\otimes_S A$ the Frobenius element and 
$[A:S]_E=\sum_i x_iy_i\in Z(A)$ the $E$-index. A Frobenius extension $A/S$ of
$\7F$-algebras is called strongly separable iff $E(1)=1$ and $[A:S]_E\in K^*1$. 
\edefin

We are interested in the case where $\7F$ is a field and $A$ is finite dimensional over
$\7F$. Then $A/\7F$ is a Frobenius extension iff $A$ is a Frobenius algebra, cf.\
\cite[Proposition 4.8]{kad}. In this case $E=\phi$. If $\{x_i\}$ is any basis of $A$ then
$\{y_i\}$ satisfies (\ref{e-fr}) iff it is the dual basis: $\phi(y_ix_j)=\delta_{ij}$.

\bprop \label{p-frob1}
A Frobenius algebra $A$ in $\7F$-Vect is canonical in the sense of Definition 
\ref{d-canon} iff the Frobenius extension $A/\7F$ is strongly separable modulo
renormalization.
\eprop
\prf It is obvious that the morphism $v'\circ v$ is given by $c\mapsto c\phi(1)$.
Now, the Frobenius property (\ref{W3}) implies 
$\Delta(x)=\Delta(1)(1\otimes x)=(x\otimes 1)\Delta( 1)$ and therefore
$m\Delta(x)=(m\Delta( 1))x=x(m\Delta( 1))$. Thus the morphism $w'\circ w=m\circ\Delta$ 
is given by multiplication with the central element $m\Delta(1)$. 
We will show that $m\Delta(1)=[A:\7F]_\phi$. Thus, if $A/\7F$ is strongly
separable then $v'\circ v\in\7F^*\id_\11$ and $w'\circ w\in \7F^*\id_A$. The converse
holds since $\phi(1)\ne 0$ allows to renormalize such that $\phi(1)=1$.

Let $\{x_i\}$ be a basis of $A$ and $\{y_i\}$ dual in the sense
$\phi(y_ix_j)=\delta_{ij}$. Then $\sum_i x_i\otimes y_i\in A\otimes A$ is the Frobenius
element. For $a,b\in A$ we compute 
\[ \left( (\Phi_1\otimes\Phi_1)(\Delta( 1))\right) (a\otimes b)=
   \left(\hat{\Delta}_2\Phi_1( 1)\right) (a\otimes b)=\Phi_1( 1)(ba)=\phi(ba)=
   \sum_i \phi(ax_i)\phi(by_i), \]
where we used (\ref{e-fr}). Thus $\Delta(1)$ equals the Frobenius element and 
$m\Delta( 1)=\sum_i x_iy_i=[A:\7F]_\phi$. 
\qed 

\brem In the commutative case Frobenius algebras satisfying the above equivalent
conditions were called `superspecial' in \cite{q}. Note furthermore that semisimplicity of
$A$ is equivalent to the weaker condition of invertibility of $\sum_i x_iy_i$ (proven in
\cite{abrams0} for the commutative and in \cite{q} for the symmetric case). Thus canonical
Frobenius algebras are semisimple.
\erem

It is well known \cite{ls} that every finite dimensional Hopf algebra over a field $\7F$
is a Frobenius algebra. Our aim in the remainder of this subsection is to clarify when
these Frobenius algebras are canonical. We recall some well known facts. For any finite
dimensional Hopf algebra $H$ one can prove \cite{ls} that the subspaces  
$I_L(H)=\{ y\in H\ | \ xy=\ve(x)y\ \ \forall x\in H\}$ and  
$I_R(H)=\{ y\in H\ | \ yx=\ve(x)y\ \ \forall x\in H\}$ are one dimensional satisfy
$S(I_L(H))=I_R(H)$. Furthermore, every non-zero $\varphi_L\in I_L(\hat{H})$ and
$\varphi_R\in I_R(\hat{H})$ is a non-degenerate functional on $H$.

In view of Theor.\ \ref{t-abrams}, both $\varphi_L$ and $\varphi_R$ give rise to coalgebra
structures $(H,\tilde{\Delta}_{L/R},\tilde{\ve}_{L/R})$ on the vector space $H$ and
therefore to Frobenius algebras 
$\5Q_{L/R}=(H,m,\eta,\tilde{\Delta}_{L/R},\tilde{\ve}_{L/R})$. (We denote the
Frobenius coproduct by $\tilde{\Delta}_{L/R}$ to avoid confusion with the Hopf algebra 
coproduct $\Delta$ of $H$.) 

\bprop Let $H$ be a finite dimensional Hopf algebra with non-zero right integrals 
$\Lambda\in I_R(H), \varphi\in I_R(\hat{H})$. Then the following are equivalent:
\begin{itemize}
\item[(i)] The Frobenius algebra $\5Q_L=(H,m, \eta,\tilde{\Delta}_L,\tilde{\ve}_L)$ in
$\7F$-Vect is canonical in the sense of Definition \ref{d-canon}.
\item[(ii)] $\langle\varphi,1\rangle\ne 0$ and $\langle\ve,\Lambda\rangle\ne 0$.
\item[(iii)] $H$ is semisimple and cosemisimple.
\end{itemize}
\eprop
\prf (ii)$\Leftrightarrow$(iii). By \cite{ls}, $H$ is semisimple iff $\ve(\Lambda)\ne 0$ and 
cosemisimple iff $\varphi(1)\ne 0$. \\
(i)$\Leftrightarrow$(ii). By \cite[Proposition 6.4]{kad}, a Frobenius system for $H/\7F$
is given by the triple $(\varphi,S^{-1}(\Lambda_{(2)}),\Lambda_{(1)})$. (This is to say
that the Frobenius element is given by 
$\sum_i x_i\otimes y_i=(S^{-1}\otimes\id)\Delta^\op(\Lambda)$.)
But then it is obvious that
$[H:\7F]_\varphi=S^{-1}(\Lambda_{(2)})\Lambda_{(1)}=\ve(\Lambda)1$. 
Thus $m\circ\tilde{\Delta}_L=[H:\7F]_f=\ve(\Lambda)1$, and the equivalence
(i)$\Leftrightarrow$(ii) follows from Proposition \ref{p-frob1}.
\qed

\brem In view of $S(I_L)=I_R$ and $S(1)=1$ we have 
$\langle\varphi_L,1\rangle\ne 0\ \Leftrightarrow\ \langle\varphi_R,1\rangle\ne 0$. Thus it
does not matter where the integrals appearing in condition (ii) are left or right
integrals. Similarly, in (i) we can write $\5Q_R$ instead of $\5Q_L$. In the canonical
case these substitutions are vacuous since semisimple Hopf algebras are unimodular, i.e.\
$I_L=I_R$.
\erem

\bcoro Let $H$ be a finite dimensional Hopf algebra over $\7F$ and $\5Q$ the corresponding
Frobenius algebra (in $\7F$-Vect). Let $\Lambda,\varphi$ be both either left or right
integrals in $H, \hat{H}$, respectively. Then there is the following identity of numerical
invariants of $\5Q$ and $H$:
\begin{equation} \label{inv}
   v'\circ w'\circ w\circ v=  \frac{\langle\varphi, 1\rangle\langle\ve,\Lambda\rangle}
   {\langle\varphi,\Lambda\rangle} \ \in\End(\11)\equiv\7F.
\end{equation}
Whenever this number is nonzero $H$ is semisimple and cosemisimple and (\ref{inv})
coincides with $\dim H\cdot 1_\7F$.
\ecoro
\prf Eq.\ (\ref{inv}) follows from the above computation of $v'\circ v$ and $w'\circ
w$. If (\ref{inv}) is non-zero then $H$ is semisimple and cosemisimple, and by \cite{eg2}
the antipode is involutive. By \cite[Theorem 2.5]{lar2}, (\ref{inv}) coincides with
$tr(S^2)$ and therefore with $\dim H\cdot\11_\7F$.
\qed

By the preceding result semisimple and cosemisimple Hopf algebras provide examples of
canonical Frobenius algebras in $\7F$-Vect. By Remark \ref{r-frob}.4 they are not
irreducible. We will now show how one can associate a canonical and irreducible
Frobenius algebra with a Hopf algebra which is semisimple and cosemisimple.


\subsection{Hopf algebras: Frobenius algebras in $H-\mod$} \label{ss-Hopf}
In the theory of quantum groups the following result is known as `strong left invariance'
(for $b=1$ or $c=1$ it reduces to left invariance:
$(\id\otimes\varphi)(\Delta(b))=\varphi(b)1$), but it is also true for all finite
dimensional Hopf algebras. We include the proof since we are not aware of a convenient
reference. 

\blemma \label{l-sli}
Let $H$ be a finite dimensional Hopf algebra and $\varphi\in I_L(\hat{H})$. Then
\[ (\id\otimes\varphi)((1\otimes c)\Delta(b))
   =(\id\otimes\varphi)((S\otimes\id)(\Delta(c))(1\otimes b))
   \quad\forall b,c\in H. \]
\elemma
\prf For every $x,y\in H$ we have
\[ x\otimes y=\sum_i (u_i\otimes 1)\Delta(v_i) \quad\mbox{with}\quad 
   \sum_i u_i\otimes v_i=(x\otimes 1)(S\otimes\id)(1\otimes y)=xS(y_{(1)})\otimes y_{(2)},
\]
as is verified by a trivial computation. Using this representation and left invariance of
$\varphi$ we have 
$(\id\otimes\varphi)(x\otimes y)=\sum_i u_i(\id\otimes\varphi)\Delta(v_i)=\sum_i
 u_i\varphi(v_i)=(\id\otimes\varphi)(\sum_i u_i\otimes v_i)$.
Applying this to $x\otimes y=(1\otimes c)\Delta(b)=b_{(1)}\otimes cb_{(2)}$ we find
\bean \sum_i u_i\otimes v_i &=& b_{(1)}S((cb_{(2)})_{(1)})\otimes (cb_{(2)})_{(2)} 
  = b_{(1)}S(c_{(1)}b_{(2)})\otimes c_{(2)}b_{(3)} \\
  &=& b_{(1)}S(b_{(2)})S(c_{(1)})\otimes c_{(2)}b_{(3)} = S(c_{(1)})\otimes c_{(2)}b \eean
and therefore 
\[ (\id\otimes\varphi)((1\otimes c)\Delta(b))=(\id\otimes\varphi)(S(c_{(1)})\otimes c_{(2)}b)
   = (\id\otimes\varphi)((S\otimes\id)(\Delta(c)(1\otimes b)), \]
as desired.
\qed

\bprop \label{p-mtilde}
Let $H$ be a finite dimensional Hopf algebra and $\varphi_L\in I_L(\hat{H})$. 
We write $\hat{m}$ for the multiplication of $\hat{H}$ and define $\2F: H\rarr\hat{H}$ by
$\2F(a)(\cdot)=\varphi(\cdot a)$. Then the map 
$\tilde{m}=\2F^{-1}\hat{m}(\2F\otimes\2F): H\otimes H\rarr H$ satisfies
\begin{itemize}
\item[(i)]
$\DS \varphi(c\tilde{m}(a\otimes b))=(\varphi\otimes\varphi)(\Delta(c)(a\otimes b))
=(\varphi\otimes\varphi)((1\otimes c)(S^{-1}\otimes\id)(\Delta(b))(a\otimes 1))
=(\varphi\otimes\varphi)((c\otimes 1)(1\otimes S^{-1}(b))\Delta(a))
\quad\forall a,b,c\in H$. 
\item[(ii)]
$\DS \tilde{m}(\Delta(c)x)=c\tilde{m}(x) \quad \forall c\in H, x\in H\otimes H$.
\item[(iii)]
$\DS \tilde{m}(a\otimes b)=\varphi(S^{-1}(b_{(1)})a)b_{(2)}=a_{(1)}\varphi(S^{-1}(b)a_{(2)})$.
\end{itemize}
\eprop
\prf (i) We compute
\begin{eqnarray*} \lefteqn{ \underline{\varphi(c\tilde{m}(a\otimes b))}
  = \langle \2F\tilde{m}(a\otimes b), c\rangle 
  = \langle \hat{m}(\2F(a)\otimes\2F(b)),c \rangle } \\
  && = \langle \2F(a)\otimes\2F(b),\Delta(c) \rangle 
   = \underline{(\varphi\otimes\varphi)(\Delta(c)(a\otimes b))} \\
  && = \varphi\left( S^{-1}[(\id\otimes\varphi)((S\otimes\id)(\Delta(c))(1\otimes b))]a\right) \\
  && = \varphi\left( S^{-1}[(\id\otimes\varphi)((1\otimes c)\Delta(b))]a\right) \\
  && =\underline{(\varphi\otimes\varphi)((1\otimes c)(S^{-1}\otimes\id)(\Delta(b))(a\otimes 1))} \\
  && = (\varphi\otimes\varphi)((1\otimes c)\sigma((S\otimes\id)(\Delta(S^{-1}(b))))(a\otimes 1)) \\
  && = (\varphi\otimes\varphi)((c\otimes 1)(S\otimes\id)(\Delta(S^{-1}(b)))(1\otimes a)) \\
  && = \varphi(c[(\id\otimes\varphi)((S\otimes\id)(\Delta(S^{-1}(b)))(1\otimes a))]) \\
  && = \varphi(c[(\id\otimes\varphi)((1\otimes S^{-1}(b))\Delta(a))]) \\
  && = \underline{(\varphi\otimes\varphi)((c\otimes 1)(1\otimes S^{-1}(b))\Delta(a))}.
\end{eqnarray*}
In the first two lines we used the definitions of $\2F$ and $\tilde{m}$, whereas the 6th
and 11th equalities follow by application of Lemma \ref{l-sli} to the expression in square
brackets. The remaining identities result from trivial rearrangements using 
$(\varphi\otimes\varphi)=\varphi(\id\otimes\varphi)=\varphi(\varphi\otimes\id)$.

(ii) Let $x=a\otimes b$. Twofold use of the first equality in (i) yields
\[ \varphi(d\tilde{m}(\Delta (c)x)) = (\varphi\otimes\varphi)(\Delta(d)\Delta(c)x) 
 =(\varphi\otimes\varphi)(\Delta(dc)x) =\varphi(dc\tilde{m}(x)), \]
which holds for all $c,d\in H, x\in H\otimes H$. The claim now follows by non-degeneracy
of $\varphi$.

(iii) We can rewrite (i) as
\[ \varphi(c\tilde{m}(a\otimes b))
  =\varphi(c(\varphi\otimes\id)(S^{-1}\otimes\id)(\Delta(b))(a\otimes 1))) 
  =\varphi(c(\id\otimes\varphi)((1\otimes S^{-1}(b))\Delta(a))). \]
Now we appeal again to non-degeneracy of $\varphi$ and rewrite in Sweedler notation. 
\qed

\btheor \label{t-hopffrob}
Let $H$ be a finite dimensional Hopf algebra over $\7F$. Let $\Lambda, \varphi$ be
left integrals in $H$ and $\hat{H}$, respectively, normalized such that
$\langle\varphi,\Lambda\rangle=1$. 
Let $Q\in H-\mod$ be the left regular representation, viz.\ $H$ acting on itself by
$\pi_Q(a)b=ab$. The linear maps 
\bea v: & \7F\rarr Q, & c\mapsto c\Lambda, \nn\\
  v': & Q\rarr\7F, & x\mapsto \ve(x), \nn\\
  w: & Q\rarr Q\otimes Q, & x\mapsto \Delta(x), \nn\\
  w': & Q\otimes Q \rarr Q, & x\otimes y\mapsto \tilde{m}(x\otimes y)
\label{m1}
\eea
are morphisms in $H-\mod$ and $(Q,v,v',w,w')$ is an irreducible Frobenius algebra in 
$H-\mod$. It is canonical iff $H$ is semisimple and cosemisimple.
\etheor
\prf In order to show that the maps defined above are morphisms in the category
$H-\mod$ we must show that they intertwine the $H$-actions. This follows from the
following diagrams, where $\pi_1=\ve$ is the tensor unit:
\newarrow{Mapsto} |--->
\[
\begin{diagram} c & \rMapsto^{v} & c\Lambda \\
  \dMapsto^{\pi_1(z)} && \dMapsto_{\pi_Q(z)} \\
  \ve(z)c & \rMapsto_{v} & \ve(z)c\Lambda & =cz\Lambda.
\end{diagram}
\quad\quad
\begin{diagram} x & \rMapsto^{v'} & \ve(x) \\
  \dMapsto^{\pi_Q(z)} && \dMapsto_{\pi_1(z)} \\
  zx & \rMapsto_{v'} & \ve(zx) & =\ve(x)\ve(z).
\end{diagram}
\]

\[
\begin{diagram} x & \rMapsto^{w} & \Delta(x) \\
  \dMapsto^{\pi_Q(z)} && \dMapsto_{\pi_{Q\otimes Q}(z)} \\
  zx & \rMapsto_{w} & \Delta(zx) & =\Delta(z)\Delta(x).
\end{diagram}
\quad\quad
\begin{diagram} x\otimes y & \rMapsto^{w'} & \tilde{m}(x\otimes y) \\
   \dMapsto^{\pi_{Q\otimes Q}(z)} && \dMapsto_{\pi_Q(z)} \\
  \Delta(z)(x\otimes y) & \rMapsto_{w'} & \tilde{m}(\Delta(z)(x\otimes y))=z\tilde{m}(x\otimes y) 
\end{diagram}
\]
Commutativity of the lower right diagram follows from Proposition \ref{p-mtilde} (ii).

The equations $v'\otimes\id_Q\circ w=\id_Q\otimes v'\circ w=\id_Q$ and
$w\otimes\id_Q\circ w=\id_Q\otimes w\circ w$ are obvious since $(Q,w,v')$ coincides with
the coalgebra structure of $H$. That $w':Q^2\rarr Q$ is associative is evident in view of 
$\tilde{m}=\2F^{-1}\hat{m}(\2F\otimes\2F)$ and associativity of $\hat{m}$. Furthermore,
$\2F(\Lambda)(a)=\varphi(a\Lambda)=\ve(a)\varphi(\Lambda)=\ve(a)$, thus
$\2F(\Lambda)=\ve=1_{\hat{H}}$ is the unit for $\tilde{m}$ and $(Q,v,w')$ is a monoid. 

Applying $\ve$ to Proposition \ref{p-mtilde} (iii) we obtain
$(\ve\tilde{m})(a\otimes b)=\varphi(S^{-1}(b)a)$. Comparing with the formulae for
$\tilde{m}$ we find 
\begin{eqnarray*} \tilde{m}(a\otimes b) &=&
   \varphi(S^{-1}(b_{(1)})a)b_{(2)} = (\ve\tilde{m})(a\otimes b_{(1)})b_{(2)} \\
  &=& a_{(1)}\varphi(S^{-1}(b)a_{(2)}) = a_{(1)}(\ve\tilde{m})(a_{(2)}\otimes b),
\end{eqnarray*}
or in diagrams
\[
\begin{tangle}
\hstep\unit\step[.3]\obj{\ve}\step[1.2]\id\\
\hh\step[-.8]\mobj{\tilde{m}}\step[.8]\hcd\step\id\\
\hh\id\step\hcu\mobj{\Delta}\\
\object{a}\step[1.5]\object{b}
\end{tangle}
\quad=\quad
\begin{tangle}
\hcd\mobj{\tilde{m}}\\
\hh\id\step\id\\
\object{a}\step\object{b}
\end{tangle}
\quad=\quad
\begin{tangle}
\id\step[1.5]\unit\step[.3]\obj{\ve}\\
\hh\id\step\hcd\mobj{\tilde{m}}\\
\hh\step[-.8]\mobj{\Delta}\step[.8]\hcu\step\id\\
\hstep\object{a}\step[1.5]\object{b}
\end{tangle}
\]
Using the first of these equalities twice to compute
\[
\begin{tangle}
\hh\step[-.8]\mobj{\Delta}\step[.8]\hcu\\
\hh\step[-.8]\obj{\tilde{m}}\step[.8]\hcd
\end{tangle}
\quad=\quad
\begin{tangle}
\hh\step[1.5]\hcu\\
\hstep\unit\step[1.5]\id\\
\hh\hcd\step\id\\
\hh\id\step\hcu
\end{tangle}
\quad=\quad
\begin{tangle}
\hstep\unit\step[1.5]\id\step[1.5]\id\\
\hh\hcd\step\id\step[1.5]\id\\
\hh\id\step\hcu\step[1.5]\id\\
\id\step[1.5]\cu\\
\end{tangle}
\quad=\quad
\begin{tangle}
\hcd\step\id\\
\id\step\hcu
\end{tangle}
\]
we have proven one of the Frobenius conditions (\ref{W3}) and the other one follows
in the same vein.

That the Frobenius algebra $\5Q$ is irreducible follows from the obvious isomorphism of
vector spaces $\Hom_{H-\mod}(\11,Q)\cong I_L$ together with $\dim I_L=1$. Finally, we
compute 
\[ (\tilde{m}\Delta)(a)=\tilde{m}(a_{(1)}\otimes a_{(2)})
   =a_{(1)}\varphi(S^{-1}(a_{(3)})a_{(2)})=a\varphi(1). \]
Thus $\5Q$ is canonical iff $\ve(\Lambda)\ne 0$ and $\varphi(1)\ne 0$, which is the case
iff $H$ is semisimple and cosemisimple. 
\qed

\bexam Let $H=\7F(G)$ be the algebra of $\7F$-valued functions on a finite group $G$ with
the usual Hopf algebra structure. With $H=\mbox{span}\{ \delta_g, g\in G\}$, where
$\delta_g(h)=\delta_{g,h}$, the integrals $\Lambda\in H, \varphi\in\hat{H}$ are 
\bean \Lambda &=& \delta_e, \\
   \langle \varphi, \delta_g\rangle &=& 1. \eean
Then we find
\[ \langle \2F(\delta_g),\delta_h\rangle= \langle \varphi, \delta_g\delta_h\rangle=
\delta_{g,h}\langle \varphi, \delta_g\rangle= \delta_{g,h} \]
and thus $\2F(\delta_g)=u_g$, where $\{u_g, g\in G\}$ is the usual basis in
$\hat{H}=\7FG$. We obtain
\bean \11_Q &=& \delta_e=\Lambda, \\
  m_Q(\delta_g,\delta_h) &=& \delta_{gh}, \eean
and thus the Frobenius algebra $Q$ associated to $H=\7F(G)$ by our prescription coincides with
the one given in \cite[A.4.5]{q}. In a similar fashion one sees that the Frobenius algebra
associated with $\7FG$ has the coalgebra structure of $\7FG$ and the algebra structure of
$\7F(G)$ under the correspondence $u_g\leftrightarrow \delta_g$. 
\eexam

\brem Our original proof of Theorem \ref{t-hopffrob} has improved considerably as a
consequence of discussions with L. Tuset. In the joint work \cite{mue-t} we examine to
which extent the above results carry over to  not necessarily finite dimensional algebraic 
quantum groups. If an algebraic quantum group $(A,\Delta)$ is discrete, there exists a monoid 
$(\pi_L,\tilde{m},\tilde{\eta})$ in the category $\mbox{Rep}(A,\Delta)$ of non-degenerate 
$*$-representations. If $(A,\Delta)$ is compact (=unital) then there is a comonoid 
$(\pi_L,\tilde{\Delta},\tilde{\ve})$ in a $\mbox{Rep}(A,\Delta)$. Both the monoid and comonoid 
structures exist only if $(A,\Delta)$ is finite dimensional, in which case they coincide with 
those considered above. Therefore, in the infinite dimensional situation one does not obtain a 
Frobenius algebra in $\mbox{Rep}(A,\Delta)$, but a `regularized' version of the Frobenius 
relation (\ref{W3}) can still be proven.
\erem


\subsection{Morita equivalence of $H-\mod$ and $\hat{H}-\mod$} \label{ss-Hopf2}
In this subsection $\7F$ is an arbitrary algebraically closed field.
If $H$ is a finite dimensional semisimple and cosemisimple Hopf algebra over $\7F$,
Theorem \ref{t-hopffrob} gives rise to a canonical Frobenius algebra $\5Q$ in
$H-\mod$. Applying Theorem \ref{main0}, we obtain a Morita context $\2E$, and it is
natural to ask what can be said about the tensor category $\2B=\END_\2E(\6B)$, which is
Morita equivalent to $H-\mod$ by construction. We may and will assume $\2E$ to be strict,
i.e.\ a 2-category. The aim of this subsection is to prove the following.

\btheor \label{t-h-hhat}
Let $H$ be a finite dimensional semisimple and cosemisimple Hopf algebra over an
algebraically closed field $\7F$ and let $\5Q$ be the associated canonical Frobenius
algebra in $H-\mod$. If $\2E$ is as in Theorem \ref{main0} and $\2B=\HOM_\2E(\6B,\6B)$
then we have the equivalence $\2B\stackrel{\otimes}{\simeq}\hat{H}-\mod$ of spherical
tensor categories.
\etheor

The theorem will be an easy consequence of the more general Theorem \ref{t-depth2}. 
In view of Definition \ref{Morita} we obtain the remarkable

\bcoro \label{c-hopf} 
Let $H$ be a finite dimensional semisimple and cosemisimple Hopf algebra. Then we have
the weak monoidal Morita equivalence $H-\mod\approx\hat{H}-\mod$ of spherical tensor
categories.
\ecoro

We begin with a semisimple spherical $\7F$-linear Morita context $\2E$. (Recall that we
require $\11_\6A$ and $\11_\6B$ to be simple.) We denote $A=\End(J\oj)$, $B=\End(\oj J)$,
$C=\End(J\oj J)$, and we write $Tr_A, Tr_B, Tr_C$ instead of 
$Tr_{J\oj},Tr_{\oj J},Tr_{J\oj J}$. We define a linear map, the `Fourier transform', by
\[ \2F: A\rarr B, \quad \quad
\begin{tangle}
\hh\id\step\id\\
\frabox{a}\\
\hh\id\step\id\\
\object{J}\step\object{\oj}
\end{tangle}
\quad\quad\mapsto\quad
\begin{tangle}
\hh\id\step\id\step\hcoev\mobj{\ol{\ve}(\oj)}\\
\hh\id\step\frabox{a}\step\id\\
\hh\hev\step\id\step\id\\
\hstep\object{\ve(\oj)}\step[1.5]\object{\oj}\step\object{J}
\end{tangle}
\]
and $\hat{\2F}: B\rarr A$ is defined by the same diagram with the obvious changes. The
Fourier transforms are clearly invertible. Furthermore, we define `antipodes' $S: A\rarr
A, \hat{S}: B\rarr B$ by  
$S=\hat{\2F}\circ\2F,\ \hat{S}=\2F\circ\hat{\2F}$. The antipodes are easily seen to be
antimultiplicative: $S(ab)=S(b)S(a)$ and analogously for $\hat{S}$. As a consequence of
axiom (3) in Definition \ref{pivotal}, we have $S\circ S=\id_A$,
$\hat{S}\circ\hat{S}=\id_A$. Using the Fourier transforms we define `convolution products'
on $A$ and $B$ by $a\star b=\2F^{-1}(\2F(a)\2F(b))$ for $a,b\in A$, and similarly for
$B$. One easily verifies 
\[ a\star b= \quad
\begin{tangle}
\id\step\coev\step\id\\
\hh\frabox{a}\Step\frabox{b}\\
\id\step\ev\step\id\\
\object{J}\step[4]\object{\oj}
\end{tangle}
\]
The Fourier transform further allows to define a bilinear form 
$\langle\cdot,\cdot\rangle: A\otimes B\rarr\7F$ by 
$\langle a,b\rangle=d(J)^{-1}\,Tr_A(a\2F^{-1}(b))$. Since $\2F$ is bijective and $Tr_A$ is
non-degenerate, this bilinear form establishes a duality between $A$ and $B$. One verifies
\[ d(J)\,\langle a,b\rangle= Tr_B(\hat{\2F}^{-1}(a)b)=Tr_J \quad
\begin{tangle} 
\hh\hcoev\step\id\\
\hh\id\step\frabox{b}\\
\hh\id\step\id\step\id\\
\hh\frabox{a}\step\id\\
\hh\id\step\hev
\end{tangle}
\quad = \ Tr_{\oj}\quad
\begin{tangle} 
\hh\hcoev\step\id\\
\hh\id\step\frabox{a}\\
\hh\id\step\id\step\id\\
\hh\frabox{b}\step\id\\
\hh\id\step\hev
\end{tangle}
\]
For later use we remark that with $a,b\in A$ we have 
$\langle a,\2F(b)\rangle=d(J)^{-1}\,Tr_A(ab)=d(J)^{-1}\,Tr_A(ba)=\langle b,\2F(a)\rangle$.
The duality between $A$ and $B$ enables us to define coproducts 
$\Delta: A\rarr A\otimes A,\ \hat{\Delta}: B\rarr B\otimes B$ by
\bean \langle \Delta(a),x\otimes y\rangle &=&\langle a,xy\rangle, \quad a\in A, x,y\in B, \\
  \langle a\otimes b,\hat{\Delta}(x) \rangle &=&\langle ab,x\rangle, \quad a,b\in A, x\in B. 
\eean
Associativity of $\hat{m}$ ($m$) implies coassocativity of $\Delta$ ($\hat{\Delta}$). With
\[ \ve(a)=\langle a,1\rangle, \quad \hat{\ve}(x)=\langle 1,x\rangle, \quad a\in A,\ x\in B \]
it is clear that $(A,\Delta,\ve)$ and $(B,\hat{\Delta},\hat{\ve})$ are coalgebras. We note
that $\ve(1)=\langle 1,\hat{1}\rangle=d(J)^{-1}\,Tr_J \id_J=1$, which explains the
normalization of $\langle\cdot,\cdot\rangle$.

The above considerations are valid without further assumptions on the Morita context. In
order to establish $A, B$ as mutually dual Hopf algebras it remains to show that the maps 
$\Delta,\hat{\Delta},\ve,\hat{\ve}$ are multiplicative and that the antipodes are
coinverses. It is here that further assumptions are needed.

\bdefin A semisimple $\7F$-linear Morita context $\2E$ has `depth two' if every simple
summand of $J\oj J\in\Hom(\6B,\6A)$ appears as a simple summand of $J$. $\2E$ is called
irreducible if the distinguished 1-morphism $J: \6B\rarr\6A$ is simple. 
\edefin

If $\2E$ is irreducible and has depth two then $J\oj J$ is a multiple of $J$.
Here we restrict ourselves to the irreducible depth two case, which is all we need to
prove Theorem \ref{t-h-hhat}. Note that we do not assume that $J$ and $\oj$ generate
$\2E$, as is the case in subfactor theory. In a depth two Morita context where this is the
case, every simple $\6B-\6A$-morphism is isomorphic to $J$. For results on the reducible
depth two case -- which leads to finite quantum groupoids -- see \cite{szl2}, where,
however, not all proofs are given.

\begin{figure}
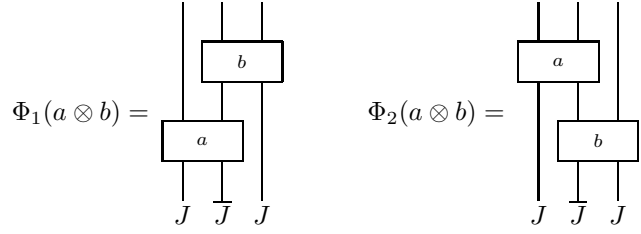

\[ \Phi_1(a\otimes b)= \quad
\begin{tangle}
\hh\id\step\id\step\id\\
\hh\id\step\frabox{b}\\
\hh\id\step\id\step\id\\
\hh\frabox{a}\step\id\\
\hh\id\step\id\step\id\\
\object{J}\step\object{\oj}\step\object{J}
\end{tangle}
\quad\quad\quad\quad
\Phi_2(a\otimes b)= \quad
\begin{tangle}
\hh\id\step\id\step\id\\
\hh\frabox{a}\step\id\\
\hh\id\step\id\step\id\\
\hh\id\step\frabox{b}\\
\hh\id\step\id\step\id\\
\object{J}\step\object{\oj}\step\object{J}
\end{tangle}
\]
\caption{$\Phi_1(a\otimes b)$ and $\Phi_2(a\otimes b)$}
\label{f1}\end{figure}

\blemma \label{l-weyl}
In addition to the above assumptions, let $\2E$ be irreducible of depth two. Then 
\begin{enumerate}
\item The maps 
\bean
  \Phi_1: A\otimes B\rarr C, && a\otimes b\mapsto \id_\oj\otimes b\mcirc a\otimes\id_J, \\
  \Phi_2: A\otimes B\rarr C, && a\otimes b\mapsto a\otimes\id_J\mcirc \id_\oj\otimes b,
\eean
depicted in Figure \ref{f1}, are bijections.
\item For all $a\in A, b\in B$ we have
\begin{equation} \label{e-abc}
   \Phi_1(a\otimes b)=\langle a_{(2)},b_{(1)}\rangle\ \Phi_2(a_{(1)}\otimes b_{(2)}). 
\end{equation}
\end{enumerate}
\elemma
\prf 1. Since $\2E$ is semisimple and $J$ is the only simple $\6B-\6A$-morphism up to
isomorphism, the composition 
$\otimes: \Hom(J\oj J,J)\otimes \Hom(J,J\oj J)\rarr\End(J\oj J)$ is an isomorphism.
Combining this with the isomorphisms $\End(J\oj)\cong\Hom(J\oj J,J)$, etc., provided by
the spherical structure, this easily implies that $\Phi_1,\Phi_2$ are isomorphisms.

2. For $a,e\in A$ we define $x\in B$ by
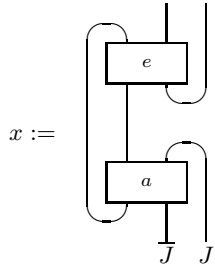
\[ x:=\quad
\begin{tangle}
\hh\hcoev\step\id\step\id\\
\hh\id\step\frabox{e}\step\id\\
\hh\id\step\id\step\hev\\
\hh\id\step\id\step\hcoev\\
\hh\id\step\frabox{a}\step\id\\
\hh\hev\step\id\step\id\\
\Step\object{\oj}\step\object{J}
\end{tangle}
\]
and claim that $x=d(J)^{-1}\,Tr_A(ea_{(1)})\2F(a_{(2)})$. To prove this, we compute
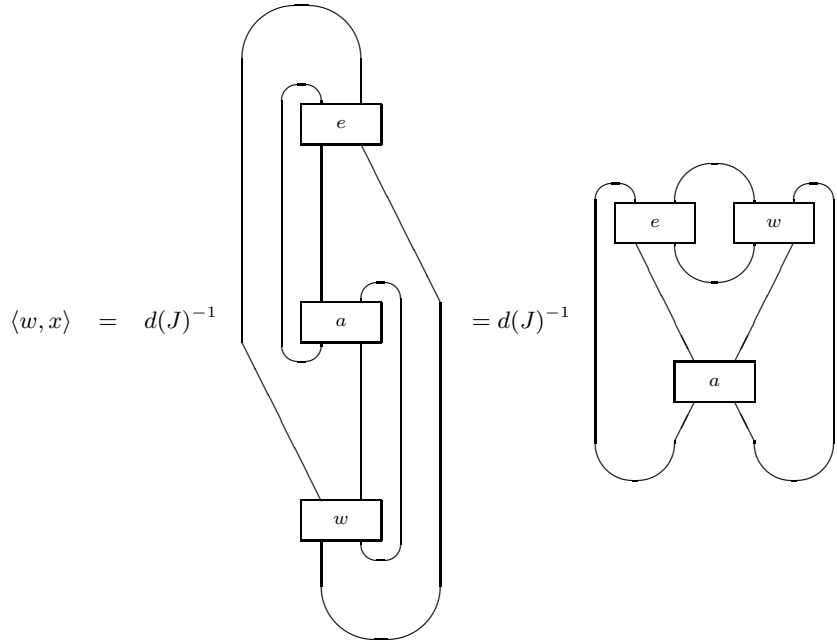
\bean \langle w,x\rangle &=& d(J)^{-1}\quad
\begin{tangle}
\mcoev\\
\hh\id\step\hcoev\step\id\\
\hh\id\step\id\step\frabox{e}\step\\
\id\step\id\step\id\step\d\\
\id\step\id\step\id\step\hcoev\d\\
\hh\id\step\id\step\frabox{a}\step\id\step\id\\
\d\hev\step\id\step\id\step\id\\
\step\d\step\id\step\id\step\id\\
\hh\Step\frabox{w}\step\id\step\id\\
\hh\step\step\id\step\hev\step\id\\
\Step\mev
\end{tangle}
\quad=d(J)^{-1}\quad
\begin{tangle}
\hcoev\step\coev\step\hcoev\\
\hh\id\step\frabox{e}\Step\frabox{w}\step\id\\
\id\step\d\ev\dd\step\id\\
\hh\id\Step\d\step\dd\Step\id\\
\hh\id\step[2.5]\frabox{a}\step[2.5]\id\\
\hh\id\Step\dd\step\d\Step\id\\
\ev\Step\ev
\end{tangle} \\
 &=& d(J)^{-1}\,Tr_A((e\star w)a)=\langle a,\2F(e\star w)\rangle
  =\langle a,\2F(e)\2F(w)\rangle  \\
 &=&\langle a_{(1)},\2F(e)\rangle\langle a_{(2)},\2F(w)\rangle
  =d(J)^{-1}\,Tr_A(ea_{(1)})\langle w,\2F(a_{(2)})\rangle.
\eean
The equality of the two diagrams follows from a simple computation using the axioms of a
spherical category, which we omit. The claim now holds by non-degeneracy of the pairing
$\langle\cdot,\cdot\rangle$. Using our formula for $x$ we find
\bea 
Tr_C \quad 
\begin{tangle}
\hh\id\step\id\step\id\\
\hh\id\step\frabox{f}\\
\hh\id\step\id\step\id\\
\hh\frabox{e}\step\id\\
\hh\id\step\id\step\id\\
\hh\id\step\frabox{b}\\
\hh\id\step\id\step\id\\
\hh\frabox{a}\step\id\\
\hh\id\step\id\step\id\\
\object{J}\step\object{\oj}\step\object{J}
\end{tangle}
\quad &=& Tr_C\quad     
\begin{tangle}
\hh\id\step\id\step\id\\
\hh\frabox{e}\step\id\\
\hh\id\step\id\step\id\\
\hh\id\step\frabox{b}\\
\hh\id\step\id\step\id\\
\hh\frabox{a}\step\id\\
\hh\id\step\id\step\id\\
\hh\id\step\frabox{f}\\
\hh\id\step\id\step\id\\
\object{J}\step\object{\oj}\step\object{J}
\end{tangle}
\quad=\quad
\begin{tangle}
\hcoev\step\coev\step\hcoev\\
\hh\id\step\frabox{x}\Step\frabox{b}\step\id\\
\id\step\d\ev\dd\step\id\\
\hh\id\Step\d\step\dd\Step\id\\
\hh\id\step[2.5]\frabox{f}\step[2.5]\id\\
\hh\id\Step\dd\step\d\Step\id\\
\ev\Step\ev
\end{tangle} \label{e-uvw}\\
&=& d(J)^{-1}\,Tr_A(ea_{(1)}) \quad
\begin{tangle}
\mcoev\\
\hh\id\step\hcoev\step\id\step\hcoev\\
\hh\frabox{a_{(2)}}\step\id\step\frabox{b}\step\id\\
\id\step\d\hev\dd\step\id\\
\hh\id\Step\frabox{f}\Step\id\\
\ev\step\ev
\end{tangle} 
\quad = Tr_A(ea_{(1)})\langle a_{(2)}\hat{\2F}(f),b\rangle \nn\\
 &=& Tr_A(ea_{(1)})\langle a_{(2)},b_{(1)}\rangle \langle\hat{\2F}(f),b_{(2)}\rangle \nn\\
 &=& d(J)^{-1}\,Tr_A(ea_{(1)}) Tr_B(fb_{(2)})\langle a_{(2)},b_{(1)}\rangle. \nn
\eea

Since $J$ is simple, the partial trace $Tr_J(a)\in\End\,\oj$ of $a\in\End\,J\oj$ is given
by $d(J)^{-1}Tr_A(a)\id_\oj$, which implies
\[ \langle a_{(2)},b_{(1)}\rangle \ Tr_C \quad
\begin{tangle}
\hh\id\step\id\step\id\\
\hh\id\step\frabox{f}\\
\hh\id\step\id\step\id\\
\hh\frabox{e}\step\id\\
\hh\id\step\id\step\id\\
\hh\frabox{a_{(1)}}\step\id\\
\hh\id\step\id\step\id\\
\hh\id\step\frabox{b_{(2)}}\\
\hh\id\step\id\step\id\\
\object{J}\step\object{\oj}\step\object{J}
\end{tangle}
\quad =d(J)^{-1}\,\langle a_{(2)},b_{(1)}\rangle Tr_A(ea_{(1)}) Tr_B(fb_{(2)}). 
\]
Comparing this with (\ref{e-uvw}), we see that multiplying both sides of (\ref{e-abc}) with
$\Phi_1(e\otimes f)=\id_\oj\otimes f\mcirc e\otimes\id_J$ and taking traces we obtain an
identity for all $e\in A,f\in B$. In view of 2.\ and the non-degeneracy of the trace in a
spherical category, we conclude that (\ref{e-abc}) holds for all $a\in A,b\in B$.
\qed

\bprop \label{p-depth2}
Let $\2E$ be a Morita context which is semisimple, irreducible and has depth two.
Let $\ve,\hat{\ve},\Delta,\hat{\Delta}$ be defined as above. Then
\begin{enumerate}
\item $\ve, \tilde{\ve}$ are multiplicative. 
\item $\Delta, \tilde{\Delta}$ are multiplicative. 
\item $S,\hat{S}$ are coinverses, i.e.\ 
$m(S\otimes\id)\Delta=m(\id\otimes S)\Delta=\eta\ve$, etc.
\item $A$ and $B$ are semisimple Hopf algebras in duality, and $C$ is the Weyl algebra of
$A$, cf.\ e.g.\ \cite{nill}.
\end{enumerate}
\eprop
\prf 1. Let $a,b\in A$. Since $J$ is simple, $\Hom(\11_\6A,J\oj)$ is one-dimensional and
we have 
\[ a\circ\ve(J)=d(J)^{-1}(\ol{\ve(J)}\circ a\circ\ve(J))\ve(J)
  =\langle a,\hat{1}\rangle\ve(J) =\ve(a)\ve(J). \]
(There should be no danger of confusion between the duality morphisms
$\ve(J), \ol{\ve(J)}$, which are part of the spherical structure, and the counits
$\ve,\hat{\ve}$ of $A$ and $B$.) Thus
$\ve(ab)\ve(J)=(ab)\ve(J)=a(b\ve(J))=a\ve(b)\ve(J)=\ve(a)\ve(b)\ve(J)$, and $\ve$ is
multiplicative.

2. Let $a,b\in A$. Using Lemma \ref{l-weyl} we compute 
\bea
\begin{tangle} 
\hh\id\step\id\step\id\\
\hh\id\step\frabox{b}\\
\hh\id\step\id\step\id\\
\hh\frabox{a}\step\id\\
\hh\id\step\hev
\end{tangle}
& = &\langle a_{(2)},b_{(1)}\rangle \quad
\begin{tangle} 
\hh\id\step\id\step\id\\
\hh\frabox{a_{(1)}}\step\id\\
\hh\id\step\id\step\id\\
\hh\id\step\frabox{b_{(2)}}\\
\hh\id\step\hev
\end{tangle}
\quad =
\langle a_{(2)},b_{(1)}\rangle \hat{\ve}(b_{(2)})\quad
\begin{tangle} 
\hh\id\step\id\step\id\\
\hh\frabox{a_{(1)}}\step\id\\
\hh\id\step\hev
\end{tangle} \nn \\
& =&
\langle a_{(2)},b\rangle \quad
\begin{tangle} 
\hh\id\step\id\step\id\\
\hh\frabox{a_{(1)}}\step\id\\
\hh\id\step\hev
\end{tangle}
\label{e1} \eea
In an entirely analogous fashion one shows
\begin{equation} \label{e2}
\begin{tangle} 
\hh\hcoev\step\id\\
\hh\id\step\frabox{b}\\
\hh\id\step\id\step\id\\
\hh\frabox{a}\step\id\\
\hh\id\step\id\step\id\\
\end{tangle}
\quad\quad = \langle a,b_{(1)}\rangle \quad
\begin{tangle} 
\hh\hcoev\step\id\\
\hh\id\step\frabox{b_{(2)}}\\
\hh\id\step\id\step\id\\
\end{tangle}
\end{equation}
Let now $a,b\in A, c,d\in B$. We compute
\bean
\langle \Delta(ab),c\otimes d\rangle=
\langle ab,cd\rangle &=& d(J)^{-1}\,Tr\quad
\begin{tangle}
\hh\hcoev\step\id\\
\hh\id\step\frabox{c}\\
\hh\id\step\id\step\id\\
\hh\id\step\frabox{d}\\
\hh\id\step\id\step\id\\
\hh\frabox{a}\step\id\\
\hh\id\step\id\step\id\\
\hh\frabox{b}\step\id\\
\hh\id\step\hev
\end{tangle}
\quad=d(J)^{-1}\,\langle a_{(2)},d_{(1)}\rangle \ Tr_J \quad
\begin{tangle}
\hh\hcoev\step\id\\
\hh\id\step\frabox{c}\\
\hh\id\step\id\step\id\\
\hh\frabox{a_{(1)}}\step\id\\
\hh\id\step\id\step\id\\
\hh\id\step\frabox{d_{(2)}}\\
\hh\id\step\id\step\id\\
\hh\frabox{b}\step\id\\
\hh\id\step\hev
\end{tangle} \\
&=& d(J)^{-1}\,\langle a_{(2)},d_{(1)}\rangle \langle a_{(1)},c_{(1)}\rangle 
    \langle b_{(2)},d_{(2)}\rangle\ Tr_J \quad 
\begin{tangle}
\hh\hcoev\step\id\\
\hh\id\step\frabox{c_{(2)}}\\
\hh\id\step\id\step\id\\
\hh\frabox{b_{(1)}}\step\id\\
\hh\id\step\hev
\end{tangle} \\
 &=& \langle a_{(2)},d_{(1)}\rangle \langle a_{(1)},c_{(1)}\rangle 
  \langle b_{(2)},d_{(2)}\rangle \langle b_{(1)},c_{(2)}\rangle  \\
 &=& \langle a_{(1)}b_{(1)},c\rangle\langle a_{(2)}b_{(2)},d\rangle 
  =\langle \Delta(a)\Delta(b),c\otimes d\rangle.
\eean
(The first and sixth equality hold by definition of $\Delta$ and $\hat{\Delta}$,
respectively. The second and fifth are just the definition of
$\langle\cdot,\cdot\rangle$. The third equality follows by Lemma \ref{l-weyl} and the
fourth is due to (\ref{e1}) and (\ref{e2}).) Since this holds for all $c,d$, we
conclude by duality that $\Delta(ab)=\Delta(a)\Delta(b)$, as desired.

3. Appealing once more to Lemma \ref{l-weyl}, we have
\[ Tr_J \quad
\begin{tangle}
\hh\id\step\hcoev\\
\hh\id\step\frabox{b}\\
\hh\id\step\id\step\id\\
\hh\frabox{a}\step\id\\
\hh\hev\step\id
\end{tangle}
\quad = \langle a_{(2)},b_{(1)}\rangle\ Tr_J \quad
\begin{tangle}
\hh\id\step\hcoev\\
\hh\frabox{a_{(1)}}\step\id\\
\hh\id\step\id\step\id\\
\hh\id\step\frabox{b_{(2)}}\\
\hh\hev\step\id
\end{tangle}
\quad\quad \forall a\in A, b\in B.
\]
The left hand side equals 
\[ \ve(a)\, Tr_J \quad
\begin{tangle}
\hh\id\step\hcoev\\
\hh\id\step\frabox{b}\\
\hh\hev\step\id\\
\end{tangle}
\quad = d(J) \ve(a)\hat{\ve}(b).
\]
The right hand side equaling 
$d(J)\langle a_{(2)},b_{(1)}\rangle\langle a_{(1)},S(b_{(2)})\rangle$, we obtain
\[ \ve(a)\hat{\ve}(b)=\langle a_{(2)},b_{(1)}\rangle\langle a_{(1)},\hat{S}(b_{(2)})\rangle 
  =\langle a,\hat{S}(b_{(2)})b_{(1)}\rangle, \]
which is equivalent to $\hat{S}(b_{(2)})b_{(1)}=\hat{\ve}(b)\hat{1}$. Since $\hat{S}$ is
involutive, we find $\hat{S}(b_{(1)})b_{(2)}=\hat{\ve}(b)\hat{1}$. The other identities
are proved similarly.

4. The first statement just summarizes our results so far, and the second is obvious by 
definition of the Weyl algebra \cite{nill}.
\qed

Given a semisimple irreducible depth two Morita context, the preceding theorem provides us with 
a pair $A,B$ of mutually dual Hopf algebras. It remains to relate these Hopf algebras to the
categorical structure. Here we use Tannaka theory for Hopf algebras together with a result
form \cite{mue-t}. The strategy is (i) to construct a faithful tensor functor 
$E: \2A\rarr\mbox{Vect}_\7F$, (ii) to deduce that $\2A$ is monoidally equivalent to $H-\mod$
for some Hopf algebra $H$ and (iii) to prove that $H\cong A$.

\btheor \label{t-depth2}
Let $\2E$ be a semisimple spherical irreducible depth two Morita context. Consider the
full tensor subcategories 
\[ \2A_0\subset\2A=\END_\2E(\6A), \quad\quad \2B_0\subset\2B=\END_\2E(\6B) \]
consisting of the tensor powers of $J\oj$ and $\oj J$, respectively, and their retracts. 
Then we have the equivalences 
\[ \2A_0\stackrel{\otimes}{\simeq} A^\cop-\mod, \quad\quad
   \2B_0\stackrel{\otimes}{\simeq} B^\cop-\mod \]
of spherical tensor categories, where $A, B$ are the Hopf algebras constructed in
Proposition \ref{p-depth2}.
\etheor
\prf Let $Q=J\oj$ and consider the Frobenius algebra $(Q,v,v',w,w')$ with 
$v=\ve(J),\ w'=\id_J\otimes\ol{\ve}(\oj)\otimes\id_\oj$ etc. In particular, 
$(Q,w',v)$ is a monoid in $\2A_0$. By the depth 2 property, we have we have 
$\oj X\cong d(X)\oj$ for all $X\in\2A_0$, i.e.\ there are morphisms $r_i: \oj\rarr\oj X$,
$r'_i: \oj X\rarr\oj$ satisfying the usual conditions. Thus the morphisms 
$s_i=\id_J\otimes r_i: Q=J\oj\rarr J\oj X=QX$ establish an isomorphism $QX\cong d(X)Q$. One
easily verifies that the $s_i$ are $Q$-module morphisms, i.e.\ satisfy
$s_i\circ w'=w'\otimes\id_X\mcirc id_Q\otimes s_i$, and similarly the $s'_i=\id_J\otimes r'_i$.
These facts imply that the functor $E: \2A_0\rarr\mbox{Vect}_\7F$ defined by 
$E(X)=\Hom(\11,Q\otimes X)$ and $E(s)\phi=(\id_Q\otimes s)\circ\phi$ for $s: X\rarr Y$ is
a faithful (strong) tensor functor, where the isomorphisms $d_{X,Y}: E(X)\boxtimes
E(Y)\rarr E(X\otimes Y)$ are given by 
$d_{X,Y}(\phi\boxtimes\psi)=w'\otimes\id_{X\otimes Y}\mcirc\id_Q\otimes\phi\otimes\id_Y\mcirc\psi$
for $\phi\in E(X), \psi\in E(Y)$. See \cite[Sect.\ 3]{mue-t} for the details. It then
follows from Tannaka theory there exists a finite dimensional Hopf algebra $H$ and an
equivalence $F: \2A_0\rarr H-\mod$ such that $K\circ F=E$, where $K: H-\mod\rarr\7F$-Vect
is the forgetful functor. Here $H=\mbox{Nat}\,E$ is the $\7F$-algebra of natural
transformations from $E$ to itself. Consider the map $\alpha$ which to $a\in\End\,Q=A$
associates the family $\alpha(a)=\{ \alpha(a)_X\in \End\,E(X),\ X\in\2A_0\}$, where
$\alpha(a)_X\phi=(a\otimes\id_X)\mcirc\phi$ for $\phi\in E(X)$. It is obvious that
$\alpha(a)\in\mbox{Nat}\,E=:H$ and that $\alpha: A\rarr H$ is an algebra
homomorphism. Semisimplicity of $\2A_0$ and finite dimensionality of $H$ imply that
$\alpha$ is an isomorphism, which we now suppress. It remains to show that the coproduct
$\Delta'$ of $H$ provided by Tannaka theory coincides with the one constructed in
Proposition \ref{p-depth2}. By definition of $\Delta'(a)=a_{(1)}\otimes a_{(2)}$, the diagram
\[ \begin{diagram}
  E(X)\otimes E(Y) & \rTo^{d_{X,Y}} & E(X\otimes Y) \\
  \dTo^{a_{(1),X}\otimes a_{(2),Y}} & & \dTo_{a_{X\otimes Y}} \\
  E(X)\otimes E(Y) & \rTo_{d_{X,Y}} & E(X\otimes Y)
\end{diagram} \]
commutes for all $a\in A$ and $X,Y\in\2A_0$. In view of the definition of $d_{X,Y}$, this
is equivalent to $a\mcirc w'=w'\mcirc a_{(2)}\otimes a_{(1)}$ for all $a\in A=\End\,Q$.
For arbitrary $b,c\in A$, this implies
\[ Tr_Q (a\mcirc w'\mcirc b\otimes c\mcirc w)=Tr_Q(w\mcirc a_{(2)}\otimes a_{(1)}\mcirc
   b\otimes c\mcirc w). \]
Consistent with previous terminology we write $b\star c=w'\mcirc b\otimes c\mcirc w\in A$
for $b,c\in A$, and the fact that $Q$ contains $1$ with multiplicity one implies
$Tr_Q(b\star c)=d(J)^{-1}Tr_Q(b)Tr_Q(c)$. Therefore, 
\[ Tr_Q(a(b\star c))=d(J)^{-1}Tr_Q(a_{(2)}b)Tr_Q(a_{(1)}c)\quad\mbox{where}\quad
   a_{(1)}\otimes a_{(2)}=\Delta'(a). \]
On the other hand, the definition $\langle\Delta(a),b\otimes c\rangle=\langle a,bc\rangle$
of $\Delta$ as given above satisfies
\[ Tr_Q(a(b\star c))=d(J)^{-1}Tr_Q(a_{(1)}b)Tr_Q(a_{(2)}c) \quad\mbox{where}\quad
   a_{(1)}\otimes a_{(2)}=\Delta(a). \]
Thus $\Delta'=\Delta^\cop$, and we are done.
\qed

We briefly recall some facts concerning the (left) regular representation $Q_l\in A-\mod$
of a semisimple Hopf algebra $A$. We have $Q_l\cong\oplus_X d(X)X$, where the $X$ are the
irreducible representations, and therefore $\dim\Hom(X,Q_l)=d(X)$ for all simple
$X$. Furthermore, the regular representation is absorbing: $X\otimes Q_l\cong Q_l\otimes
X\cong d(X)Q_l$ for every $X\in A-\mod$. \\

\noindent {\it Proof of Theorem \ref{t-h-hhat}.} 
By \cite{bw2}, the category $H-\mod$ is spherical and by the coherence theorem
\cite{bw2} we may consider $H-\mod$ as strict monoidal and strict spherical. By Theorem
\ref{t-hopffrob} we have a canonical and irreducible Frobenius algebra $\5Q$ in $H-\mod$,
which we can normalize such that $\lambda_1=\lambda_2$. Since $\5Q$ is irreducible, by
Proposition \ref{p-simple} the same is true for the Morita context $\2E$ of Theorem
\ref{main0}. By Theorem \ref{t-spherical2} there thus is a spherical structure on $\2E$ 
extending that of $\2A$. The claim now follows from Theorem \ref{t-depth2} and the fact
$H^{\op,\cop}\cong H$, provided we can show that $\2E$ has depth 2. 

By definition of $\2E$, every $Y:\6B\rarr\6A$ is a retract of $XJ$ for some
$X\in\END(\6A)\simeq\2A$. By semisimplicity is thus sufficient to show that $XJ$ is a
multiple of $J$ for every simple $X\in\2A$. We have
\bean \Hom(J,XJ) &\cong& \Hom(J\oj,X)\cong\Hom_\2A(Q,X), \\
   \End(XJ) &\cong& \Hom(XJ\oj,X)\cong\Hom_\2A(XQ,X).
\eean
By the properties of $Q$ recalled above, we have $\dim\Hom(J,XJ)=d(X)$, implying that
$XJ$ contains $J$ with multiplicity $d(X)$, thus $\End(XJ)$ contains the matrix algebra
$M_{d(X)}(\7F)$. In view of $\dim\End(XJ)=\dim\Hom(XQ,X)=\dim\Hom(d(X)Q,X)=d(X)^2$ we
conclude $\End(XJ)\cong M_{d(X)}(\7F)$ and therefore $XJ\cong d(X)J$ as desired.
\qed

\brem
If $\2E$ in Proposition \ref{p-depth2} is a $*$-bicategory then $A, B$ come with canonical
$*$-operations. It is then not difficult to show that $\ve, \Delta, S$ are
$*$-homomorphisms, thus $A$ and $B$ are Hopf $*$-algebras. (E.g., the property 
$\ol{\ve(J)}=\ve(J)^*$ immediately implies $\ol{\ve(a)}=\ve(a^*)$.) In the Theorems
\ref{t-depth2} and \ref{t-h-hhat} we then have equivalences of tensor $*$-categories. We
omit the proofs. 
\erem


\subsection{Subfactors} \label{ss-subfact}
The entire analysis of this paper is motivated by the mathematical structures which are
implicit in subfactor theory. In this subsection we make the link between subfactor
theory and our categorical setting explicit, shedding light on both subjects. The
main aim of this section is in fact to improve the communication between subfactor
theorists and category minded people, the only new result being Theorem \ref{main3}. We
begin with a very brief definition of the notions we will use. For everything else see any
textbook on von Neumann algebras, e.g., \cite{sun,takes,stzs} and subfactors 
\cite{vfr-su,ek2}. 

A von Neumann algebra (vNa) is a unital subalgebra $M\subset\2B(\2H)$ of the algebra
$\2B(\2H)$ of bounded operators on some Hilbert space $\2H$ which is closed w.r.t.\ the
hermitian conjugation $x\mapsto x^*$ and w.r.t.\ weak convergence. Equivalently, by von
Neumann's double commutant theorem a vNa is a set $M\subset\2B(\2H)$ which is closed under
conjugation and satisfies $M''=M$, where $S'=\{ x\in\2B(\2H)\ | \ xy=yx \ \forall y\in S\}$ 
is the commutant of $S$. A factor is a vNa with trivial center ($M\cap M'=\7C\11$) and if
$N, M$ are factors such $N\subset M$ then $N$ is called a subfactor. (By abuse of notation
`subfactor' occasionally refers to the inclusion $N\subset M$.) Every factor $M$ is of one
of the types I, II or III, where $M$ is of type I iff $M\simeq \2B(\2K)$ for some
Hilbert space $\2K$. (Every finite dimensional factor is of type I.)
If $N\subset M$ are both of type I then also $M\cap N'$ is of type I and 
$M\simeq N\ol{\otimes} (M\cap N')$. Under this isomorphism the embedding $N\hookrightarrow
M$ becomes $x\mapsto x\ol{\otimes} \11$, and nothing more of interest is to be said. In
our discussion of the remaining cases we restrict ourselves to vNas on a separable Hilbert
space which simplifies the definitions. Then a factor $M$ is of type III iff every
orthogonal projection $e=e^2=e^*\in M$ is the range of some isometry $v\in M$, i.e.\
$vv^*=e,\ v^*v=\11$. A factor $M$ which is neither type I or type III is of type II,
of which there are two subclasses: II$_1$ and II$_\infty$. A factor is of type II$_1$ iff
it admits a tracial state (trace, for short) $tr: M\rarr\7C$, i.e.\ a weakly continuous
linear functional which is positive ($A>0\impl tr\,A>0$), normalized ($tr(\11)=1$) and
vanishes on commutators. (It follows that every isometry in a type II$_1$ factor must be
unitary.) A type II$_\infty$ factor is isomorphic to the tensor product of some II$_1$
factor with $\2B(\2H)$ where $\dim \2H=\aleph_0$. The tensor product of any factor with a
type III factor is of type III. 

In the literature on subfactors the focus has been on type II$_1$ factors, which are
technically easiest to deal with thanks to the existence of a trace, cf.\ the textbooks
\cite{vfr-su,ek2}. Yet, in our discussion we concentrate on the type III case, the
technical aspects of which have been clarified in particular in the work of Longo
\cite{lo1, lo3}. Those aspects of subfactor theory \cite{lo3, iz2} which most directly
inspired the present investigation and \cite{mue10} were in fact done in the type III
setting. Anyway, by Popa's results \cite{popa} the classifications of amenable inclusions
of hyperfinite type II$_1$ and III$_1$ factors amount to the same thing.

The following is implicit in much of the literature on type III subfactors and explicit
in \cite[Section 7]{lro}. 
\bdefin We denote by $\2T$ the 2-category whose objects are type III factors with
separable predual. The 1-morphisms are normal unital $*$-homomorphisms with the obvious 
composition. For parallel 1-morphisms $\rho, \sigma: M\rarr N$ the 2-morphisms are given by
\[ \Hom_\2T(\rho,\sigma)= \{ s\in N\ | \ s\rho(x)=\sigma(x)s\ \forall x\in M\}. \]
The vertical composition of 2-morphisms is multiplication in $N$ and the horizontal
composite 
\[
M\ctwo{\rho_1}{\sigma_1}{s}%
N\ctwo{\rho_2}{\sigma_2}{t}%
O
\]
is given by 
\[ s\times t= t\rho_2(s)=\sigma_2(s)t: \rho_2\rho_1\rarr\sigma_2\sigma_1. \]
\edefin

\blemma The 2-category $\2T$ has direct sums of 1-morphisms idempotent 2-morphisms
split. All identity 1-morphisms are simple.
\elemma
\prf Let $\rho,\sigma: M\rarr N$. Pick an orthogonal projection $e\in N$, put $f=1-e$ and
choose isometries $p,q\in N$ such that $pp^*=e, qq^*=f$. Then 
\[ \eta(\cdot)=p\rho(\cdot)p^*+q\sigma(\cdot)q^* \]
is a direct sum. Let $\rho: M\rarr N$ and $e=e^2=e^*\in\End_\2T(\rho)\subset N$. Picking
an isometry $q\in N$ such that $qq^*=e$ and setting 
\[ \sigma(\cdot)=q\rho(\cdot)q^* \]
one obviously has $q\in\Hom_\2T(\sigma,\rho)$. The last claim follows from
factoriality. 
\qed

\brem This lemma fails for finite factors, which is why one works with bimodules in
the type II$_1$ case.
\erem

Since the kernel of a normal $*$-homomorphism is a closed two-sided ideal and a type III
factor with separable predual is simple, all 1-morphisms in $\2T$ are injective. Thus a
morphism $\rho: N\rarr M$ provides an isomorphism of $N$ with the subalgebra
$\rho(N)\subset M$. Now Longo's main result in \cite{lo1} can be rephrased as  follows:

\bprop Let $\rho: N\rarr M$ be a 1-morphism in $\2T$. The $\rho$ has a two-sided adjoint
$\ol{\rho}: M\rarr N$ in $\2T$ iff the index $[M:\rho(N)]$ is finite. 
\eprop
\prf We identify $N$ with the subalgebra $\rho(N)\cong N$ of $M$, but we insist on writing
an embedding map $\iota: N\hookrightarrow M$. Longo proves the following: the index
$[M:N]$ of $N\subset M$ is finite iff there is a triple $(\gamma, v, w)$  where $\gamma$
is a normal $*$-endomorphism of $M$ with that $\rho(M)\subset N$ such that there are
isometries $v\in M, w\in N$ satisfying
\begin{equation} \label{p1}
   v\in\Hom_\2T(\id_M,\gamma),\quad w\in\Hom_\2T(\id_N,\gamma\restr N), 
\end{equation}
\begin{equation} \label{p2}
   v^*w=\gamma(v^*)w=c\11. 
\end{equation}
Then $c=[M:N]^{-1/2}$, and we refer to \cite{lo1} for the original definition of the
index $[M:N]$. In order to translate Longo's result into categorical language we write
$\ol{\iota}=\gamma$ and observe that $\gamma$ maps $M$ into $N$, so that only
$\iota\circ\ol{\iota}$ gives an morphism of $M$. On the other hand we see 
$\gamma\restr N=\ol{\iota}\circ\iota$. Now (\ref{p1}) becomes 
\[ v\in\Hom_\2T(\id_M,\iota\ol{\iota}),\quad w\in\Hom_\2T(\id_N,\ol{\iota}\iota) \]
and $v^*, w^*$ are morphisms in the opposite directions. Finally in view of the definition
of the horizontal composition of 2-morphisms in $\2T$ the equations (\ref{p2}) and their
$*$-conjugates are -- up to a numerical factor which can be absorbed in $v$ or $w$ -- the
four triangular equations which make $\iota, \ol{\iota}$ two-sided duals. If we normalize
$v, w$ such that $v^*v=w^*w=[M:N]^{1/2}$ then $d(\iota)=d(\ol{\iota})=[M:N]^{1/2}$. 
\qed

\blemma Every inclusion of type III factors with finite index defines a  2-*-category
$\2T_{N\subset M}$ which is a Morita context. The dimension of $\2F$ (which is
well-defined by Proposition \ref{eqofdims}) is finite iff the subfactor has finite depth. 
\label{subf}\elemma 
\prf $\2T_{N\subset M}$ is the subcategory of $\2T$ whose objects are $\{N, M\}$ and whose
1-morphisms are generated by $\iota, \ol{\iota}$. More precisely, $\Hom_{\2T_{N\subset
M}}(N,M)$ is the replete full subcategory of $\Hom_\2T(N,M)$ whose morphisms are retracts
of some $\iota(\ol{\iota}\iota)^N,\ N\in \7N\cup\{ 0\}$, and similarly for the other
categories of 1-morphisms. In this category the functors $-\otimes\iota$ etc. are clearly
dominant, thus $\2T_{N\subset M}$ is a Morita context for the tensor categories
$\End_{\2T_{N\subset M}}(N), \End_{\2T_{N\subset M}}(M)$. The $*$-involution obviously it
the $*$-operation of the algebras. The last claim follows from Proposition \ref{eqofdims} and
the definition according to which $N\subset M$ has finite depth iff the powers of
$\iota\ol{\iota}$ contain finitely many simple $M-M$ morphisms up to equivalence.
\qed 

\brem 1. In subfactor theory the dimension of $\2F$ is called the global index
(as opposed to the index $[M:N]=d(Q)$). 

2. $\2T_{N\subset M}$ having a $*$-structure it can be made into a spherical category,
though not in a completely unique way. See Section \ref{stars}. 
\erem

Now, given an endomorphism $\gamma$ of finite index of a type III factor $M$ it is
natural to ask whether there is a subfactor $N\subset M$ such that
$\gamma=\iota\ol{\iota}$, where $\ol{\iota}$ is a two-sided conjugate of the embedding
morphism $\iota$. The answer, given in \cite{lo3}, is positive iff there are isometries 
\[ v\in\Hom(\id_M,\gamma),\quad w\in\Hom(\gamma,\gamma^2) \]
satisfying
\bea w^2 &=& \gamma(w)w, \label{L1}\\
   ww^* &=& \gamma(w^*)w, \label{L2}\\
   v^*w &=& \gamma(v)^*w=c\11,\ c\in\7C^*. \label{L3}\eea
(It turned out that (\ref{L2}) is redundant, cf.\ \cite{lro}.) Then the subfactor is given
by  
\begin{equation} N= w^* \gamma(M) w. \label{n}\end{equation}

The equations (\ref{L1}-\ref{L3}) together with the requirement that $v, w$ are isometries
are, of course, saying precisely that $(\gamma, v, v^*, w, w^*)$ is a canonical Frobenius
algebra in $\End M$. In a sense,  this entire paper is about finding a categorical
analogue for the simple formula (\ref{n}), which turned out to be more tedious than one
might expect. This is precisely due to the fact that as seen above subfactor theory comes
with a rich and beautiful inherent categorical structure which we had to model by Theorem
\ref{main0}. The reward for our work is the following result. 

\btheor Let $M$ be a type III factor, let $\gamma\in\End(M)$ satisfy (\ref{L1}-\ref{L3})
and let $\2A$ be the replete full subcategory with subobjects of $\End(M)$ generated by
$\gamma$. Let $\2T_{N\subset M}$ be the bicategory associated with the subfactor $N\subset
M$, where $N$ is given by (\ref{n}). (Obviously, $\2A=\HOM_\2T(\6A,\6A)$.) If $\2E$ is
obtained from $(\2A, \gamma)$ by Theorem \ref{main0} then have an equivalence of
bicategories $\2E\simeq\2T_{N\subset M}$. 
\label{main3}\etheor 
\prf By Lemma \ref{subf} the 2-category $\2T_{N\subset M}$ is a Morita context. Thus the
claim follows directly from Proposition \ref{main2}. 
\qed

\brem 1.
The importance of this theorem is that it allows us to compute the 2-category
$\2T_{N\subset M}$ (up to equivalence) from the data $(\2A, \gamma)$ without explicitly
working with subfactors. In the case where $\2A$ is the subcategory generated by a
canonical object $\gamma$ in some $\End M$ this may seem a rather complicated detour. Yet,
we have gained two things. One one hand we see that the bicategory associated with a
subfactor with finite index is a structure which appears also in other contexts. Equally 
important is the fact that our constructions work for arbitrary tensor categories
$\2A$ which are not subcategories of some $\End M$ generated by one object $\gamma$, in
fact for arbitrary (algebraically closed) ground field. This will be exploited in
\cite{mue10}, the results of which seem hard to prove without our machinery. 

2. The results of Section \ref{ss-Hopf2} are related to subfactor theory in a very direct way
whenever the Hopf algebra $H$ is a finite dimensional $C^*$-Hopf algebra. (This means $H$
is a complex multi-matrix algebra and $\Delta, \ve$ respect the natural $*$-operation.)
As shown in \cite{yaman} every finite dimensional $C^*$-Hopf algebra $H$ admits an action
on a type II$_1$ factor $M$. This action is outer, i.e. $(M^H)'\cap M=\7C\11$. It has long
been known \cite{ocn1} that in this situation the Jones extension $M_1$ of the subfactor
$M^H\subset M$ carries an outer action of $\hat{H}$. Together with the well-known material
in the present section this provides a von Neumann algebraic proof of the weak Morita
equivalence $H-\mod\approx\hat{H}-\mod$. From the perspective of this paper this proof is,
however, quite unsatisfactory. On the one hand it is restricted to $C^*$-Hopf algebras,
on the other it is rather indirect since it involves infinite dimensional von Neumann
algebras. 
\erem


\section{Morita Invariance of State Sum Invariants} \label{s-topology}
In this section we give an interesting and non-trivial application of our notion of weak
monoidal Morita equivalence to the study of triangulation (or state sum) invariants of
closed 3-manifolds. We begin with a very brief sketch of the works which are relevant to
our discussion, apologizing to everyone whose contribution is being glossed over. 

In \cite{tv} Turaev and Viro used the 6j-symbols of the quantum group $SU_q(2)$ to define
a numerical invariant $TV_{SU_q(2)}(M,T)$ for any closed 3-manifold together with a
triangulation $T$. They went on to prove that it does not depend on $T$ and thus gives
rise to a topological invariant $TV_{SU_q(2)}(M)$. In \cite{t} this construction was
generalized to an invariant $TV(M,\2C)$ associated with any modular category
$\2C$. (Modular categories are braided ribbon categories satisfying a certain
non-degeneracy condition.) Recently S.\ Gelfand and Kazhdan \cite{gk} and Barrett and
Westbury \cite{bw2} defined triangulation invariants for 3-manifolds on the basis of
certain tensor categories which are not required to be braided. In our discussion we focus
on the invariant of Barrett and Westbury, which we call $BW(M,\2C)$, since it is based on
spherical categories and therefore close in spirit to our work. With appropriate
normalizations one has $TV(M,\2C)=BW(M,\2C)$ if $\2C$ is modular. The utility of the
notion of weak monoidal Morita equivalence is now illustrated by the following

\btheor \label{appl}
Let $\2A, \2B$ be (strict) semisimple spherical categories with simple unit and finitely
many simple objects. If $\2A, \2B$ are weakly monoidally Morita equivalent and 
$\dim\2A\ne 0$ then we have $BW(M,\2A)=BW(M,\2B)$ for all closed orientable 3-manifolds $M$. 
\etheor

\brem \label{ex1}
1. Before we sketch the proof of this result we point out that it resolves a (minor) puzzle
concerning the BW invariant. In \cite{ku} Kuperberg had defined a 3-manifold invariant
$Ku(M, H)$ for every finite dimensional Hopf algebra $H$ over an algebraically closed
field $\7F$ which is involutive ($S^2=\id$) and whose characteristic does not divide the
dimension of $H$. (That these conditions are equivalent to semisimplicity of $H$ and
$\hat{H}$ was not known then.) In \cite{bw3} it was proven that the invariant $BW$ is a
generalization of $Ku$ in the sense that $Ku(M,H)=BW(M,H-\mod)$, again assuming
appropriate normalizations. For Kuperberg's invariant it had been known that
$Ku(M,H)=Ku(M,\hat{H})$, but this becomes obscure when it is expressed in terms of the
invariant BW. This puzzle is resolved by Theorem \ref{appl} together with Corollary
\ref{c-hopf}, according to which $H-\mod$ and $\hat{H}-\mod$ are weakly monoidally Morita
equivalent.

2. The preceding application of weak monoidal Morita equivalence is not really new in that
the result $BW(M,H-\mod)=BW(M,\hat{H}-\mod)$ can be derived from the connection between
the invariants $BW$ and $Ku$. A less obvious example is provided in the companion paper
\cite{mue10}. There we prove that the center $\2Z(\2C)$ \cite{ma1,js2, str2}, which is the
categorical version of Drinfel'd's quantum double, of a semisimple spherical category with
non-zero dimension is again spherical and semisimple (and modular in the sense of
Turaev). Furthermore, we prove the weak monoidal Morita equivalence
$\2Z(\2C)\approx\2C\boxtimes\2C^\op$, which by the above theorem implies
\bean \lefteqn{ BW(M,\2Z(\2C))=BW(M,\2C\boxtimes\2C^\op)=BW(M,\2C)\cdot BW(M,\2C^\op)} \\
  && = BW(M,\2C)\cdot BW(-M,\2C) \quad
   (=|BW(M,\2C)|^2 \ \ \mbox{if $\2C$ is a $*$-category}).
\eean
The relation between a category $\2C$ and its quantum double being quite non-trivial we
are not aware of a simpler proof of this equality.
\erem

\noindent{\it Sketch of Proof.} The proof relies strongly on ideas of Ocneanu which,
unfortunately, found expression only in the unpublished (and unfinished) manuscript 
\cite{ocn3}. Therefore the more complete accounts \cite{ek3,kosu} are very useful. In
\cite{ocn3} Ocneanu defined a triangulation invariant $Oc(M, A\subset B)$ of 3-manifolds
departing from an inclusion $A\subset B$ of type II$_1$ factors with finite index and
finite depth. We recall from Section \ref{ss-subfact} that subfactors $A\subset B$ with
finite index give rise to a Morita context $\2E$, whose dimension is finite iff the
subfactor has finite depth. ($\2E$ is given by bimodules associated with the subfactor or, 
alternatively, by $*$-algebra homomorphisms in the type III case.) Ocneanu's invariant
is easily seen to depend only on $\2E$ and not on other, in particular analytic properties
of the subfactor. Furthermore, as is quite evident from \cite{ek3}, it generalizes to any 
spherical (or $*$-) Morita context of finite, non-zero dimension over an algebraically
closed field. As to the definition of the invariant we only say that one chooses a 
triangulation $T$ with directed edges and an assignment $\5V\in\{\6A,\6B\}^V$ of labels
$\{\6A,\6B\}$ to the vertices $V$. Then to every edge of $T$ one assigns an isomorphism
class of 1-morphisms in $\HOM_\2E(\6X,\6Y)$, where $\6X, \6Y$ are the labels attached to
the initial and terminal vertices of the edge. $Oc(M,\2E,T,\5V)$ is now defined as the
sum over the edge labelings of a product of 6j-symbols. Note that there is no
summation over the labeling $\5V$ of the vertices! In fact it is shown in \cite{ocn3,ek3}
that $Oc(M,\2E,T,\5V)$ depends neither on the labeling $\5V$ (for fixed triangulation
$T$) nor on $T$. If one labels all vertices of $T$ with $\6A$ one finds that
$Oc(M,\2E,T)=BW(M,\END_\2E(\6A),T)$, i.e.\ the invariant reduces to the invariant of
Barrett and Westbury for the spherical category $\END_\2E(\6A)$. Similarly, by labeling
all vertices with $\6B$ one obtains $Oc(M,\2E,T)=BW(M,\END_\2E(\6A),T)$. By independence
of $BW$ and $Oc$ of the triangulation one concludes 
\begin{equation} Oc(M,\2E)=BW(M,\END_\2E(\6A))=BW(M,\END_\2E(\6B)). \label{mi}
\end{equation} 
Our claim thus follows if we take $\2E$ to be a Morita context for $\2A\approx\2B$.
\qed

\brem 1. The argument sketched above should of course be spelled out in more detail. In
particular, this requires a construction of the TQFT associated with the invariant 
$BW(\cdot,\2A)$. (In doing so extreme care is required when the tensor category $\2A$
contains simple objects which are self-dual and pseudo-real. Unfortunately, this is
neglected in the bulk of the literature on the subject, with the notable exception of
\cite{t} and the remarks in \cite{bw1}.) We hope to do this in a future part of this
series. 

2. Let us emphasize the lesson we draw from the above considerations. In view of
(\ref{mi}) the invariant $Oc(M,\2E)$ is already determined by considerably smaller amounts
of data, as contained in the tensor categories $\END_\2E(\6A)$ or $\END_\2E(\6B)$. 
(Therefore the observation \cite{bw1,gk} that the Turaev-Viro invariant generalizes to
tensor categories without braiding could have been made already by the authors of
\cite{ocn3, ek3}.) Despite the greater generality of \cite{bw1,gk} there is a lasting
significance of the invariant $Oc$ which clearly goes beyond \cite{bw1,gk}, viz.\
precisely the Morita invariance of the invariant $BW$ which we pointed out above.
\erem


\section{Discussion and Outlook} \label{s-last}
In various places we have already mentioned closely related works by other authors. We
summarize these references and comment on several other recent works. The relation between
classical Frobenius algebras and Frobenius algebras in $\7F$-Vect is due to Quinn \cite{q}
and Abrams \cite{abrams0,abrams}. The literature on Frobenius algebras in categories other
than Vect is quite small but has begun to grow recently. As mentioned earlier, canonical 
Frobenius algebras in $C^*$-categories (`Q-systems') were first considered in \cite{lro}, 
motivated by subfactor theory. In an algebraic-topological context commutative Frobenius
algebras in symmetric tensor categories appear in \cite{strick}. The relation between
Frobenius algebras and two-sided duals (only in ${\cal CAT}$, though) is hinted at in
\cite{khov} but not developed very far. The discussion in \cite[Section 3.3]{szl1} has
some relations to our work, and \cite[Section 3]{szl2} has some overlap with our
Subsection \ref{ss-Hopf2}. Note, however, that in these references most proofs are
omitted, and the discussion in \cite{szl2} is limited to $C^*$-bicategories.
In \cite{ko}, module categories of `rigid $\2C$-algebras' 
in braided tensor categories are considered with the aim of categorifying the
considerations on modular invariants in \cite{bek}. In view of Proposition \ref{alt}, rigid
$\2C$-algebras are nothing but Frobenius algebras. (Since the Frobenius algebras are
assumed to be commutative (i.e.\ $w'\circ c(Q,Q)=w'$) only `type I' modular invariants are
covered by this analysis.) Similar matters are considered in somewhat greater generality in
\cite{fs}, where Frobenius algebras in the sense of Definition \ref{d-Frob} appear
explicitly, influenced by a talk of the author. A recent construction by Yamagami
\cite{yamag3} bears some relation to our construction of the bicategory $\2E$. Given a
tensor category $\2A$ and a full subcategory $\2A_0\simeq H-\mod$, he constructs a
bicategory which seems to be equivalent to our $\2E$ in the special case of Subsection
\ref{ss-Hopf}, where the Frobenius algebra arises from the regular representation of a
semisimple and cosemisimple Hopf algebra $H$. 

When the present work was essentially finished we learned that the definition of the
bicategory $\2E$ via (bi)modules was discovered independently by Yamagami. Furthermore, in
a very interesting sequel \cite{ostr} of \cite{ko}, Ostrik considers modules $\2M$ over a
tensor category $\2C$. (Note that, given a bicategory $\2E$ with $\6A,\6B\in\obj\2E$, the
categories $\HOM_\2E(\6A,\6B), \HOM_\2E(\6B,\6A)$ are left and right, respectively,
modules over $\2A=\END_\2E(\6A)$.) He shows that every module over $\2C$ arises from an
algebra $A$ in $\2C$ and constructs a dual algebra $\2C^*$. The latter is equivalent to
our $\2B=\END_\2E(\6A)$. He also considers applications to modular invariants and succeeds
in phrasing most results of \cite{bek} in categorical terms, albeit without many proofs.

Results like Theorem \ref{appl} lead us to believe that all existing (and future)
applications of subfactor theory to low dimensional topology `factor through category
theory' -- as is by now well known for the knot invariants of Jones and HOMFLY. More
generally, we are convinced that essentially all algebraic aspects and results of
subfactor theory (at finite index) permit generalization to a considerably wider
categorical setting. This is further vindicated by the subsequent parts of this series
whose main results we briefly outline.

As already mentioned, in part II \cite{mue10} we prove that the center $\2Z(\2C)$ of a
finite semisimple spherical tensor category $\2C$ of non-zero dimension is weakly 
monoidally Morita equivalent to $\2C\boxtimes\2C^\op$. Furthermore, it is a modular
category in the sense of \cite{t}. In view of the relation \cite{ka} between the
categorical and the Hopf algebraic quantum double this should be interpreted as a
generalization of the fact \cite{eg1} that quantum doubles of nice Hopf algebras have 
modular representation categories. 

In Part III \cite{mue15} we will consider the bicategory $\tilde{\2E}$ mentioned in
Remark \ref{r-alt}. It will be shown to satisfy the assumptions of Theorem \ref{t-uniq},
which implies its equivalence with $\2E$. In particular, we have $\2A\approx\2B$ iff 
there exists a Frobenius algebra $\5Q$ in $\2A$ such that 
$\2B\stackrel{\otimes}{\simeq}\5Q$-mod-$\5Q$. As mentioned in Remark \ref{r-alt}, the
latter implies the equivalance of the braided tensor categories
$\2Z_1(\2A)\simeq\2Z_1(\2B)$. It is natural to ask whether the converse is true.

The programme of identifying the connections between subfactor theory (at finite index)
and category theory is certainly vindicated by the applications to the classification of
modular invariants, cf.\ \cite{ko,fs,ostr}, and to topology, as considered in Section
\ref{s-topology} and \cite{mue10}. The rapprochement of subfactor theory and `mainstream'
mathematics which this foreshadows will undoubtedly be helpful also in the classification
programme of subfactors.

\vspace{1cm}
\noindent{\it Acknowledgments.} During the long time of preparation of this work I was
financially supported by the European Union, the NSF and the NWO and hosted by the 
universities ``Tor Vergata'' and ``La Sapienza'', Rome, the IRMA, Strasbourg, the School
of Mathematical Sciences, Tel Aviv, the MSRI, Berkeley, and the Korteweg-de Vries
Institute, Amsterdam, to all of which I wish to express my sincere gratitude. 

Parts of the results of this paper and of \cite{mue10} were presented at an early stage at
the conference {\it Category Theory 99} at Coimbra, July 1999, at the conference {\it
$C^*$-algebras and tensor categories} at Cortona, August 1999, and at the workshop {\it
Quantum groups and knot theory} at Strasbourg, September 1999. More recently, I reported
on them at the conferences on `Modular invariance, ADE, subfactors and geometry of
moduli spaces' in Kyoto, November 2000 and on `Operator algebras and mathematical physics'
in Constanta, July 2001. On these and other occasions I received a lot of response and
encouragement. The following is a long but incomplete list of people whom I want to thank
for their interest and/or useful conversations:
J.\ Baez, J.\ Bernstein, A.\ Brugui\`{e}res, D.\ E.\ Evans, J.\ Fuchs, F.\ Goodman, 
M.\ Izumi, V. F. R. Jones, C.\ Kassel, Y.\ Kawahigashi, G. Kuperberg, N. P. Landsman,
R.\ Longo, G.\ Maltsiniotis, G.\ Masbaum, J.\ E.\ Roberts, N.\ Sato, V.\ Turaev, L. Tuset,
P.\ Vogel, A.\ Wassermann, H.\ Wenzl and S. Yamagami.


\begin{thebibliography}{99}
\bibitem{abrams0} L. Abrams: Two-dimensional topological quantum field theories and
   Frobenius algebras. J. Knot Th. Ramif. {\bf 5}, 569-587 (1996).
\bibitem{abrams} L. Abrams: Modules, comodules and cotensor products over Frobenius
   algebras. J. Alg. {\bf 219}, 201-213 (1999).
\bibitem{bw1} J. W. Barrett \& B. W. Westbury: Invariants of piecewise-linear 3-manifolds.
   Trans. Amer. Math. Soc. {\bf 348}, 3997-4022 (1996).
\bibitem{bw2} J. W. Barrett \& B. W. Westbury: Spherical categories. 
  Adv. Math. {\bf 143}, 357-375 (1999).
\bibitem{bw3} J. W. Barrett \& B. W. Westbury: The equality of 3-manifold invariants.
   Math. Proc. Cambr. Phil. Soc. {\bf 118}, 503-510 (1995).
\bibitem{bass} H. Bass: {\it Algebraic K-Theory}. Benjamin Inc., 1968.
\bibitem{benab} J. Benabou: Introduction to bicategories. {\it In}: Reports of the Midwest
  Category Seminar. LNM 47, pp.\ 1--77. Springer Verlag, 1967.
\bibitem{bek} J. B\"ockenhauer, D. E. Evans \& Y. Kawahigashi: On $\alpha$-induction,
   chiral generators and modular invariants for subfactors. Commun. Math. Phys. {\bf 208}, 429-487 (1999).
   Chiral structure of modular invariants for subfactors. Commun. Math. Phys. {\bf 210}, 733-784 (2000).
\bibitem{bor} F. Bor\c{c}eux: {\it Handbook of categorical algebra I. Basic category
   theory.} Cambridge University Press, 1994.
\bibitem{cp} V. Chari \& A. Pressley: {\it A guide to quantum groups}. Cambridge
   University Press, 1995. 
\bibitem{connes} A. Connes: {\it Noncommutative geometry}. Academic Press, 1994.
\bibitem{dr6} S. Doplicher \& J. E. Roberts: A new duality theory for compact groups. 
   Invent. Math. {\bf 98}, 157-218 (1989).
\bibitem{eg1} P. Etingof \& S. Gelaki: Some properties of finite dimensional semi-simple
   Hopf algebras. Math. Res. Lett. {\bf 5}, 191-197 (1998).
\bibitem{eg2} P. Etingof \& S. Gelaki: On finite-dimensional semisimple and cosemisimple
   Hopf algebras in positive characterstic. Int. Math. Res. Not. {\bf 1998}, 187-195.
\bibitem{ek3} D. E. Evans \& Y. Kawahigashi: From subfactors to 3-dimensional topological
   quantum field theories and back. A detailed account of Ocneanu's theory. 
   Intern. J. Math. {\bf 6}, 537-558 (1995).
\bibitem{ek2} D. E. Evans \& Y. Kawahigashi: {\it Quantum symmetries on operator 
   algebras}. Oxford University Press, 1998.
\bibitem{fgsv} J. Fuchs, A. Ch. Ganchev, K. Szlach\'{a}nyi \& P. Vecserny\'{e}s:
   $S_4$ symmetry of $6j$ symbols and Frobenius-Schur indicators in rigid monoidal $C^*$
   categories. Journ. Math. Phys. {\bf 40}, 408-426 (1999).
\bibitem{fs} J. Fuchs \& C. Schweigert: Category theory for conformal boundary conditions.
   {\tt math.CT/0106050}.
\bibitem{gr} P. Gabriel \& A. V. Roiter: {\it Representations of finite-dimensional
   algebras}. Springer Verlag, 1992. 
\bibitem{gk} S. Gelfand \& D. Kazhdan: Invariants of three-dimensional manifolds.
   Geom. Funct. Anal. {\bf 6}, 268-300 (1996).
\bibitem{glr} P. Ghez, R. Lima \& J. E. Roberts: $W^*$-categories. Pac. J. Math.
   {\bf 120}, 79-109 (1985).
\bibitem{gray} J. W. Gray: {\it Formal category theory: Adjointness for 2-categories}. LNM 
   391. Springer Verlag, 1974. 
\bibitem{iz2} M. Izumi: The structure of sectors associated with Longo-Rehren 
   inclusions I. General theory, Commun. Math. Phys. {\bf 213}, 127-179 (2000).
\bibitem{vfr1} V. F. R. Jones: Index for subfactors. Invent. Math. {\bf 72}, 1-25 (1983).
\bibitem{vfr3} V. F. R. Jones: {\it Subfactors and knots}. CBMS Regional Conference Series
   in Mathematics, 80. American Mathematical Society, Providence, RI, 1991.
\bibitem{vfr-su} V. F. R. Jones \& V. S. Sunder: {\it Introduction to subfactors}. Cambridge
  University Press, 1997. 
\bibitem{js2} A. Joyal \& R. Street: Tortile Yang-Baxter operators in tensor categories.
   J. Pure Appl. Alg. {\bf 71}, 43-51 (1991).
\bibitem{js1} A. Joyal \& R. Street: Braided tensor categories. Adv. Math. {\bf 102}, 
   20-78 (1993).
\bibitem{kad} L. Kadison: {\it New examples of Frobenius extensions}. University Lecture
   Series \#14, AMS, 1999.
\bibitem{ka} C. Kassel: {\it Quantum groups}. Springer Verlag, 1995.
\bibitem{ks} G. M. Kelly \& R. Street: Review of the elements of 2-categories. LNM 420,
   pp.\ 75-103. Springer Verlag, 1974.
\bibitem{khov} M. Khovanov: A functor-valued invariant of tangles. {\tt math.QA/0103190}.
\bibitem{ko} A. Kirillov Jr. \& V. Ostrik: On $q$-analog of McKay correspondence and ADE
   classification of $\widehat{{\mathfrak sl}_2}$ conformal field theories. 
   {\tt math.QA/0101219}.
\bibitem{kosu} V. Kodiyalam \& V. S. Sunder: {\it Topological quantum field theories from
   subfactors}. Chapman \& Hall/CRC Research Notes in Mathematics, 2001.
\bibitem{kos} H. Kosaki: Extension of Jones' theory on index to arbitrary factors.
   Journ. Funct. Anal. {\bf 66}, 123-140 (1986).
\bibitem{ku} G. Kuperberg: Involutory Hopf algebras and 3-manifold
   invariants. Intern. J. Math. {\bf 2}, 41-66 (1991).
\bibitem{lar2} R. G. Larson \& D. E. Radford: Finite dimensional cosemisimple Hopf
   algebras in characteristic zero are semisimple. J. Alg. {\bf 117}, 267-289 (1988).
\bibitem{ls} R. G. Larson \& M. E. Sweedler: An associative orthogonal bilinear form 
   for Hopf algebras. Amer. J. Math. {\bf 91}, 75-94 (1969).
\bibitem{lo1} R. Longo: Index of subfactors and statistics of quantum fields I \& II.
   Commun. Math. Phys. {\bf 126}, 217-247 (1989) \& {\bf 130}, 285-309 (1990). 
\bibitem{lo3} R. Longo: A duality for Hopf algebras and for subfactors I. Commun. Math. Phys. 
   {\bf 159}, 133-150 (1994).
\bibitem{lre} R. Longo \& K.-H. Rehren: Nets of subfactors. 
   Rev. Math. Phys. {\bf 7}, 567-597 (1995).
\bibitem{lro} R. Longo \& J. E. Roberts: A theory of dimension. K-Theory {\bf 11}, 
   103-159 (1997).
\bibitem{mack} M. Mackaay: Spherical 2-categories and 4-manifold invariants. 
   Adv. Math. {\bf 143}, 288-348 (1999).
\bibitem{cwm} S. Mac Lane: {\it Categories for the working mathematician}. 2nd ed. 
   Springer Verlag, 1998.
\bibitem{ma1} S. Majid: Representations, duals and quantum doubles of monoidal 
  categories. Rend. Circ. Mat. Palermo Suppl. {\bf 26}, 197-206 (1991).
\bibitem{mue6} M. M\"uger: Galois theory for braided tensor categories and the modular
   closure. Adv. Math. {\bf 150}, 151-201 (2000).
\bibitem{mue10} M. M\"uger: From subfactors to categories and topology II. The quantum
  double of tensor categories and subfactors. {\tt math.CT/0111205}.
\bibitem{mue15} M. M\"uger: From subfactors to categories and topology III. 
   In preparation. 
\bibitem{mue-t} M. M\"uger \& L. Tuset: Representations of algebraic quantum groups and
   reconstruction theorems for tensor categories II. Regular representations and embedding
   theorems. In preparation. 
\bibitem{neuchl} M. Neuchl: Representation theory of Hopf categories. Dissertation. Munich
   university. Available from {\tt http://www.mathematik.uni-muenchen.de/$\sim$neuchl/}.
\bibitem{nill} F. Nill: Weyl algebras, Fourier transformations and integrals on
   finite-dimensional Hopf algebras. Rev. Math. Phys. {\bf 6}, 149-166 (1994).
\bibitem{ocn1} A. Ocneanu: Quantized group string algebras and Galois theory for
   algebras. {\it In:} D. E. Evans \& M. Takesaki (eds.): {\it Operator algebras and
   applications}, Vol.\ 2. London Math. Soc. Lect. Notes 136. 
   Cambridge University Press, 1988.
\bibitem{ocn3} A. Ocneanu: An invariant coupling between 3-manifolds and subfactors, with
   connections to topological and conformal quantum field theory. Unpublished
  manuscript. Ca.\ 1991.
\bibitem{ostr} V. Ostrik: Module categories, weak Hopf algebras and modular invariants. 
   {\tt math.QA/0111139}.
\bibitem{par0} B. Pareigis: Non-additive ring and module theory II. C-categories,
   C-functors, and C-morphisms. Publ. Math. Debrecen {\bf 24}, 351-361 (1977).
\bibitem{par} B. Pareigis: Non-additive ring and module theory III. Morita equivalences.
   Publ. Math. Debrecen {\bf 25}, 177-186 (1978); Morita equivalence of module categories
   with tensor products.  Commun. Alg. {\bf 9}, 1455-1477 (1981).
\bibitem{popa} S. Popa: Classification of subfactors and their endomorphisms. CBMS
   Regional Conference Series in Mathematics, 86. AMS, 1995.
\bibitem{q} F. Quinn: Lectures on axiomatic topological quantum field theory. {\it In:}
  Geometry and quantum field theory (Park City, UT, 1991), 323--453, IAS/Park City
   Math. Ser., 1, Amer. Math. Soc.,  Providence, RI, 1995. 
\bibitem{sch} P. Schauenburg: The monoidal center construction and bimodules. 
   J. Pure Appl. Alg. {\bf 158}, 325-346 (2001).
\bibitem{stzs} \c{S}. Str\v{a}til\v{a} \& L. Zsid\'{o}: {\it Lectures on von Neumann
   algebras}. Abacus Press, 1979. 
\bibitem{str1} R. Street: The formal theory of monads. 
   J. Pure Appl. Alg. {\bf 2}, 149-168 (1972).
\bibitem{str2} R. Street: The quantum double and related constructions. J. Pure Appl. Alg.
   {\bf 132}, 195-206 (1998).
\bibitem{strick} N. P. Strickland: $K(N)$-local duality for finite groups and groupoids.
   Topology {\bf 39}, 733-772 (2000). 
\bibitem{sun} V. S. Sunder: {\it An invitation to von Neumann algebras}. 
   Springer Verlag, 1987.
\bibitem{szl1} K. Szlach\'{a}nyi: Finite quantum groupoids and inclusions of finite
   type. To appear in: R. Longo (ed.): {\it Mathematical Physics in Mathematics and
   Physics: Quantum and Operator Algebraic Aspects}. Fields Institute Communications,
   2001. {\tt math.QA/0011036}.
\bibitem{szl2} K. Szlach\'{a}nyi: Weak Hopf algebra symmetries of $C^*$-algebra
   inclusions. {\tt math.QA/0101005}.
\bibitem{szym} W. Szyma\'{n}ski: Finite index subfactors and Hopf algebra crossed
   products. Proc. AMS {\bf 120}, 519-528 (1994).
\bibitem{takes} M. Takesaki: {\it Theory of operator algebras I}. Springer Verlag, 1979.
\bibitem{t} V. G. Turaev: {\it Quantum invariants of knots and 3-manifolds}. 
   Walter de Gruyter, 1994.
\bibitem{tv} V. G. Turaev \& O. Y. Viro: State sum invariants of 3-manifolds and quantum
   6j-symbols, Topology {\bf 31}, 865-902 (1992).
\bibitem{w} H. Wenzl: $C^*$-Tensor categories from quantum groups. J. Amer. Math. Soc. 
   {\bf 11}, 261-282 (1998).
\bibitem{xu1} F. Xu: New braided endomorphisms from conformal inclusions. 
   Commun. Math. Phys. {\bf 192}, 349-403 (1998).
\bibitem{yamag1} S. Yamagami: On unitary representation theories of compact quantum
   groups. Commun. Math. Phys. {\bf 167}, 509-529 (1995).
\bibitem{yamag3} S. Yamagami: Tannaka duals in semisimple tensor categories. 
   {\tt math.CT/0106065}.
\bibitem{yamag2} S. Yamagami: Frobenius reciprocity in tensor categories; Frobenius
   duality in $C^*$-tensor categories. Preprints, available at \\
   {\tt http://suuri.sci.ibaraki.ac.ap/$\sim$yamagami/}.
\bibitem{yamag4} S. Yamagami: Frobenius algebras in tensor categories. Manuscript.
\bibitem{yaman} T. Yamanouchi: Construction of an outer action of a finite-dimensional Kac
   algebra on the AFD factor of type II$_1$. Int. J. Math. {\bf 4}, 1007-1045 (1993).
\bibitem{y} D. N. Yetter: Framed tangles and a theorem of Deligne on braided deformations
   of tannakian categories. Contemp. Math. {\bf 134}, 325-350 (1992).
\end{thebibliography}
\end{document}